\documentclass[a4paper,12pt]{article}

\usepackage{amsmath}
\usepackage{amssymb}
\usepackage{amsthm}
\usepackage{amscd}

\newtheorem{theorem}{Theorem}[section]
\newtheorem{lemma}[theorem]{Lemma}
\newtheorem{coro}[theorem]{Corollary}
\newtheorem{proposition}[theorem]{Proposition}
\theoremstyle{definition}
\newtheorem{example}[theorem]{Example}

\newtheorem{rules}[theorem]{Rule}

\def\Aff{\mbox{\sl Aff}}
\def\bt{\tilde{b}}

\def\geh{\mathfrak{g}}
\def\La{\Lambda}
\def\la{\lambda}
\def\ol#1{\overline{#1}}

\def\cd{\cdots}
\def\ot{\otimes}

\def\veps{\varepsilon}
\def\vphi{\varphi}

\def\T{{\mathbb T}}

\def\Z{{\mathbb Z}}
\def\Zn{\Z_{\ge0}}
\def\etd{\tilde{e}}
\def\ftd{\tilde{f}}
\def\et#1{\tilde{e}_{#1}}
\def\ft#1{\tilde{f}_{#1}}
\newcommand{\DOMINOnormal}[2]{%
	\put(0,0){\line(1,0){1}}
	\put(0,1){\line(1,0){1}}
	\put(0,2){\line(1,0){1}}
	\put(0,0){\line(0,1){2}}
	\put(1,0){\line(0,1){2}}
	\put(0,0){\makebox(1,1){$#2$}}
	\put(0,1){\makebox(1,1){$#1$}}
}

\newcommand{\DOMINOlong}[2]{%
	\put(0,0){\line(1,0){2}}
	\put(0,1){\line(1,0){2}}
	\put(0,2){\line(1,0){2}}
	\put(0,0){\line(0,1){2}}
	\put(2,0){\line(0,1){2}}
	\put(0,0){\makebox(2,1){$#2$}}
	\put(0,1){\makebox(2,1){$#1$}}
}

\newcommand{\HAKOnormal}[1]{%
	\put(0,0){\line(1,0){1}}
	\put(0,1){\line(1,0){1}}
	\put(0,0){\line(0,1){1}}
	\put(1,0){\line(0,1){1}}
	\put(0,0){\makebox(1,1){$#1$}
	}
}

\newcommand{\HAKOlong}[1]{%
	\put(0,0){\line(1,0){2}}
	\put(0,1){\line(1,0){2}}
	\put(0,0){\line(0,1){1}}
	\put(2,0){\line(0,1){1}}
	\put(0,0){\makebox(2,1){$#1$}}
}

\newcommand{\BeforeCOLUMNINSERTIONnn}[3]{%
	\put(0,0){\HAKOnormal{#3}}
	\put(1,0){\makebox(1,1){$\to$}}
	\put(2,0){\DOMINOnormal{#1}{#2}
	}
}

\newcommand{\BeforeCOLUMNINSERTIONnl}[3]{%
	\put(0,0){\HAKOnormal{#3}}
	\put(1,0){\makebox(1,1){$\to$}}
	\put(2,0){\DOMINOlong{#1}{#2}
	}
}

\newcommand{\BeforeCOLUMNINSERTIONll}[3]{%
	\put(0,0){\HAKOlong{#3}}
	\put(2,0){\makebox(1,1){$\to$}}
	\put(3,0){\DOMINOlong{#1}{#2}
	}
}

\newcommand{\AfterCOLUMNINSERTIONnn}[3]{%
	\put(0,0){\DOMINOnormal{#2}{#3}}
	\put(2.5,0){\makebox(1,1){$\to$}}
	\put(1.5,0){\HAKOnormal{#1}}
}

\newcommand{\AfterCOLUMNINSERTIONln}[3]{%
	\put(0,0){\DOMINOlong{#2}{#3}}
	\put(3.5,0){\makebox(1,1){$\to$}}
	\put(2.5,0){\HAKOnormal{#1}}
}

\newcommand{\AfterCOLUMNINSERTIONll}[3]{%
	\put(0,0){\DOMINOlong{#2}{#3}}
	\put(4.5,0){\makebox(1,1){$\to$}}
	\put(2.5,0){\HAKOlong{#1}}
}

\newcommand{\ToBOX}[2]{%
	\setlength{\unitlength}{5.5mm}
	\begin{picture}(9,1.5)(-2,0)
	\put(-1,0){\makebox(0,1)[r]{\bfseries{#1}}}
	\put(0,0){
		\put(0,0){\HAKOnormal{#2}}
		\put(1,0){\makebox(1,1){$\to$}}
		\put(2,0){\makebox(1,1){$\emptyset$}}
		}
	\put(3.5,0){\makebox(1,1){$=$}}
	\put(5,0){\HAKOnormal{#2}}
	\end{picture}
}

\newcommand{\ToDOMINO}[5]{%
	\setlength{\unitlength}{5.5mm}
	\begin{picture}(9,2.5)(-2,0)
	\put(-1,0){\makebox(0,2)[r]{\bfseries{#1}}}
	\put(0,0.5){
		\put(0,0){\HAKOnormal{#3}}
		\put(1,0){\makebox(1,1){$\to$}}
		\put(2,0){\HAKOnormal{#2}}
		}
	\put(3.5,0.5){\makebox(1,1){$=$}}
	\put(5,0){\DOMINOnormal{#4}{#5}}
	\end{picture}
}

\newcommand{\ToYOKODOMINO}[5]{%
	\setlength{\unitlength}{5.5mm}
	\begin{picture}(11,1.5)(-2,0)
	\put(-1,0){\makebox(0,1)[r]{\bfseries{#1}}}
	\put(0,0){
		\put(0,0){\HAKOnormal{#3}}
		\put(1,0){\makebox(1,1){$\to$}}
		\put(2,0){\HAKOnormal{#2}}
		}
	\put(3.5,0){\makebox(1,1){$=$}}
	\put(5,0){
		\put(0,0){\HAKOnormal{#5}}
		\put(2.5,0){\makebox(1,1){$\to$}}
		\put(1.5,0){\HAKOnormal{#4}}
		}
	\end{picture}
}

\newcommand{\FromYOKODOMINO}[5]{%
	\setlength{\unitlength}{5.5mm}
	\begin{picture}(11,1.5)(-2,0)
	\put(-1,0){\makebox(0,1)[r]{\bfseries{#1}}}
	\put(0,0){
		\put(0,0){\HAKOnormal{#3}}
		\put(2.5,0){\makebox(1,1){$\to$}}
		\put(1.5,0){\HAKOnormal{#2}}
		}
	\put(3.5,0){\makebox(1,1){$=$}}
	\put(5,0){
		\put(0,0){\HAKOnormal{#5}}
		\put(1,0){\makebox(1,1){$\to$}}
		\put(2,0){\HAKOnormal{#4}}
		}
	\end{picture}
}

\newcommand{\ToHOOKnn}[7]{%
	\setlength{\unitlength}{5.5mm}
	\begin{picture}(11,2.5)(-2,0)
	\put(-1,0){\makebox(0,2)[r]{\bfseries{#1}}}
	\put(0,0){\BeforeCOLUMNINSERTIONnn{#2}{#3}{#4}}
	\put(3.5,0){\makebox(1,2){$=$}}
	\put(5,0){\AfterCOLUMNINSERTIONnn{#5}{#6}{#7}}
	\end{picture}
}

\newcommand{\ToHOOKln}[7]{%
	\setlength{\unitlength}{5.5mm}
	\begin{picture}(13,2.5)(-2,0)
	\put(-1,0){\makebox(0,2)[r]{\bfseries{#1}}}
	\put(0,0){\BeforeCOLUMNINSERTIONnn{#2}{#3}{#4}}
	\put(3.5,0){\makebox(1,2){$=$}}
	\put(5,0){\AfterCOLUMNINSERTIONln{#5}{#6}{#7}}
	\end{picture}
}

\newcommand{\ToHOOKll}[7]{%
	\setlength{\unitlength}{5.5mm}
	\begin{picture}(13,2.5)(-2,0)
	\put(-1,0){\makebox(0,2)[r]{\bfseries{#1}}}
	\put(0,0){\BeforeCOLUMNINSERTIONnn{#2}{#3}{#4}}
	\put(3.5,0){\makebox(1,2){$=$}}
	\put(5,0){\AfterCOLUMNINSERTIONll{#5}{#6}{#7}}
	\end{picture}
}

\newcommand{\ToHOOKllnn}[7]{%
	\setlength{\unitlength}{5.5mm}
	\begin{picture}(13,2.5)(-2,0)
	\put(-1,0){\makebox(0,2)[r]{\bfseries{#1}}}
	\put(0,0){\BeforeCOLUMNINSERTIONll{#2}{#3}{#4}}
	\put(5.5,0){\makebox(1,2){$=$}}
	\put(7,0){\AfterCOLUMNINSERTIONnn{#5}{#6}{#7}}
	\end{picture}
}

\newcommand{\FromHOOKnn}[7]{%
	\setlength{\unitlength}{5.5mm}
	\begin{picture}(11,2.5)(-2,0)
	\put(-1,0){\makebox(0,2)[r]{\bfseries{#1}}}
	\put(0,0){\AfterCOLUMNINSERTIONnn{#2}{#3}{#4}}
	\put(3.5,0){\makebox(1,2){$=$}}
	\put(5,0){\BeforeCOLUMNINSERTIONnn{#5}{#6}{#7}}
	\end{picture}
}

\newcommand{\FromHOOKnl}[7]{%
	\setlength{\unitlength}{5.5mm}
	\begin{picture}(13,2.5)(-2,0)
	\put(-1,0){\makebox(0,2)[r]{\bfseries{#1}}}
	\put(0,0){\AfterCOLUMNINSERTIONnn{#2}{#3}{#4}}
	\put(3.5,0){\makebox(1,2){$=$}}
	\put(5,0){\BeforeCOLUMNINSERTIONnl{#5}{#6}{#7}}
	\end{picture}
}

\newcommand{\FromHOOKll}[7]{%
	\setlength{\unitlength}{5.5mm}
	\begin{picture}(13,2.5)(-2,0)
	\put(-1,0){\makebox(0,2)[r]{\bfseries{#1}}}
	\put(0,0){\AfterCOLUMNINSERTIONnn{#2}{#3}{#4}}
	\put(3.5,0){\makebox(1,2){$=$}}
	\put(5,0){\BeforeCOLUMNINSERTIONll{#5}{#6}{#7}}
	\end{picture}
}

\newcommand{\FromHOOKllnn}[7]{%
	\setlength{\unitlength}{5.5mm}
	\begin{picture}(13,2.5)(-2,0)
	\put(-1,0){\makebox(0,2)[r]{\bfseries{#1}}}
	\put(0,0){\AfterCOLUMNINSERTIONll{#2}{#3}{#4}}
	\put(5.5,0){\makebox(1,2){$=$}}
	\put(7,0){\BeforeCOLUMNINSERTIONnn{#5}{#6}{#7}}
	\end{picture}
}

\newcommand{\fbx}[1]{\fbox{\rule{0mm}{4mm}\parbox[t]{2.4mm}{$ #1$}}}
\newcommand{\bbx}[2]{\overbrace{\fbox{$\vphantom{\ol{1}} #1 \cd #1$}}^{#2}}
\newcommand{\mapright}[1]{%
  \smash{\mathop{%
   \hbox to 1cm{\rightarrowfill}}\limits^{#1}}}

\makeatletter
\@addtoreset{equation}{section}
\makeatother

\title{Combinatorial $R$ matrices for a family \\
of crystals : $B^{(1)}_n$, $D^{(1)}_{n}$, $A^{(2)}_{2n}$ 
and $D^{(2)}_{n+1}$ cases}
\author{
Goro Hatayama\thanks{
Institute of Physics, University of Tokyo, Komaba, Tokyo 153-8902, Japan},
Atsuo Kuniba,$\hspace{-1.2mm}^*$\\
Masato Okado\thanks{Department of Informatics and Mathematical Science,
Graduate School of Engineering Science,
Osaka University,
Toyonaka, Osaka 560-8531,
Japan}
and Taichiro Takagi\thanks{
Department of Applied Physics, National Defense Academy,
Kanagawa 239-8686, Japan}}
\date{}
\begin{document}
\maketitle
\begin{abstract}
For coherent families of crystals of affine Lie algebras
of type $B^{(1)}_n$, $D^{(1)}_{n}$, $A^{(2)}_{2n}$ and $D^{(2)}_{n+1}$ 
we describe the combinatorial $R$ matrix using
column insertion algorithms for $B,C,D$ Young tableaux.
This is a continuation of \cite{HKOT}.
\end{abstract}
\section{Introduction}
\label{sec:intro}
A combinatorial $R$ matrix 
is the $q = 0$ limit of the quantum $R$ matrix for a quantum affine algebra
$U_q(\geh)$,
where $q$ is the deformation parameter and $q=1$ means non-deformed.
It is defined on the tensor
product of two affine crystals $\Aff(B)\ot\Aff(B')$ 
(See Section \ref{sec:crystals} for notations), and consists of
an isomorphism and an energy function.
It was first introduced in \cite{KMN1} for the {\em homogeneous} case
where one has $B=B'$.
In this case the isomorphism is trivial.
The energy function was used to 
describe the path realization of the crystals of highest weight
representations of quantum affine algebras.
The definition of the energy function 
was extended in \cite{NY} to the {\em inhomogeneous} case,
i.e. $B\neq B'$, to study the charge of the 
Kostka-Foulkes polynomials \cite{Ma,LS,KR}.

In \cite{KKM} the theory of
coherent families of perfect crystals
was developed for quantum affine algebras
of type $A^{(1)}_n,A^{(2)}_{2n-1}, A^{(2)}_{2n}$,
$B^{(1)}_n, C^{(1)}_n$, $D^{(1)}_n$ and $D^{(2)}_{n+1}$.
An element of these crystals is written as
an array of nonnegative integers
and an explicit description of the
energy functions is given
in terms of piecewise linear functions
of its entries for the homogeneous case.
Unfortunately this description is not applicable to 
the inhomogeneous cases.
The purpose of this paper is to give an explicit
description of the isomorphism and energy function
for the inhomogeneous cases.

The main tool of our description 
is an insertion algorithm, that is a certain procedure on Young tableaux.
Insertion algorithm itself had been invented
in the context of the Robinson-Schensted correspondence \cite{F}
long before the crystal basis theory was initiated, and 
subsequently generalized in, e.g. \cite{Ber,P,Su}.
As far as $A_n^{(1)}$ crystals are concerned, 
the isomorphisms and 
energy functions were obtained 
in terms of usual (type $A$)
Young tableaux and insertion algorithms thereof \cite{S,SW}.
In contrast,
no similar description for the combinatorial $R$ matrix had been made
for other quantum affine algebras, since
an insertion algorithm suitable for the $B,C,D$ tableaux
given in \cite{KN} was known only recently \cite{B1,B2,L}.
In the previous work \cite{HKOT} the authors gave a description for 
type $C_n^{(1)}$ and $A_{2n-1}^{(2)}$.
There we used the $\mathfrak{sp}$-version of semistandard tableaux
defined in \cite{KN} 
and the column insertion algorithm presented in \cite{B1} on these tableaux.
In this paper we study the remaining types,
$A^{(2)}_{2n}, D^{(2)}_{n+1}, B^{(1)}_n$ and $D^{(1)}_n$.
We use $\mathfrak{sp}$- and $\mathfrak{so}$- versions of semistandard tableaux
defined in \cite{KN} 
and the column insertion algorithms presented in \cite{B1,B2}
on these tableaux.

The layout of this paper is as follows.
In Section \ref{sec:crystals} we give a brief review of the basic notions
in the theory of crystals and give the definition of combinatorial $R$ matrix.
We first give the description for types
$A^{(2)}_{2n}$ and $D^{(2)}_{n+1}$ in Section \ref{sec:Cx}.
In Sections \ref{subsec:twista} and \ref{subsec:twistd}
we recall the definitions of crystal $B_l$ for
type $A^{(2)}_{2n}$ and $D^{(2)}_{n+1}$ respectively,
and give a description of its elements in terms of one-row tableaux.
We introduce the map $\omega$ from these crystals to the crystal
of type $C^{(1)}_n$, hence the procedure is reduced to that of
the latter case which we have already developed in \cite{HKOT}.
In Section \ref{subsec:ccis} we list up elementary operations of 
column insertions and their inverses
for type $C$ tableaux with at most two rows.
In Section \ref{subsec:ruleCx} we give the main theorem
and give the description of the isomorphism and energy function
for type $A^{(2)}_{2n}$ and $D^{(2)}_{n+1}$, and
in Section \ref{subsec:exCx} we give examples.
The $B^{(1)}_n$ and $D^{(1)}_n$ cases are treated in Section \ref{sec:bd}.
The layout is parallel to Section \ref{sec:Cx}. 
In Sections \ref{subsec:cib} and \ref{subsec:cid}, however, we also
prove the column bumping lemmas 
(Lemmas \ref{lem:bcblxx} and \ref{lem:cblxx}) for $B$ and $D$
tableaux, since a route in the tableau made from inserted letters 
(bumping route) has some importance in the main theorem.
\section{\mathversion{bold}Crystals and combinatorial $R$ matrix}
\label{sec:crystals}
Let us recall basic notions in the theory of crystals.

See \cite{KMN1,KKM} for details.
Let $I=\{0,1,\cdots,n\}$ be the index set.
Let $B$ be a $P_{\scriptstyle \mbox{\scriptsize \sl cl}}$-weighted 
crystal, i.e. $B$ is a finite set equipped with the {\em crystal structure} that is given by the
maps $\etd_i$ and $\ftd_i$ from $B\sqcup \{0\}$ to $B \sqcup \{0\}$
and maps $\veps_i$ and $\vphi_i$ from $B$ to $\Zn$.
It is always assumed that $\etd_i 0 = \ftd_i 0 = 0$ and $\ftd_i b = b'$ means
$\etd_i b' = b$.

The crystal $B$ is identified with a colored oriented graph ({\em crystal graph})
if one draws an arrow as  $b \stackrel{i}{\rightarrow} b'$ for $\ftd_i b = b'$.
Such an arrow is called $i$-arrow. Pick any $i$ and neglect all the $j$-arrows
with $j\neq i$. One then finds that all the connected 
components are {\em strings} of finite lengths,
i.e. there is no loop or branch.
Fix a string and take any node $b$ in the string.
Then the maps $\veps_i(b),\vphi_i(b)$ have the following meaning.
Along the string you can go forward by $\vphi_i(b)$ steps to an end
following the arrows and backward by $\veps_i(b)$ steps against the arrows.

Given two crystals $B$ and $B'$, let $B \otimes B'$ be a crystal 
defined as follows. As a set it is identified with $B\times B'$.
The actions of the operators $\et{i},\ft{i}$ 
on $B\ot B'$ are given by
\begin{eqnarray*}
\et{i}(b\ot b')&=&\left\{
\begin{array}{ll}
\et{i} b\ot b'&\mbox{ if }\vphi_i(b)\ge\veps_i(b')\\
b\ot \et{i} b'&\mbox{ if }\vphi_i(b) < \veps_i(b'),
\end{array}\right. \\
\ft{i}(b\ot b')&=&\left\{
\begin{array}{ll}
\ft{i} b\ot b'&\mbox{ if }\vphi_i(b) > \veps_i(b')\\
b\ot \ft{i} b'&\mbox{ if }\vphi_i(b)\le\veps_i(b').
\end{array}\right. 
\end{eqnarray*}
Here $0\ot b'$ and $b\ot 0$ should be understood as $0$.
All crystals $B$ and the tensor products of them $B\ot B'$ are connected 
as a graph.

Let $\mbox{\sl Aff} (B)=\left\{ z^d b | b \in B,\,d \in \Z \right\}$ be an
affinization of $B$ \cite{KMN1},
where $z$ is an indeterminate.
The crystal
$\Aff (B)$ is equipped with the crystal structure, where
the actions of $\et{i},\ft{i}$ are defined as
$\et{i}\cdot z^d b=z^{d+\delta_{i0}}(\et{i}b),\,
\ft{i}\cdot z^d b=z^{d-\delta_{i0}}(\ft{i}b)$. 
The {\em combinatorial $R$ matrix} 
is given by
\begin{eqnarray*}
R\;:\;\Aff(B)\ot\Aff(B')&\longrightarrow&\Aff(B')\ot\Aff(B)\\
z^d b\ot z^{d'} b'&\longmapsto&z^{d'+H(b\ot b')}\bt'\ot z^{d-H(b\ot b')}\bt,
\end{eqnarray*}
where $\iota (b\ot b') = \bt'\ot\bt$ under the isomorphism $\iota: B\ot B'
\stackrel{\sim}{\rightarrow}B'\ot B$. 
$H(b\ot b')$ is called the 
{\em energy function} and determined up to a global additive constant by
\[
H(\et{i}(b\ot b'))=\left\{%
\begin{array}{ll}
H(b\ot b')+1&\mbox{ if }i=0,\ \vphi_0(b)\geq\veps_0(b'),\ 
\vphi_0(\bt')\geq\veps_0(\bt),\\
H(b\ot b')-1&\mbox{ if }i=0,\ \vphi_0(b)<\veps_0(b'),\ 
\vphi_0(\bt')<\veps_0(\bt),\\
H(b\ot b')&\mbox{ otherwise},
\end{array}\right.
\]
since $B\ot B'$ is connected.
By definition $\iota$ satisfies $\et{i} \iota = \iota \et{i}$ and 
$\ft{i} \iota = \iota \ft{i}$ on $B \ot B'$.
The definition of the energy function
assures the intertwining property of
$R$, i.e. $\et{i} R = R \et{i}$ and $\ft{i} R = R \ft{i}$
on $\Aff(B)\ot\Aff(B')$.
In the remaining part of this paper we do not 
stick to the formalism on $\Aff(B)\ot\Aff(B')$
and rather treat the isomorphism and
energy function separately.
\section{\mathversion{bold}$U_q'(A_{2n}^{(2)})$ and $U_q'(D_{n+1}^{(2)})$ 
crystal cases}
\label{sec:Cx}
\subsection{\mathversion{bold}Definitions : $U_q'(A_{2n}^{(2)})$ case}
\label{subsec:twista}
Given a positive integer $l$,
we consider a $U_q'(A_{2n}^{(2)})$ crystal denoted by $B_l$,
that is defined
in \cite{KKM}.
$B_l$'s are the crystal bases
of the irreducible finite-dimensional representations
of the quantum affine algebra $U_q'(A_{2n}^{(2)})$.
As a set $B_{l}$ reads
$$
B_{l} = \left\{(
x_1,\ldots, x_n,\overline{x}_n,\ldots,\overline{x}_1) \Biggm|
x_i, \overline{x}_i \in \Z_{\ge 0},
\sum_{i=1}^n(x_i + \overline{x}_i) \in \{l,l-1,\ldots ,0\} \right\}.
$$
%
For its crystal structure see \cite{KKM}.
$B_{l}$ is isomorphic to
$\bigoplus_{0 \leq j \leq l} B(j \La_1)$ as
a $U_q(C_n)$ crystal, where $B(j \La_1)$ is
the crystal associated with the irreducible representation of $U_q(C_n)$
with highest weight $j \Lambda_1$.
The $U_q(C_n)$ crystal $B(j \La_1)$ has a description in terms of
semistandard $C$ tableaux \cite{KN}.
The entries are $1,\ldots ,n$ and $\ol{1}, \ldots ,\ol{n}$ with the
total order:
\begin{displaymath}
1 < 2 < \cd < n < \ol{n} < \cd < \ol{2} < \ol{1}.
\end{displaymath}
For an element $b$ of $B(j\La_1)$ let us 
denote by $\mathcal{T}(b)$ the 
tableau associated with $b$.
Thus for
$b= (x_1, \ldots, x_n, \overline{x}_n,\ldots,\overline{x}_1) \in
B(j \La_1)$
the tableau $\mathcal{T}(b)$ is depicted by
\begin{equation}
\label{eq:tabtwistax}
\mathcal{T}(b)=\overbrace{\fbox{$\vphantom{\ol{1}} 1 \cd 1$}}^{x_1}\!
\fbox{$\vphantom{\ol{1}}\cd$}\!
\overbrace{\fbox{$\vphantom{\ol{1}}n \cd n$}}^{x_n}\!
\overbrace{\fbox{$\vphantom{\ol{1}}\ol{n} \cd \ol{n}$}}^{\ol{x}_n}\!
\fbox{$\vphantom{\ol{1}}\cd$}\!
\overbrace{\fbox{$\ol{1} \cd \ol{1}$}}^{\ol{x}_1}.
\end{equation}
The length of this one-row tableau is equal to $j$, namely
$\sum_{i=1}^n(x_i + \overline{x}_i) =j$.
Here and in the remaining part of this paper we denote
$\overbrace{\fbox{$\vphantom{\ol{1}}i$} \fbox{$\vphantom{\ol{1}}i$} \fbox{$\vphantom{\ol{1}}\cd$} \fbox{$\vphantom{\ol{1}}i$}}^{x}$ by
\par\noindent
\setlength{\unitlength}{5mm}
\begin{picture}(22,3)(-6,0)
\put(0,0){\makebox(10,3)
{$\overbrace{\fbox{$\vphantom{\ol{1}} i \cd i$}}^{x}$ or
more simply by }}
\put(10,0.5){\line(1,0){3}}
\put(10,1.5){\line(1,0){3}}
\put(10,0.5){\line(0,1){1}}
\put(13,0.5){\line(0,1){1}}
\put(10,0.5){\makebox(3,1){$i \cdots i$}}
\put(10,1.5){\makebox(3,1){${\scriptstyle x}$}}
\put(13,0){\makebox(1,1){.}}
\end{picture}
\par\noindent
To describe our rule for the combinatorial $R$ matrix
we shall depict the elements of
$B_{l}$ 
by one-row tableaux with length $2l$.
We do this by duplicating each letters and then by supplying
pairs of $0$ and $\ol{0}$.
Adding $0$ and $\ol{0}$ into the set of entries of the
tableaux, we assume
the total order $0 < 1 < \cd <  \ol{1} < \ol{0}$.

Let us introduce the map $\omega$ from the $U_q'(A_{2n}^{(2)})$ 
crystal $B_l$ to the $U_q'(C_{n}^{(1)})$ crystal\footnote{Here we adopted 
the notation $B_{2l}$ that we have used in the previous work \cite{HKOT}. 
This $B_{2l}$ was originally denoted by $B_l$ in \cite{KKM}.}
$B_{2l}$.
This $\omega$ sends
$b= (x_1, \ldots, x_n, \overline{x}_n,\ldots,\overline{x}_1)$ to
$\omega (b) = (2x_1, \ldots, 2x_n ,$ $2\overline{x}_n ,\ldots,2\overline{x}_1)$.
On the other hand let us introduce the symbol
$\T (b')$ for 
a $U_q'(C_{n}^{(1)})$ crystal
element $b' \in B_{l'}$ \cite{HKOT},
that represents a one-row tableau with length $l'$.
Putting these two symbols together we have
\begin{equation}
\label{eq:tabtwistaxx}
\T (\omega(b))=\bbx{0}{x_\emptyset}\!
\overbrace{\fbox{$\vphantom{\ol{1}} 1 \cd 1$}}^{2x_1}\!
\fbox{$\vphantom{\ol{1}}\cd$}\!
\overbrace{\fbox{$\vphantom{\ol{1}}n \cd n$}}^{2x_n}\!
\overbrace{\fbox{$\vphantom{\ol{1}}\ol{n} \cd \ol{n}$}}^{2\ol{x}_n}\!
\fbox{$\vphantom{\ol{1}}\cd$}\!
\overbrace{\fbox{$\ol{1} \cd \ol{1}$}}^{2\ol{x}_1}\!
\overbrace{\fbox{$\ol{0} \cd \ol{0}$}}^{\ol{x}_\emptyset},
\end{equation}
where $x_\emptyset = \overline{x}_\emptyset 
= l-\sum_{i=1}^n (x_i + \overline{x}_i)$.
We shall use this tableau in our description of the 
combinatorial $R$ matrix (Rule \ref{rule:Cx}).

{}From now on we shall devote ourselves to describing several important
properties of the map $\omega$.
Our goal here is Lemma \ref{lem:1} that
our description of the combinatorial $R$ matrix relies on.
For this purpose
we also use the symbol $\omega$ for the following map for $\et{i},\ft{i}$.
\begin{align*}
\omega (\etd_i ) & = (\etd_i')^{2 - \delta_{i,0}} ,\\
\omega (\ftd_i ) & = (\ftd_i')^{2 - \delta_{i,0}}.
\end{align*}
Hereafter 
we attach prime $'$ to the notations for
the $U_q'(C_{n}^{(1)})$ crystals, e.g. $\etd_i', \varphi_i'$ and so on.
\begin{lemma}
\label{lem:01}
\begin{align*}
\omega (\etd_i b) & = \omega (\etd_i) \omega (b), \\
\omega (\ftd_i b) & = \omega (\ftd_i) \omega (b),
\end{align*}
i.e. the $\omega$ commutes with actions of the operators on $B_l$.
\end{lemma}
%
\noindent
Let us give a proof in $\etd_0$ case.
(For the other $\etd_i$s (and also for $\ftd_i$s) the proof is similar.)
\begin{proof}
Let $b=(x_1,\ldots,\ol{x}_l) \in B_l$ ($U_q'(A_{2n}^{(2)})$ crystal).
We have \cite{KKM}
\begin{displaymath}
\etd_0 b = 
\begin{cases}
(x_1-1,x_2,\ldots,\ol{x}_1) & \text{if $x_1 \geq \ol{x}_1+1$,} \\
(x_1,\ldots,\ol{x}_2,\ol{x}_1+1) & \text{if $x_1 \leq \ol{x}_1$.}
\end{cases}
\end{displaymath}
This means that
\begin{displaymath}
\omega (\etd_0 b) = 
\begin{cases}
(2x_1-2,2x_2,\ldots,2\ol{x}_1) & \text{if $2x_1 \geq 2\ol{x}_1+2$,} \\
(2x_1,\ldots,2\ol{x}_2,2\ol{x}_1+2) & \text{if $2x_1 \leq 2\ol{x}_1$.}
\end{cases}
\end{displaymath}
On the other hand, 
let $b'=(x'_1,\ldots,\ol{x}'_l) \in B'_{l'}$ ($U_q'(C_{n}^{(1)})$ crystal).
We have \cite{KKM}
\begin{displaymath}
\etd'_0 b' = 
\begin{cases}
(x'_1-2,x'_2,\ldots,\ol{x}'_1) & \text{if $x'_1 \geq \ol{x}'_1+2$,} \\
(x'_1-1,x'_2,\ldots,\ol{x}'_2,\ol{x}'_1+1) & \text{if $x'_1 = \ol{x}'_1+1$,} \\
(x'_1,\ldots,\ol{x}'_2,\ol{x}'_1+2) & \text{if $x'_1 \leq \ol{x}'_1$.}
\end{cases}
\end{displaymath}
Thus putting $l'=2l$ and $b'=\omega (b)$ we obtain
$\omega (\etd_0) \omega (b) = \etd'_0 b' = \omega (\etd_0 b) $.
(The second choice in the above equation does not occur.)
\end{proof}
Let us denote $\omega (b_1) \ot \omega (b_2)$ by $\omega(b_1 \ot b_2)$
for $b_1 \ot b_2 \in B_l \ot B_k$.
\begin{lemma}
\label{lem:02}
\begin{align*}
\omega (\etd_i (b_1 \ot b_2)) &= 
\omega (\etd_i) \omega(b_1 \ot b_2), \\
\omega (\ftd_i (b_1 \ot b_2)) &= 
\omega (\ftd_i) \omega(b_1 \ot b_2).
\end{align*}
Namely the $\omega$ commutes with actions of 
the operators on $B_l \ot B_k$.
\end{lemma}
\begin{proof}
Let us check the latter.
Suppose we have $\vphi_i(b_1) \geq \veps_i(b_2)+1$.
Then $\ftd_i (b_1 \ot b_2) = (\ftd_i b_1) \ot b_2$.
In this case we have $\vphi'_i (\omega (b_1)) 
\geq \veps'_i (\omega (b_2))+2-\delta_{i,0}$,
since
\begin{align*}
\varphi_i'(\omega (b)) & = (2 - \delta_{i,0}) \varphi_i (b) ,\\
\varepsilon_i'(\omega (b)) & = (2 - \delta_{i,0}) \varepsilon_i (b) .
\end{align*}
Therefore we obtain
\begin{align*}
\omega (\ftd_i) \omega(b_1 \ot b_2) &=
(\ftd_i')^{2-\delta_{i,0}} (\omega (b_1) \ot \omega (b_2)) \\
&= ((\ftd_i')^{2-\delta_{i,0}} \omega (b_1)) \ot \omega (b_2) \\
&= (\omega (\ftd_i) \omega (b_1)) \ot \omega (b_2) \\
&= \omega (\ftd_i b_1) \ot \omega (b_2) \\
&= \omega (\ftd_i b_1 \ot b_2).
\end{align*}
The other case when $\vphi_i(b_1)\le\veps_i(b_2)$ is similar.
\end{proof}
Finally we obtain the following important properties
of the map $\omega$.
\begin{lemma}
\label{lem:1}
\par\noindent
\begin{enumerate}
\renewcommand{\labelenumi}{(\roman{enumi})}
\item If $b_1 \ot b_2$ is mapped to $b'_2 \ot b'_1$ under the 
isomorphism of the $U_q'(A_{2n}^{(2)})$ crystals $B_l \ot B_k \simeq B_k \ot B_l$,
then $\omega(b_1) \ot \omega(b_2)$ is mapped to 
$\omega(b'_2) \ot \omega(b'_1)$ under the 
isomorphism of the $U_q'(C_{n}^{(1)})$ crystals 
$B_{2l} \ot B_{2k} \simeq B_{2k} \ot B_{2l}$.
\item Up to a global additive constant, 
the value of the energy function $H_{B_lB_k}(b_1 \ot b_2)$ 
for the $U_q'(A_{2n}^{(2)})$ crystal $B_l \ot B_k$
is equal to the value of the energy function 
$H'_{B_{2l}B_{2k}}(\omega(b_1) \ot \omega(b_2))$ for the $U_q'(C_{n}^{(1)})$ crystal
$B_{2l} \ot B_{2k}$.
\end{enumerate}
\end{lemma}
\begin{proof}
First we consider (i).
Since the crystal graph of $B_l \ot B_k$ is connected,
it remains to check (i) for any specific element in $B_l \ot B_k \simeq
B_k \ot B_l$.
We can do it by taking
$(l,0,\ldots,0)\ot (k,0,\ldots,0) \stackrel{\sim}{\mapsto}
(k,0,\ldots,0)\ot (l,0,\ldots,0) $ as the specific element,
for which (i) certainly holds.

We proceed to (ii).
We can set
\begin{displaymath}
H_{B_l B_k}((l,0,\ldots,0)\ot (k,0,\ldots,0) ) = 
H'_{B_{2l} B_{2k}}(\omega((l,0,\ldots,0))\ot \omega((k,0,\ldots,0)) ).
\end{displaymath}
Suppose $\etd_i (b_1 \ot b_2) \ne 0$.
Recall the defining relations of 
the energy function $H_{B_l B_k}$.
\begin{displaymath}
H_{B_l B_k} (\etd_i (b_1 \ot b_2)) =
\begin{cases} 
H_{B_l B_k} (b_1 \ot b_2)+1 & 
\text{if $i=0, \varphi_0(b_1) \geq \varepsilon_0(b_2), \varphi_0(b_2') \geq \varepsilon_0(b_1')$,} \\
H_{B_l B_k} (b_1 \ot b_2)-1 & 
\text{if $i=0, \varphi_0(b_1) < \varepsilon_0(b_2), \varphi_0(b_2') < \varepsilon_0(b_1')$,} \\
H_{B_l B_k} (b_1 \ot b_2) & 
\text{otherwise.}
\end{cases}
\end{displaymath}
Claim (ii) holds if
for any $i$ and $b_1 \ot b_2$ with $\etd_i (b_1 \ot b_2) \ne 0$, we have
\begin{eqnarray}
&& H'_{B_{2l} B_{2k}}(\omega(\etd_i (b_1 \ot b_2)))-
H'_{B_{2l} B_{2k}}(\omega(b_1 \ot b_2)) \nonumber\\
&& \quad =
H_{B_l B_k} (\etd_i (b_1 \ot b_2))-
H_{B_l B_k} (b_1 \ot b_2).
\label{eq:hfuncdif}
\end{eqnarray}
The $i=0$ case is verified as follows.
Since $\omega$ commutes with crystal actions 
we have $\omega (\etd_0 (b_1 \ot b_2)) =
\omega (\etd_0) (\omega(b_1) \ot \omega(b_2)) =
\etd'_0 (\omega(b_1 \ot b_2))$.
On the other hand $\omega$ 
preserves the inequalities in the classification
conditions in the defining relations of the energy 
function, i.e. $\varphi'_0(\omega (b_1)) \geq 
\varepsilon'_0(\omega(b_2)) \Leftrightarrow 
\varphi_0(b_1) \geq \varepsilon_0(b_2)$, and so on.
Thus (\ref{eq:hfuncdif}) follows from the defining relations of the 
$H'_{B_{2l} B_{2k}}$.
The $i\neq0$ case is easier.
This completes the proof.
\end{proof}
Since $\omega$ is injective we obtain the converse of (i).
\begin{coro}
If $\omega(b_1) \ot \omega(b_2)$ is mapped to 
$\omega(b'_2) \ot \omega(b'_1)$ under the 
isomorphism of the $U_q'(C_{n}^{(1)})$ crystals 
$B_{2l} \ot B_{2k} \simeq B_{2k} \ot B_{2l}$, then
$b_1 \ot b_2$ is mapped to $b'_2 \ot b'_1$ under the 
isomorphism of the $U_q'(A_{2n}^{(2)})$ crystals $B_l \ot B_k \simeq B_k \ot B_l$.
\end{coro}
\subsection{\mathversion{bold}Definitions : $U_q'(D_{n+1}^{(2)})$ crystal case}
\label{subsec:twistd}
Given a positive integer $l$,
we consider a $U_q'(D_{n+1}^{(2)})$ crystal denoted by $B_l$
that is defined
in \cite{KKM}.
$B_l$'s are the crystal bases
of the irreducible finite-dimensional representation
of the quantum affine algebra $U_q'(D_{n+1}^{(2)})$.
As a set $B_{l}$ reads
$$
B_{l} = \left\{(
x_1,\ldots, x_n,x_{\circ},\overline{x}_n,\ldots,\overline{x}_1) \Biggm|
\begin{array}{l}
x_{\circ}=\mbox{$0$ or $1$}, \quad x_i, \overline{x}_i \in \Z_{\ge 0}\\
x_{\circ}+\sum_{i=1}^n(x_i + \overline{x}_i) \in \{l,l-1,\ldots ,0\} 
\end{array}
\right\}.
$$
For its crystal structure see \cite{KKM}.
$B_{l}$ is isomorphic to
$\bigoplus_{0 \leq j \leq l} B(j \La_1)$ as a $U_q(B_n)$ crystal, where $B(j \La_1)$ is
the crystal associated with the irreducible representation of $U_q(B_n)$
with highest weight $j \Lambda_1$.
The $U_q(B_n)$ crystal $B(j \La_1)$ has a description in terms of semistandard $B$-tableaux \cite{KN}.
The entries are $1,\ldots ,n$, $\ol{1}, \ldots ,\ol{n}$ and $\circ$ with the
total order:
\begin{displaymath}
1 < 2 < \cd < n < \circ < \ol{n} < \cd < \ol{2} < \ol{1}.
\end{displaymath}
In this paper we use the symbol $\circ$ for the member of entries of
the semistandard $B$ tableaux that is conventionally denoted by $0$.
For $b= (x_1, \ldots, x_n, x_{\circ}, \overline{x}_n,\ldots,\overline{x}_1) \in
B(j \La_1)$
the tableau $\mathcal{T}(b)$ is depicted by
\begin{displaymath}
\mathcal{T}(b)=\overbrace{\fbox{$\vphantom{\ol{1}} 1 \cd 1$}}^{x_1}\!
\fbox{$\vphantom{\ol{1}}\cd$}\!
\overbrace{\fbox{$\vphantom{\ol{1}}n \cd n$}}^{x_n}\!
\overbrace{\fbox{$\vphantom{\ol{1}}\hphantom{1}\circ\hphantom{1}$}}^{x_{\circ}}\!
\overbrace{\fbox{$\vphantom{\ol{1}}\ol{n} \cd \ol{n}$}}^{\ol{x}_n}\!
\fbox{$\vphantom{\ol{1}}\cd$}\!
\overbrace{\fbox{$\ol{1} \cd \ol{1}$}}^{\ol{x}_1}.
\end{displaymath}
The length of this one-row tableau is equal to $j$, namely
$x_{\circ}+\sum_{i=1}^n(x_i + \overline{x}_i) =j$.

To describe our rule for the combinatorial $R$ matrix
we shall depict the elements of $B_{l}$ 
by one-row $C$ tableaux with length $2l$.
We introduce the map $\omega$ from the $U_q'(D_{n+1}^{(2)})$ 
crystal $B_l$ to the $U_q'(C_{n}^{(1)})$ crystal $B_{2l}$.
$\omega$ sends
$b= (x_1, \ldots, x_n, \overline{x}_n,\ldots,\overline{x}_1)$ to
$\omega (b) = (2x_1, \ldots, 2x_{n-1},$ $2x_n + x_{\circ}, 2\overline{x}_n + x_{\circ},2\overline{x}_{n-1},\ldots,2\overline{x}_1)$.
By using the symbol $\T$
introduced in the previous subsection
the tableau $\T (\omega(b))$ is depicted by
\begin{displaymath}
\T (\omega (b))=\bbx{0}{x_\emptyset}\!
\overbrace{\fbox{$\vphantom{\ol{1}} 1 \cd 1$}}^{2x_1}\!
\fbox{$\vphantom{\ol{1}}\cd$}\!
\overbrace{\fbox{$\vphantom{\ol{1}}n \cd n$}}^{2x_n+x_{\circ}}\!
\overbrace{\fbox{$\vphantom{\ol{1}}\ol{n} \cd \ol{n}$}}^{2\ol{x}_n+x_{\circ}}\!
\fbox{$\vphantom{\ol{1}}\cd$}\!
\overbrace{\fbox{$\ol{1} \cd \ol{1}$}}^{2\ol{x}_1}\!
\overbrace{\fbox{$\ol{0} \cd \ol{0}$}}^{\ol{x}_\emptyset},
\end{displaymath}
where $x_\emptyset = \overline{x}_\emptyset 
= l-x_{\circ}-\sum_{i=1}^n (x_i + \overline{x}_i)$.

Our description of the combinatorial $R$ matrix (Theorem \ref{th:main1}) 
is based on the following lemma.
\begin{lemma}
\label{lem:2}
\par\noindent
\begin{enumerate}
\renewcommand{\labelenumi}{(\roman{enumi})}
\item If $b_1 \ot b_2$ is mapped to $b'_2 \ot b'_1$ under the 
isomorphism of the $U_q'(D_{n+1}^{(2)})$ crystals $B_l \ot B_k \simeq B_k \ot B_l$,
then $\omega(b_1) \ot \omega(b_2)$ is mapped to $\omega(b'_2) \ot 
\omega(b'_1)$ under the 
isomorphism of the $U_q'(C_{n}^{(1)})$ crystals 
$B_{2l} \ot B_{2k} \simeq B_{2k} \ot B_{2l}$.
\item Up to a global additive constant, 
the value of the energy function $H_{B_lB_k}(b_1 \ot b_2)$ 
for the $U_q'(D_{n+1}^{(2)})$ crystal $B_l \ot B_k$
is equal to the value of the energy function 
$H'_{B_{2l}B_{2k}}(\omega(b_1) \ot \omega(b_2))$ for the $U_q'(C_{n}^{(1)})$ crystal
$B_{2l} \ot B_{2k}$.
\end{enumerate}
\end{lemma}
\begin{proof}
To distinguish the notations
we attach prime $'$ to the notations for
the $U_q'(C_{n}^{(1)})$ crystals.
Then we have
\begin{align*}
\varphi_i'(\omega (b)) & = (2 - \delta_{i,0}- \delta_{i,n}) \varphi_i (b) ,\\
\varepsilon_i'(\omega (b)) & = (2 - \delta_{i,0}- \delta_{i,n}) \varepsilon_i (b) .
\end{align*}
Let us define the action of $\omega$ on the operators by
\begin{align*}
\omega (\etd_i ) & = (\etd_i')^{2 - \delta_{i,0}- \delta_{i,n}} ,\\
\omega (\ftd_i ) & = (\ftd_i')^{2 - \delta_{i,0}- \delta_{i,n}}.
\end{align*}
By repeating an argument similar to that in the previous subsection
we obtain formally the same assertions of Lemmas
\ref{lem:01}, \ref{lem:02} and \ref{lem:1}.
This completes the proof.
\end{proof}
Since $\omega$ is injective we obtain the converse of (i).
\begin{coro}
If $\omega(b_1) \ot \omega(b_2)$ is mapped to $\omega(b'_2) \ot 
\omega(b'_1)$ under the 
isomorphism of the $U_q'(C_{n}^{(1)})$ crystals 
$B_{2l} \ot B_{2k} \simeq B_{2k} \ot B_{2l}$, then
$b_1 \ot b_2$ is mapped to $b'_2 \ot b'_1$ under the 
isomorphism of the $U_q'(D_{n+1}^{(2)})$ crystals $B_l 
\ot B_k \simeq B_k \ot B_l$.
\end{coro}
\subsection{\mathversion{bold} Column insertion and inverse insertion for
$C_n$}
\label{subsec:ccis}

Set an alphabet $\mathcal{X}=\mathcal{A} \sqcup \bar{\mathcal{A}},\,
\mathcal{A}=\{0, 1,\dots,n\}$ and
$\bar{\mathcal{A}}=\{\overline{0}, \overline{1},\dots,\overline{n}\}$,
with the total order
$0 < 1 < 2 < \dots < n < \overline{n} < \dots < 
\overline{2} < \overline{1} < \overline{0}$.\footnote{We also introduce $0$ and $\ol{0}$ in the 
alphabet. Compare with \cite{B1,KN}.}
%
\subsubsection{Semi-standard $C$ tableaux}
\label{subsubsec:ssct}
Let us consider a {\em semistandard $C$ tableau}
made by the letters from this alphabet.
We follow \cite{KN} for its definition.
We present the definition here, but restrict ourselves to
the special cases that are sufficient for our purpose.
Namely we consider only
those tableaux that have no more than two rows
in their shapes. 
Thus they have the forms as
\begin{equation}
\label{eq:pictureofsst}
\setlength{\unitlength}{5mm}
\begin{picture}(6,1.5)(1.5,0)
\put(1,0){\line(1,0){5}}
\put(1,1){\line(1,0){5}}
\put(1,0){\line(0,1){1}}
\put(2,0){\line(0,1){1}}
\put(3,0){\line(0,1){1}}
\put(5,0){\line(0,1){1}}
\put(6,0){\line(0,1){1}}
\put(1,0){\makebox(1,1){$\alpha_1$}}
\put(2,0){\makebox(1,1){$\alpha_2$}}
\put(3,0){\makebox(2,1){$\cd$}}
\put(5,0){\makebox(1,1){$\alpha_j$}}
\end{picture}
\mbox{or}
\setlength{\unitlength}{5mm}
\begin{picture}(12,2.5)
\put(2,0){\line(1,0){5}}
\put(2,1){\line(1,0){10}}
\put(2,2){\line(1,0){10}}
\put(2,0){\line(0,1){2}}
\put(3,0){\line(0,1){2}}
\put(4,0){\line(0,1){2}}
\put(6,0){\line(0,1){2}}
\put(7,0){\line(0,1){2}}
\put(9,1){\line(0,1){1}}
\put(11,1){\line(0,1){1}}
\put(12,1){\line(0,1){1}}
\put(2,0){\makebox(1,1){$\beta_1$}}
\put(2,1){\makebox(1,1){$\alpha_1$}}
\put(3,0){\makebox(1,1){$\beta_2$}}
\put(3,1){\makebox(1,1){$\alpha_2$}}
\put(4,0){\makebox(2,1){$\cd$}}
\put(4,1){\makebox(2,1){$\cd$}}
\put(6,0){\makebox(1,1){$\beta_i$}}
\put(6,1){\makebox(1,1){$\alpha_i$}}
\put(7,1){\makebox(2,1){$\alpha_{i+1}$}}
\put(9,1){\makebox(2,1){$\cd$}}
\put(11,1){\makebox(1,1){$\alpha_j$}}
\end{picture},
\end{equation}
and the letters inside the boxes should obey the following conditions:
\begin{equation}
\label{eq:notdecrease}
\alpha_1 \leq \cdots \leq \alpha_j,\quad \beta_1 \leq \cdots \leq \beta_i,
\end{equation}
\begin{equation}
\label{eq:strictincrease}
\alpha_a < \beta_a,
\end{equation}
\begin{equation}
\label{eq:zerozerobar}
(\alpha_a,\beta_a) \ne (0,\ol{0}),
\end{equation}
\begin{equation}
\label{eq:absenceofxxconf}
(\alpha_a,\alpha_{a+1},\beta_{a+1}) \ne (x,x,\ol{x}),\quad
(\alpha_a,\beta_{a},\beta_{a+1}) \ne (x,\ol{x},\ol{x}).
\end{equation}
Here we assume $1 \leq x \leq n$.
The last conditions (\ref{eq:absenceofxxconf}) are referred to as the
absence of the $(x,x)$-configurations.
\subsubsection{\mathversion{bold}Column insertion for 
$C_n$ \cite{B1}}
\label{subsubsec:insc}
We give a list of column insertions on  
semistandard $C$ tableaux that are sufficient for our 
purpose (Rule \ref{rule:Cx}).
First of all let us explain the relation between the {\em insertion} and 
the {\em inverse insertion}.
Since we are deliberately avoiding the occurrence of the
{\em bumping sliding transition}
(\cite{B1}), the situation is basically the same as that for the usual
tableau case (\cite{F}, Appendix A.2).
Namely, when a letter $\alpha$ was inserted into the tableau $T$,
we obtain a new tableau $T'$ whose shape is one more box larger than
the shape of $T$.
If we know the location of the new box we can reverse the insertion
process to retrieve the original tableau $T$ and letter $\alpha$.
This is the inverse insertion process.
These processes go on column by column.
Thus, from now on let us pay our attention to a particular column $C$ in the
tableau.
Suppose we have inserted a letter $\alpha$ into $C$.
Suppose then we have
obtained a column $C'$ and a letter $\alpha'$
bumped out from the column.
If we inversely insert the letter $\alpha'$ into the column
$C'$, we retrieve the original column $C$ and letter $\alpha$.

For the alphabet $\mathcal{X}$, we follow the convention
that Greek letters $ \alpha, \beta, \ldots $ belong to
$\mathcal{X}$ while Latin
letters $x,y,\ldots$ (resp. $\overline{x},\overline{y},\ldots$) belong to
$\mathcal{A}$ (resp. $\bar{\mathcal{A}}$).
The pictorial equations in the list should be interpreted as follows.
(We take up two examples.)
\begin{itemize}
	\item In Case B0, the letter\footnote{By abuse of notation 
we identify a letter with the one-box tableau having the letter in it.}
$\alpha$ is inserted into the column
with only one box that has letter $\beta$ in it.
The $\alpha$ is set in the box and the $\beta$ is bumped out to
the right-hand column.
	\item In Case B1, the letter $\beta$ is inserted into the column
with two boxes that have letters $\alpha$ and $\gamma$ in them.
The $\beta$ is set in the lower box and the $\gamma$ is bumped out to
the right-hand column.
\end{itemize}
Other equations should be interpreted in a similar way\footnote{This interpretation 
is also eligible for the lists of
type $B$ and $D$ cases in Sections \ref{subsubsec:insb} and 
\ref{subsubsec:insd}.}.
We note that there is no overlapping case in the list.
Note also that it does not exhaust all patterns of the column insertions
that insert a letter into a column with at most two boxes.
For instance it does not cover the case of insertion
\begin{math}
\setlength{\unitlength}{3mm}
\begin{picture}(3,2)(0,0.3)
\multiput(0,0)(1,0){2}{\line(0,1){1}}
\multiput(0,0)(0,1){2}{\line(1,0){1}}
\put(0,0){\makebox(1,1){${\scriptstyle \ol{2}}$}}
\put(1,0){\makebox(1,1){${\scriptstyle \rightarrow}$}}
\multiput(2,0)(1,0){2}{\line(0,1){2}}
\multiput(2,0)(0,1){3}{\line(1,0){1}}
\put(2,0){\makebox(1,1){${\scriptstyle 2}$}}
\put(2,1){\makebox(1,1){${\scriptstyle 1}$}}
\end{picture}
\end{math}.
In Rule \ref{rule:Cx} we do not encounter such a case.
\par
\noindent
\ToBOX{A0}{\alpha}%
	\raisebox{1.25mm}{,}

\noindent
\ToDOMINO{A1}{\alpha}{\beta}{\alpha}{\beta}%
	\raisebox{4mm}{if $\alpha < \beta$ ,}

\noindent
\ToYOKODOMINO{B0}{\beta}{\alpha}{\beta}{\alpha}%
	\raisebox{1.25mm}{if $\alpha \le \beta$,}

\noindent
\ToHOOKnn{B1}{\alpha}{\gamma}{\beta}{\gamma}{\alpha}{\beta}%
	\raisebox{4mm}{if $\alpha < \beta \leq \gamma$ and $(\alpha,\gamma) \ne (x,\overline{x})$,}

\noindent
\ToHOOKnn{B2}{\beta}{\gamma}{\alpha}{\beta}{\alpha}{\gamma}
	\raisebox{4mm}{if $\alpha \leq \beta < \gamma$ and $(\alpha,\gamma) \ne (x,\overline{x})$,}

\noindent
\ToHOOKll{B3}{x}{\overline{x}}{\beta}{\overline{x\!-\!1}}{x\!-\!1}{\beta}%
	\raisebox{4mm}{if $x \le \beta \le \overline{x}$ and $x\ne 0$,}

\noindent
\ToHOOKln{B4}{\beta}{\overline{x}}{x}{\beta}{x\!+\!1}{\overline{x\!+\!1}}%
	\raisebox{4mm}{if $x < \beta < \overline{x}$ and $x\ne n$.}

\noindent
\subsubsection{Column insertion and $U_q (\ol{\geh})$ crystal morphism}
In this subsection we illustrate the relation between column
insertion and the crystal morphism that was given by Baker \cite{B1,B2}.
A crystal morphism is a (not necessarily one-to-one)
map between two crystals
that commutes with the actions of crystals.
See, for instance \cite{KKM} for a precise definition.
A $U_q (\ol{\geh})$ crystal morphism is a morphism
that commutes with the actions of 
$\etd_i$ and $\ftd_i$ for $i \ne 0$.
For later use we also include 
semistandard $B$ and $D$ tableaux in our discussion,
therefore we assume $\ol{\geh} =B_n, C_n$ or $D_n$.
See section \ref{subsubsec:ssbt} (resp.~\ref{subsubsec:ssdt})
for the definition of semistandard $B$ (resp.~$D$) tableaux.

Let $T$ be a semistandard $B$, $C$ or $D$ tableau.
For this $T$ we denote by $w(T)$ the
Japanese reading word of $T$, i.e. $w(T)$ is a sequence of letters that is created by
reading all letters on $T$ from the rightmost column
to the leftmost one,
and in each column, from top to bottom.
For instance,
\par\noindent
\setlength{\unitlength}{5mm}
\begin{picture}(22,1.5)(-5,0)
\put(0,0){\makebox(1,1){$w($}}
\put(1,0){\line(1,0){5}}
\put(1,1){\line(1,0){5}}
\put(1,0){\line(0,1){1}}
\put(2,0){\line(0,1){1}}
\put(3,0){\line(0,1){1}}
\put(5,0){\line(0,1){1}}
\put(6,0){\line(0,1){1}}
\put(1,0){\makebox(1,1){$\alpha_1$}}
\put(2,0){\makebox(1,1){$\alpha_2$}}
\put(3,0){\makebox(2,1){$\cd$}}
\put(5,0){\makebox(1,1){$\alpha_j$}}
\put(6,0){\makebox(6,1){$)=\alpha_j \cd \alpha_2 \alpha_1,$}}
\end{picture}
\par\noindent
\setlength{\unitlength}{5mm}
\begin{picture}(22,2.5)(-2,0)
\put(0,0){\makebox(2,2){$w \Biggl($}}
\put(2,0){\line(1,0){5}}
\put(2,1){\line(1,0){10}}
\put(2,2){\line(1,0){10}}
\put(2,0){\line(0,1){2}}
\put(3,0){\line(0,1){2}}
\put(4,0){\line(0,1){2}}
\put(6,0){\line(0,1){2}}
\put(7,0){\line(0,1){2}}
\put(9,1){\line(0,1){1}}
\put(11,1){\line(0,1){1}}
\put(12,1){\line(0,1){1}}
\put(2,0){\makebox(1,1){$\beta_1$}}
\put(2,1){\makebox(1,1){$\alpha_1$}}
\put(3,0){\makebox(1,1){$\beta_2$}}
\put(3,1){\makebox(1,1){$\alpha_2$}}
\put(4,0){\makebox(2,1){$\cd$}}
\put(4,1){\makebox(2,1){$\cd$}}
\put(6,0){\makebox(1,1){$\beta_i$}}
\put(6,1){\makebox(1,1){$\alpha_i$}}
\put(7,1){\makebox(2,1){$\alpha_{i+1}$}}
\put(9,1){\makebox(2,1){$\cd$}}
\put(11,1){\makebox(1,1){$\alpha_j$}}
\put(12,0){\makebox(14,2){$\Biggr) =\alpha_j \cd \alpha_{i+1}
\alpha_i \beta_i \cd
\alpha_2 \beta_2 \alpha_1 \beta_1.$}}
\end{picture}
\par\noindent
Let $T$ and $T'$ be two tableaux.
We define the product tableau $T * T'$ by
\begin{displaymath}
T * T' = (\tau_1 \to \cd (\tau_{j-1} \to ( \tau_j \to T ) ) \cd )
\end{displaymath}
where
\begin{displaymath}
w(T') = \tau_j \tau_{j-1} \cd \tau_1.
\end{displaymath}
The symbol $\to$ represents the column insertions in \cite{B1,B2} which
we partly describe in sections \ref{subsubsec:insc}, 
\ref{subsubsec:insb} and \ref{subsubsec:insd}.
(Note that the author of \cite{B1,B2} uses $\leftarrow$ instead of $\to$.)

For a dominant integral weight $\lambda$ of the $\ol{\geh}$ root system,
let $B(\lambda)$ be the $U_q(\ol{\geh})$ crystal
associated with the irreducible highest weight representation $V(\lambda)$.
The elements of $B(\lambda)$ can be represented by
semistandard $\ol{\geh}$ tableaux of shape $\lambda$ \cite{KN}.
\begin{proposition}[\cite{B1,B2}]
\label{pr:morphgen}
Let $B(\mu) \ot B(\nu) \simeq \bigoplus_j B(\lambda_j)^{\oplus m_j}$
be the tensor product decomposition of crystals. Here $\la_j$'s are
distinct highest weights and $m_j(\ge1)$ is the multiplicity of $B(\la_j)$.
Forgetting the multiplicities we have the canonical morphism from
$B(\mu) \ot B(\nu)$ to $\bigoplus_j B(\lambda_j)$.
Define $\psi$ by
\begin{displaymath}
\psi(b_1 \ot b_2) = b_1 * b_2.
\end{displaymath}
Then $\psi$ gives the unique $U_q(\ol{\geh})$ crystal morphism from
$B(\mu) \ot B(\nu)$ to $\bigoplus_j B(\lambda_j)$.
\end{proposition}
\noindent
See Examples \ref{ex:morC1}, \ref{ex:morC2}, \ref{ex:morB}, \ref{ex:morD1}
and \ref{ex:morD2}.
\subsubsection{Column insertion and $U_q (C_n)$ crystal morphism}
To illustrate Proposition \ref{pr:morphgen},
let us check a morphism of the $U_q(C_2)$ crystal
$B(\Lambda_2) \ot B(\Lambda_1)$ by taking two examples.
Let $\psi$ be the map that sends
\begin{math}
\setlength{\unitlength}{3mm}
\begin{picture}(3,2)(0,0.3)
\multiput(0,0)(1,0){2}{\line(0,1){2}}
\multiput(0,0)(0,1){3}{\line(1,0){1}}
\put(0,0){\makebox(1,1){${\scriptstyle \beta}$}}
\put(0,1){\makebox(1,1){${\scriptstyle \alpha}$}}
\put(1,0.5){\makebox(1,1){${\scriptstyle \otimes}$}}
\multiput(2,0.5)(1,0){2}{\line(0,1){1}}
\multiput(2,0.5)(0,1){2}{\line(1,0){1}}
\put(2,0.5){\makebox(1,1){${\scriptstyle \gamma}$}}
\end{picture}
\end{math}
to the tableau which is made by the column insertion
\begin{math}
\setlength{\unitlength}{3mm}
\begin{picture}(3,2)(0,0.3)
\multiput(0,0)(1,0){2}{\line(0,1){1}}
\multiput(0,0)(0,1){2}{\line(1,0){1}}
\put(0,0){\makebox(1,1){${\scriptstyle \gamma}$}}
\put(1,0){\makebox(1,1){${\scriptstyle \rightarrow}$}}
\multiput(2,0)(1,0){2}{\line(0,1){2}}
\multiput(2,0)(0,1){3}{\line(1,0){1}}
\put(2,0){\makebox(1,1){${\scriptstyle \beta}$}}
\put(2,1){\makebox(1,1){${\scriptstyle \alpha}$}}
\end{picture}
\end{math}.
\begin{example}
\label{ex:morC1}
\begin{displaymath}
\begin{CD}
\setlength{\unitlength}{5mm}
\begin{picture}(3,2)(0,0.3)
\multiput(0,0)(1,0){2}{\line(0,1){2}}
\multiput(0,0)(0,1){3}{\line(1,0){1}}
\put(0,0){\makebox(1,1){$\ol{2}$}}
\put(0,1){\makebox(1,1){$2$}}
\put(1,0.5){\makebox(1,1){$\otimes$}}
\multiput(2,0.5)(1,0){2}{\line(0,1){1}}
\multiput(2,0.5)(0,1){2}{\line(1,0){1}}
\put(2,0.5){\makebox(1,1){$2$}}
\end{picture}
@>\text{$\etd_1$}>>
\setlength{\unitlength}{5mm}
\begin{picture}(3,2)(0,0.3)
\multiput(0,0)(1,0){2}{\line(0,1){2}}
\multiput(0,0)(0,1){3}{\line(1,0){1}}
\put(0,0){\makebox(1,1){$\ol{2}$}}
\put(0,1){\makebox(1,1){$1$}}
\put(1,0.5){\makebox(1,1){$\otimes$}}
\multiput(2,0.5)(1,0){2}{\line(0,1){1}}
\multiput(2,0.5)(0,1){2}{\line(1,0){1}}
\put(2,0.5){\makebox(1,1){$2$}}
\end{picture}
\\
@VV\text{$\psi$}V @VV\text{$\psi$}V \\
\setlength{\unitlength}{5mm}
\begin{picture}(2,2)(0,0.3)
\multiput(0,0)(1,0){2}{\line(0,1){2}}
\put(0,0){\line(1,0){1}}
\multiput(0,1)(0,1){2}{\line(1,0){2}}
\put(2,1){\line(0,1){1}}
\put(0,0){\makebox(1,1){$2$}}
\put(0,1){\makebox(1,1){$1$}}
\put(1,1){\makebox(1,1){$\ol{1}$}}
\end{picture}
@>\text{$\etd_1$}>>
\setlength{\unitlength}{5mm}
\begin{picture}(2,2)(0,0.3)
\multiput(0,0)(1,0){2}{\line(0,1){2}}
\put(0,0){\line(1,0){1}}
\multiput(0,1)(0,1){2}{\line(1,0){2}}
\put(2,1){\line(0,1){1}}
\put(0,0){\makebox(1,1){$2$}}
\put(0,1){\makebox(1,1){$1$}}
\put(1,1){\makebox(1,1){$\ol{2}$}}
\end{picture}
\end{CD}
\end{displaymath}
\vskip3ex
\noindent
Here the left (resp.~right) $\psi$ is given by Case B3 (resp.~B1) column insertion.
\end{example}
\begin{example}
\label{ex:morC2}
\begin{displaymath}
\begin{CD}
\setlength{\unitlength}{5mm}
\begin{picture}(3,2)(0,0.3)
\multiput(0,0)(1,0){2}{\line(0,1){2}}
\multiput(0,0)(0,1){3}{\line(1,0){1}}
\put(0,1){\makebox(1,1){$2$}}
\put(0,0){\makebox(1,1){$\ol{1}$}}
\put(1,0.5){\makebox(1,1){$\otimes$}}
\put(2,0.5){\makebox(1,1){$1$}}
\multiput(2,0.5)(1,0){2}{\line(0,1){1}}
\multiput(2,0.5)(0,1){2}{\line(1,0){1}}
\end{picture}
@>\text{$\ftd_1$}>>
\setlength{\unitlength}{5mm}
\begin{picture}(3,2)(0,0.3)
\multiput(0,0)(1,0){2}{\line(0,1){2}}
\multiput(0,0)(0,1){3}{\line(1,0){1}}
\put(0,1){\makebox(1,1){$2$}}
\put(0,0){\makebox(1,1){$\ol{1}$}}
\put(1,0.5){\makebox(1,1){$\otimes$}}
\put(2,0.5){\makebox(1,1){$2$}}
\multiput(2,0.5)(1,0){2}{\line(0,1){1}}
\multiput(2,0.5)(0,1){2}{\line(1,0){1}}
\end{picture}
\\
@VV\text{$\psi$}V @VV\text{$\psi$}V \\
\setlength{\unitlength}{5mm}
\begin{picture}(2,2)(0,0.3)
\multiput(0,0)(1,0){2}{\line(0,1){2}}
\put(0,0){\line(1,0){1}}
\multiput(0,1)(0,1){2}{\line(1,0){2}}
\put(2,1){\line(0,1){1}}
\put(0,1){\makebox(1,1){$2$}}\put(1,1){\makebox(1,1){$2$}}
\put(0,0){\makebox(1,1){$\ol{2}$}}
\end{picture}
@>\text{$\ftd_1$}>>
\setlength{\unitlength}{5mm}
\begin{picture}(2,2)(0,0.3)
\multiput(0,0)(1,0){2}{\line(0,1){2}}
\put(0,0){\line(1,0){1}}
\multiput(0,1)(0,1){2}{\line(1,0){2}}
\put(2,1){\line(0,1){1}}
\put(0,1){\makebox(1,1){$2$}}\put(1,1){\makebox(1,1){$2$}}
\put(0,0){\makebox(1,1){$\ol{1}$}}
\end{picture}
\end{CD}
\end{displaymath}
\vskip3ex
\noindent
Here the left (resp.~right) $\psi$ is given by Case B4 (resp.~B2) column insertion.
\end{example}

\subsubsection{\mathversion{bold}Inverse insertion for 
$C_n$ \cite{B1}}
\label{subsubsec:invinsc}
In this subsection we give a list of inverse column 
insertions on semistandard $C$ tableaux that are sufficient 
for our purpose (Rule \ref{rule:Cx}).
The pictorial equations in the list should be interpreted as follows.
(We take two examples.)
\begin{itemize}
	\item In Case C0, the letter $\beta$ is inversely inserted into the column
with only one box that has letter $\alpha$ in it.
The $\beta$ is set in the box and the $\alpha$ is bumped out to
the left-hand column.
	\item In Case C1, the letter $\gamma$ is inversely inserted into the column
with two boxes that have letters $\alpha$ and $\beta$ in them.
The $\gamma$ is set in the lower box and the $\beta$ is bumped out to
the left-hand column.
\end{itemize}
Other equations illustrate analogous procedures.
%
\par
\noindent
\FromYOKODOMINO{C0}{\beta}{\alpha}{\beta}{\alpha}%
	\raisebox{1.25mm}{if $\alpha \le \beta$,}

\noindent
\FromHOOKnn{C1}{\gamma}{\alpha}{\beta}{\alpha}{\gamma}{\beta}%
	\raisebox{4mm}{if $\alpha < \beta \leq \gamma$ and $(\alpha,\gamma) \ne (x,\overline{x})$,}

\noindent
\FromHOOKnn{C2}{\beta}{\alpha}{\gamma}{\beta}{\gamma}{\alpha}%
	\raisebox{4mm}{if $\alpha \leq \beta < \gamma$ and $(\alpha,\gamma) \ne (x,\overline{x})$,}

\noindent
\FromHOOKnl{C3}{\overline{x}}{x}{\beta}{x\!+\!1}{\overline{x\!+\!1}}{\beta}%
	\raisebox{4mm}{if $x < \beta < \overline{x}$ and $x\ne n$,}

\noindent
\FromHOOKll{C4}{\beta}{x}{\overline{x}}{\beta}{\overline{x\!-\!1}}{x\!-\!1}%
	\raisebox{4mm}{if $x \le \beta \le \overline{x}$ and $x\ne 0$.}


\subsection{\mathversion{bold}Main theorem : $A^{(2)}_{2n}$ and 
$D^{(2)}_{n+1}$ cases}
\label{subsec:ruleCx}
Fix $l, k \in \Z_{\ge 1}$.
Given $b_1 \otimes b_2 \in B_{l} \otimes B_{k}$,
we define an $U'_q(C^{(1)}_n)$ crystal element
$\tilde{b}_2 \otimes \tilde{b}_1 \in B_{2k} \otimes B_{2l}$
and $l',k', m \in \Z_{\ge 0}$ by the following rule.

\begin{rules}\label{rule:Cx}
\hfill\par\noindent
Set $z = \min(\sharp\,\fbx{0} \text{ in }\T(\omega(b_1)),\, 
\sharp\,\fbx{0} \text{ in }\T(\omega(b_2)))$.
Thus $\T(\omega(b_1))$ and $\T(\omega(b_2))$ can be depicted by
\begin{eqnarray*}
\T(\omega (b_1)) &=&
\setlength{\unitlength}{5mm}
\begin{picture}(10.5,1.4)(0,0.3)
\multiput(0,0)(0,1){2}{\line(1,0){10}}
\put(0,0){\line(0,1){1}}
\put(3,0){\line(0,1){1}}
\put(7,0){\line(0,1){1}}
\put(10,0){\line(0,1){1}}
\put(0,0){\makebox(3,1){$0\cdots 0$}}
\put(3,0){\makebox(4,1){$T_*$}}
\put(7,0){\makebox(3,1){$\ol{0}\cdots \ol{0}$}}
\multiput(0,0.9)(7,0){2}{\put(0,0){\makebox(3,1){$z$}}}
\end{picture},\\
\T(\omega(b_2)) &=&
\setlength{\unitlength}{5mm}
\begin{picture}(9.5,2)(0,0.3)
\multiput(0,0)(0,1){2}{\line(1,0){9}}
\multiput(0,0)(3,0){4}{\line(0,1){1}}
\put(3.9,0){\line(0,1){1}}
\put(5,0){\line(0,1){1}}
\put(0,0){\makebox(3,1){$0\cdots 0$}}
\put(3,0){\makebox(1,1){$v_{1}$}}
\put(4,0){\makebox(1,1){$\cdots$}}
\put(5,0){\makebox(1,1){$v_{k'}$}}
\put(6,0){\makebox(3,1){$\ol{0}\cdots \ol{0}$}}
\multiput(0,0.9)(6,0){2}{\put(0,0){\makebox(3,1){$z$}}}
\end{picture}.
\end{eqnarray*}
Set $l' = 2l-2z$ and $k' =2k-2z$,
hence $T_*$ is a one-row tableau with length $l'$.
Operate the column insertions for semistandard $C$ tableaux
and define
\begin{displaymath}
T^{(0)} := (v_1 \longrightarrow ( \cdots ( v_{k'-1} \longrightarrow ( 
v_{k'} \longrightarrow T_* ) ) \cdots ) ).
\end{displaymath}
It has the form:

\setlength{\unitlength}{5mm}
\begin{picture}(20,4)
\put(5,1.5){\makebox(3,1){$T^{(0)}=$}}
\put(8,1){\line(1,0){3.5}}
\put(8,2){\line(1,0){9}}
\put(8,3){\line(1,0){9}}
\put(8,1){\line(0,1){2}}
\put(11.5,1){\line(0,1){1}} 
\put(12.5,2){\line(0,1){1}} 
\put(17,2){\line(0,1){1}}
\put(12.5,2){\makebox(4.5,1){$i_{m+1} \;\cdots\; i_{l'}$}}
\put(8,1){\makebox(3,1){$\;\;i_1 \cdots i_m$}}
\put(8.5,2){\makebox(3,1){$\;\;j_1 \cdots\cdots j_{k'}$}}
\end{picture}

\noindent
where $m$ is the length of the second row, hence that of the first
row is $l'+k'-m$ ($0 \le m \le k'$).

Next we  bump out  $l'$ letters from
the tableau $T^{(0)}$ by the type $C$ reverse bumping
algorithm in section \ref{subsubsec:invinsc}.
In general, an inverse column insertion starts at a rightmost box in a row.
After an inverse column insertion we obtain a tableau which has the shape
with one box deleted, i.e. the box where we started the reverse bumping is 
removed from the original shape.
We have labeled the boxes by $i_{l'}, i_{l'-1}, \ldots, i_1$ at which we start
the inverse column insertions.
Namely, for the boxes containing $i_{l'}, i_{l'-1}, \ldots, i_1$ in the above
tableau, we do it first for $i_{l'}$ then $i_{l'-1}$ and so on.
Correspondingly, let $w_{1}$ be the first letter that is  bumped out from
the leftmost column and $w_2$ be the second and so on.
Denote by $T^{(i)}$  the resulting tableau when $w_i$ is bumped out
($1 \le i \le l'$).
Note that $w_1 \le w_2 \le \cdots \le w_{l'}$.
Now the $U'_q(C^{(1)}_n)$ crystal elements
$\tilde{b}_1 \in B_{2l}$ and $\tilde{b}_2 \in B_{2k}$ are uniquely specified by
\begin{eqnarray*}
\T(\tilde{b}_2) &=&
\setlength{\unitlength}{5mm}
\begin{picture}(9.5,1.4)(0,0.3)
\multiput(0,0)(0,1){2}{\line(1,0){9}}
\multiput(0,0)(3,0){4}{\line(0,1){1}}
\put(0,0){\makebox(3,1){$0\cdots 0$}}
\put(3,0){\makebox(3,1){$T^{(l')}$}}
\put(6,0){\makebox(3,1){$\ol{0}\cdots \ol{0}$}}
\multiput(0,0.9)(6,0){2}{\put(0,0){\makebox(3,1){$z$}}}
\end{picture},\\
\T(\tilde{b}_1) &=&
\begin{picture}(10.5,2)(0,0.3)
\multiput(0,0)(0,1){2}{\line(1,0){10}}
\multiput(0,0)(3,0){2}{\line(0,1){1}}
\multiput(4.25,0)(1.5,0){2}{\line(0,1){1}}
\multiput(7,0)(3,0){2}{\line(0,1){1}}
\put(0,0){\makebox(3,1){$0\cdots 0$}}
\put(3,0){\makebox(1.25,1){$w_{1}$}}
\put(4.25,0){\makebox(1.5,1){$\cdots$}}
\put(5.75,0){\makebox(1.25,1){$w_{l'}$}}
\put(7,0){\makebox(3,1){$\ol{0}\cdots \ol{0}$}}
\multiput(0,0.9)(7,0){2}{\put(0,0){\makebox(3,1){$z$}}}
\end{picture}.
\end{eqnarray*}
\end{rules}
\hfill
(End of the Rule)
\vskip3ex
We normalize the energy function as $H_{B_l B_k}(b_1 \otimes b_2)=0$
for 
\begin{math}
\mathcal{T}(b_1) =
\setlength{\unitlength}{5mm}
\begin{picture}(3,1.5)(0,0.3)
\multiput(0,0)(0,1){2}{\line(1,0){3}}
\multiput(0,0)(3,0){2}{\line(0,1){1}}
\put(0,0){\makebox(3,1){$1\cdots 1$}}
\put(0,1){\makebox(3,0.5){$\scriptstyle l$}}
\end{picture}
\end{math}
and
\begin{math}
\mathcal{T}(b_2) =
\setlength{\unitlength}{5mm}
\begin{picture}(3,1.5)(0,0.3)
\multiput(0,0)(0,1){2}{\line(1,0){3}}
\multiput(0,0)(3,0){2}{\line(0,1){1}}
\put(0,0){\makebox(3,1){$\ol{1}\cdots \ol{1}$}}
\put(0,1){\makebox(3,0.5){$\scriptstyle k$}}
\end{picture}
\end{math}
irrespective of $l < k$ or $l \ge k$.
Our main result for $A^{(2)}_{2n}$ and 
$D^{(2)}_{n+1}$ is
\begin{theorem}\label{th:main1}
Given $b_1 \ot b_2 \in B_l \ot B_k$, find 
the $U'_q(C^{(1)}_n)$ crystal element
$\tilde{b}_2 \ot \tilde{b}_1 \in B_{2k} \ot B_{2l}$
and $l', k', m$ by Rule \ref{rule:Cx}.
Let $\iota: B_l \ot B_k \stackrel{\sim}{\rightarrow} B_k \ot B_l$ be the isomorphism of
$U'_q(A^{(2)}_{2n})$ (or $U'_q(D^{(2)}_{n+1})$) crystal.
Then $\tilde{b}_2 \ot \tilde{b}_1$ is in the image of the
$B_k \ot B_l$ by the injective map $\omega$ and
we have
\begin{align*}
\iota(b_1\otimes b_2)& = \omega^{-1}(\tilde{b}_2 \otimes \tilde{b}_1),\\
H_{B_l B_k}(b_1 \otimes b_2) &= \min(l',k')- m.
\end{align*}
\end{theorem}
\noindent
By Theorem 3.4 of \cite{HKOT} (the corresponding theorem for the $C^{(1)}_n$ case),
one can immediately obtain this theorem using 
Lemmas \ref{lem:1}, \ref{lem:2} and their corollaries.
\subsection{Examples}
\label{subsec:exCx}
\begin{example}
Let us consider $B_3 \ot B_2 \simeq B_2 \ot B_3$ for $A^{(2)}_4$.
Let $b$ be an element of $B_3$ (resp.~$B_2$). It is depicted by a one-row
tableau $\mathcal{T} (b)$ with length 0, 1, 2 or 3 (resp.~0, 1 or 2).
\begin{displaymath}
\begin{array}{ccccccc}
\setlength{\unitlength}{5mm}
\begin{picture}(3,1)(0,0.3)
\multiput(1,0)(1,0){2}{\line(0,1){1}}
\multiput(1,0)(0,1){2}{\line(1,0){1}}
\put(1,0){\makebox(1,1){$1$}}
\end{picture}
& \otimes & 
\setlength{\unitlength}{5mm}
\begin{picture}(2,1)(0,0.3)
\multiput(0,0)(1,0){3}{\line(0,1){1}}
\multiput(0,0)(0,1){2}{\line(1,0){2}}
\put(0,0){\makebox(1,1){$2$}}
\put(1,0){\makebox(1,1){$2$}}
\end{picture}
& \stackrel{\sim}{\mapsto}  &
\setlength{\unitlength}{5mm}
\begin{picture}(2,1)(0,0.3)
\multiput(0,0)(1,0){3}{\line(0,1){1}}
\multiput(0,0)(0,1){2}{\line(1,0){2}}
\put(0,0){\makebox(1,1){$1$}}
\put(1,0){\makebox(1,1){$1$}}
\end{picture}
& \otimes & 
\setlength{\unitlength}{5mm}
\begin{picture}(3,1)(0,0.3)
\multiput(0,0)(1,0){4}{\line(0,1){1}}
\multiput(0,0)(0,1){2}{\line(1,0){3}}
\put(0,0){\makebox(1,1){$2$}}
\put(1,0){\makebox(1,1){$2$}}
\put(2,0){\makebox(1,1){$\ol{1}$}}
\end{picture}
\\
& & & & & & \\
\setlength{\unitlength}{5mm}
\begin{picture}(3,1)(0,0.3)
\multiput(0.5,0)(1,0){3}{\line(0,1){1}}
\multiput(0.5,0)(0,1){2}{\line(1,0){2}}
\put(0.5,0){\makebox(1,1){$1$}}
\put(1.5,0){\makebox(1,1){$2$}}
\end{picture}
& \otimes & 
\setlength{\unitlength}{5mm}
\begin{picture}(2,1)(0,0.3)
\multiput(0.5,0)(1,0){2}{\line(0,1){1}}
\multiput(0.5,0)(0,1){2}{\line(1,0){1}}
\put(0.5,0){\makebox(1,1){$2$}}
\end{picture}
& \stackrel{\sim}{\mapsto}  &
\setlength{\unitlength}{5mm}
\begin{picture}(2,1)(0,0.3)
\multiput(0.5,0)(1,0){2}{\line(0,1){1}}
\multiput(0.5,0)(0,1){2}{\line(1,0){1}}
\put(0.5,0){\makebox(1,1){$1$}}
\end{picture}
& \otimes & 
\setlength{\unitlength}{5mm}
\begin{picture}(3,1)(0,0.3)
\multiput(0.5,0)(1,0){3}{\line(0,1){1}}
\multiput(0.5,0)(0,1){2}{\line(1,0){2}}
\put(0.5,0){\makebox(1,1){$2$}}
\put(1.5,0){\makebox(1,1){$2$}}
\end{picture}
\\
& & & & & & \\
\setlength{\unitlength}{5mm}
\begin{picture}(3,1)(0,0.3)
\multiput(0,0)(1,0){4}{\line(0,1){1}}
\multiput(0,0)(0,1){2}{\line(1,0){3}}
\put(0,0){\makebox(1,1){$1$}}
\put(1,0){\makebox(1,1){$1$}}
\put(2,0){\makebox(1,1){$2$}}
\end{picture}
& \otimes & 
\setlength{\unitlength}{5mm}
\begin{picture}(2,1)(0,0.3)
\multiput(0,0)(1,0){3}{\line(0,1){1}}
\multiput(0,0)(0,1){2}{\line(1,0){2}}
\put(0,0){\makebox(1,1){$2$}}
\put(1,0){\makebox(1,1){$\ol{1}$}}
\end{picture}
& \stackrel{\sim}{\mapsto}  &
\setlength{\unitlength}{5mm}
\begin{picture}(2,1)(0,0.3)
\multiput(0,0)(1,0){3}{\line(0,1){1}}
\multiput(0,0)(0,1){2}{\line(1,0){2}}
\put(0,0){\makebox(1,1){$1$}}
\put(1,0){\makebox(1,1){$2$}}
\end{picture}
& \otimes & 
\setlength{\unitlength}{5mm}
\begin{picture}(3,1)(0,0.3)
\multiput(1,0)(1,0){2}{\line(0,1){1}}
\multiput(1,0)(0,1){2}{\line(1,0){1}}
\put(1,0){\makebox(1,1){$2$}}
\end{picture}
\end{array}
\end{displaymath}
Here we have picked up three samples.
One can check that they are mapped to each other under the 
isomorphism of the $U_q'(A^{(2)}_4)$ crystals by explicitly writing down
the crystal graphs of $B_3 \ot B_2$ and $B_2 \ot B_3$ .

First we shall show that the use of the
tableau $\mathcal{T} (b)$ given by (\ref{eq:tabtwistax})
is not enough for our purpose, while the less simpler
tableau $\T (\omega (b))$ given by (\ref{eq:tabtwistaxx}) suffices it.
Recall that by neglecting its zero arrows any $U_q'(A^{(2)}_4)$ crystal graph
decomposes into $U_q(C_2)$ crystal graphs.
Thus if $b_1 \ot b_2$ is mapped to $b_2' \ot b_1'$
under the isomorphism of the $U_q'(A^{(2)}_4)$ crystals,
they should also be mapped to each other under an 
isomorphism of $U_q(C_2)$ crystals.
In this example this $U_q(C_2)$ crystal isomorphism 
can be checked in terms of the tableau $\mathcal{T} (b)$ in the following way.
Given $b_1 \ot b_2$ let us construct the product tableau
$\mathcal{T} (b_1) * \mathcal{T} (b_2)$ 
according to the original insertion rule in \cite{B1} where
in particular we have
\begin{math}
(\ol{1} \longrightarrow 
\setlength{\unitlength}{5mm}
\begin{picture}(1,1)(0,0.3)
\multiput(0,0)(1,0){2}{\line(0,1){1}}
\multiput(0,0)(0,1){2}{\line(1,0){1}}
\put(0,0){\makebox(1,1){$1$}}
\end{picture}
) = \emptyset.
\end{math}\footnote{See the first footnote of subsection \ref{subsec:ccis}.}
One can see that
both sides of the
above three mappings then yield a common tableau
\begin{math}
\setlength{\unitlength}{3mm}
\begin{picture}(2,2)(0,0.5)
\multiput(0,0)(1,0){2}{\line(0,1){2}}
\put(2,1){\line(0,1){1}}
\multiput(0,1)(0,1){2}{\line(1,0){2}}
\put(0,0){\line(1,0){1}}
\put(0,1){\makebox(1,1){${\scriptstyle 1}$}}
\put(1,1){\makebox(1,1){${\scriptstyle 2}$}}
\put(0,0){\makebox(1,1){${\scriptstyle 2}$}}
\end{picture}
\end{math}.
This means that they are mapped to each other under isomorphisms of 
$U_q(C_2)$ crystals.
Thus we see that they are satisfying the necessary condition for
the isomorphism of $U_q'(A^{(2)}_4)$ crystals.
However, we also see that this method of 
constructing $\mathcal{T} (b_1) * \mathcal{T} (b_2)$ is not strong enough to determine 
the $U_q'(A^{(2)}_4)$ crystal isomorphism.
Theorem \ref{th:main1} asserts that we are able to determine
the $U_q'(A^{(2)}_4)$ crystal 
isomorphism by means of the tableau $\T (\omega (b_1))* \T (\omega (b_2))$.
Namely the above three mappings are embedded into the following
mappings in
$B_6 \otimes B_4 \simeq B_4 \otimes B_6$ for the $U'_q(C^{(1)}_2)$ 
crystals.
\begin{displaymath}
\begin{array}{ccccccc}
\setlength{\unitlength}{5mm}
\begin{picture}(6,1)(0,0.3)
\multiput(0,0)(1,0){7}{\line(0,1){1}}
\multiput(0,0)(0,1){2}{\line(1,0){6}}
\put(0,0){\makebox(1,1){$0$}}
\put(1,0){\makebox(1,1){$0$}}
\put(2,0){\makebox(1,1){$1$}}
\put(3,0){\makebox(1,1){$1$}}
\put(4,0){\makebox(1,1){$\ol{0}$}}
\put(5,0){\makebox(1,1){$\ol{0}$}}
\end{picture}
& \otimes & 
\setlength{\unitlength}{5mm}
\begin{picture}(4,1)(0,0.3)
\multiput(0,0)(1,0){5}{\line(0,1){1}}
\multiput(0,0)(0,1){2}{\line(1,0){4}}
\put(0,0){\makebox(1,1){$2$}}
\put(1,0){\makebox(1,1){$2$}}
\put(2,0){\makebox(1,1){$2$}}
\put(3,0){\makebox(1,1){$2$}}
\end{picture}
& \stackrel{\sim}{\mapsto}  &
\setlength{\unitlength}{5mm}
\begin{picture}(4,1)(0,0.3)
\multiput(0,0)(1,0){5}{\line(0,1){1}}
\multiput(0,0)(0,1){2}{\line(1,0){4}}
\put(0,0){\makebox(1,1){$1$}}
\put(1,0){\makebox(1,1){$1$}}
\put(2,0){\makebox(1,1){$1$}}
\put(3,0){\makebox(1,1){$1$}}
\end{picture}
& \otimes & 
\setlength{\unitlength}{5mm}
\begin{picture}(6,1)(0,0.3)
\multiput(0,0)(1,0){7}{\line(0,1){1}}
\multiput(0,0)(0,1){2}{\line(1,0){6}}
\put(0,0){\makebox(1,1){$2$}}
\put(1,0){\makebox(1,1){$2$}}
\put(2,0){\makebox(1,1){$2$}}
\put(3,0){\makebox(1,1){$2$}}
\put(4,0){\makebox(1,1){$\ol{1}$}}
\put(5,0){\makebox(1,1){$\ol{1}$}}
\end{picture}
\\
& & & & & & \\
\setlength{\unitlength}{5mm}
\begin{picture}(6,1)(0,0.3)
\multiput(0,0)(1,0){7}{\line(0,1){1}}
\multiput(0,0)(0,1){2}{\line(1,0){6}}
\put(0,0){\makebox(1,1){$0$}}
\put(1,0){\makebox(1,1){$1$}}
\put(2,0){\makebox(1,1){$1$}}
\put(3,0){\makebox(1,1){$2$}}
\put(4,0){\makebox(1,1){$2$}}
\put(5,0){\makebox(1,1){$\ol{0}$}}
\end{picture}
& \otimes & 
\setlength{\unitlength}{5mm}
\begin{picture}(4,1)(0,0.3)
\multiput(0,0)(1,0){5}{\line(0,1){1}}
\multiput(0,0)(0,1){2}{\line(1,0){4}}
\put(0,0){\makebox(1,1){$0$}}
\put(1,0){\makebox(1,1){$2$}}
\put(2,0){\makebox(1,1){$2$}}
\put(3,0){\makebox(1,1){$\ol{0}$}}
\end{picture}
& \stackrel{\sim}{\mapsto}  &
\setlength{\unitlength}{5mm}
\begin{picture}(4,1)(0,0.3)
\multiput(0,0)(1,0){5}{\line(0,1){1}}
\multiput(0,0)(0,1){2}{\line(1,0){4}}
\put(0,0){\makebox(1,1){$0$}}
\put(1,0){\makebox(1,1){$1$}}
\put(2,0){\makebox(1,1){$1$}}
\put(3,0){\makebox(1,1){$\ol{0}$}}
\end{picture}
& \otimes & 
\setlength{\unitlength}{5mm}
\begin{picture}(6,1)(0,0.3)
\multiput(0,0)(1,0){7}{\line(0,1){1}}
\multiput(0,0)(0,1){2}{\line(1,0){6}}
\put(0,0){\makebox(1,1){$0$}}
\put(1,0){\makebox(1,1){$2$}}
\put(2,0){\makebox(1,1){$2$}}
\put(3,0){\makebox(1,1){$2$}}
\put(4,0){\makebox(1,1){$2$}}
\put(5,0){\makebox(1,1){$\ol{0}$}}
\end{picture}
\\
& & & & & & \\
\setlength{\unitlength}{5mm}
\begin{picture}(6,1)(0,0.3)
\multiput(0,0)(1,0){7}{\line(0,1){1}}
\multiput(0,0)(0,1){2}{\line(1,0){6}}
\put(0,0){\makebox(1,1){$1$}}
\put(1,0){\makebox(1,1){$1$}}
\put(2,0){\makebox(1,1){$1$}}
\put(3,0){\makebox(1,1){$1$}}
\put(4,0){\makebox(1,1){$2$}}
\put(5,0){\makebox(1,1){$2$}}
\end{picture}
& \otimes & 
\setlength{\unitlength}{5mm}
\begin{picture}(4,1)(0,0.3)
\multiput(0,0)(1,0){5}{\line(0,1){1}}
\multiput(0,0)(0,1){2}{\line(1,0){4}}
\put(0,0){\makebox(1,1){$2$}}
\put(1,0){\makebox(1,1){$2$}}
\put(2,0){\makebox(1,1){$\ol{1}$}}
\put(3,0){\makebox(1,1){$\ol{1}$}}
\end{picture}
& \stackrel{\sim}{\mapsto}  &
\setlength{\unitlength}{5mm}
\begin{picture}(4,1)(0,0.3)
\multiput(0,0)(1,0){5}{\line(0,1){1}}
\multiput(0,0)(0,1){2}{\line(1,0){4}}
\put(0,0){\makebox(1,1){$1$}}
\put(1,0){\makebox(1,1){$1$}}
\put(2,0){\makebox(1,1){$2$}}
\put(3,0){\makebox(1,1){$2$}}
\end{picture}
& \otimes & 
\setlength{\unitlength}{5mm}
\begin{picture}(6,1)(0,0.3)
\multiput(0,0)(1,0){7}{\line(0,1){1}}
\multiput(0,0)(0,1){2}{\line(1,0){6}}
\put(0,0){\makebox(1,1){$0$}}
\put(1,0){\makebox(1,1){$0$}}
\put(2,0){\makebox(1,1){$2$}}
\put(3,0){\makebox(1,1){$2$}}
\put(4,0){\makebox(1,1){$\ol{0}$}}
\put(5,0){\makebox(1,1){$\ol{0}$}}
\end{picture}
\end{array}
\end{displaymath}
We adopted a rule that
the column insertion
\begin{math}
(\ol{1} \longrightarrow 
\setlength{\unitlength}{5mm}
\begin{picture}(1,1)(0,0.3)
\multiput(0,0)(1,0){2}{\line(0,1){1}}
\multiput(0,0)(0,1){2}{\line(1,0){1}}
\put(0,0){\makebox(1,1){$1$}}
\end{picture}
) 
\end{math}
does not vanish \cite{HKOT}.
Accordingly the both sides of
the first mapping give the tableau
\begin{math}
\setlength{\unitlength}{3mm}
\begin{picture}(6,2)(0,0.5)
\multiput(0,0)(1,0){5}{\line(0,1){2}}
\multiput(5,1)(1,0){2}{\line(0,1){1}}
\multiput(0,1)(0,1){2}{\line(1,0){6}}
\put(0,0){\line(1,0){4}}
\put(0,1){\makebox(1,1){${\scriptstyle 0}$}}
\put(1,1){\makebox(1,1){${\scriptstyle 0}$}}
\put(2,1){\makebox(1,1){${\scriptstyle 1}$}}
\put(3,1){\makebox(1,1){${\scriptstyle 1}$}}
\put(4,1){\makebox(1,1){${\scriptstyle \ol{0}}$}}
\put(5,1){\makebox(1,1){${\scriptstyle \ol{0}}$}}
\put(0,0){\makebox(1,1){${\scriptstyle 2}$}}
\put(1,0){\makebox(1,1){${\scriptstyle 2}$}}
\put(2,0){\makebox(1,1){${\scriptstyle 2}$}}
\put(3,0){\makebox(1,1){${\scriptstyle 2}$}}
\end{picture}
\end{math}.
By deleting a $0,\ol{0}$ pair, 
those of the second one give the tableau
\begin{math}
\setlength{\unitlength}{3mm}
\begin{picture}(4,2)(0,0.5)
\multiput(0,0)(1,0){3}{\line(0,1){2}}
\multiput(3,1)(1,0){2}{\line(0,1){1}}
\multiput(0,1)(0,1){2}{\line(1,0){4}}
\put(0,0){\line(1,0){2}}
\put(0,1){\makebox(1,1){${\scriptstyle 1}$}}
\put(1,1){\makebox(1,1){${\scriptstyle 1}$}}
\put(2,1){\makebox(1,1){${\scriptstyle 2}$}}
\put(3,1){\makebox(1,1){${\scriptstyle 2}$}}
\put(0,0){\makebox(1,1){${\scriptstyle 2}$}}
\put(1,0){\makebox(1,1){${\scriptstyle 2}$}}
\end{picture}
\end{math}.
Those of the third one give the tableau
\begin{math}
\setlength{\unitlength}{3mm}
\begin{picture}(6,2)(0,0.5)
\multiput(0,0)(1,0){5}{\line(0,1){2}}
\multiput(5,1)(1,0){2}{\line(0,1){1}}
\multiput(0,1)(0,1){2}{\line(1,0){6}}
\put(0,0){\line(1,0){4}}
\put(0,1){\makebox(1,1){${\scriptstyle 0}$}}
\put(1,1){\makebox(1,1){${\scriptstyle 0}$}}
\put(2,1){\makebox(1,1){${\scriptstyle 1}$}}
\put(3,1){\makebox(1,1){${\scriptstyle 1}$}}
\put(4,1){\makebox(1,1){${\scriptstyle 2}$}}
\put(5,1){\makebox(1,1){${\scriptstyle 2}$}}
\put(0,0){\makebox(1,1){${\scriptstyle 2}$}}
\put(1,0){\makebox(1,1){${\scriptstyle 2}$}}
\put(2,0){\makebox(1,1){${\scriptstyle \ol{0}}$}}
\put(3,0){\makebox(1,1){${\scriptstyle \ol{0}}$}}
\end{picture}
\end{math}.
They are distinct.
The right hand side is uniquely
determined from the left hand side.

Second let us illustrate in more detail the procedure of Rule \ref{rule:Cx}.
Take the last example.
{}From the left hand side we proceed the column insertions as follows.
\begin{align*}
\ol{1} &\rightarrow 
\setlength{\unitlength}{5mm}
\begin{picture}(6,1)(0,0.3)
\multiput(0,0)(1,0){7}{\line(0,1){1}}
\multiput(0,0)(0,1){2}{\line(1,0){6}}
\put(0,0){\makebox(1,1){$1$}}
\put(1,0){\makebox(1,1){$1$}}
\put(2,0){\makebox(1,1){$1$}}
\put(3,0){\makebox(1,1){$1$}}
\put(4,0){\makebox(1,1){$2$}}
\put(5,0){\makebox(1,1){$2$}}
\end{picture}
\quad = \quad
\setlength{\unitlength}{5mm}
\begin{picture}(6,2)(0,0.8)
\multiput(0,0)(1,0){2}{\line(0,1){2}}
\multiput(2,1)(1,0){5}{\line(0,1){1}}
\multiput(0,1)(0,1){2}{\line(1,0){6}}
\put(0,0){\line(1,0){1}}
\put(0,1){\makebox(1,1){$1$}}
\put(1,1){\makebox(1,1){$1$}}
\put(2,1){\makebox(1,1){$1$}}
\put(3,1){\makebox(1,1){$1$}}
\put(4,1){\makebox(1,1){$2$}}
\put(5,1){\makebox(1,1){$2$}}
\put(0,0){\makebox(1,1){$\ol{1}$}}
\end{picture}
\\
\ol{1} &\rightarrow 
\setlength{\unitlength}{5mm}
\begin{picture}(6,2)(0,0.8)
\multiput(0,0)(1,0){2}{\line(0,1){2}}
\multiput(2,1)(1,0){5}{\line(0,1){1}}
\multiput(0,1)(0,1){2}{\line(1,0){6}}
\put(0,0){\line(1,0){1}}
\put(0,1){\makebox(1,1){$1$}}
\put(1,1){\makebox(1,1){$1$}}
\put(2,1){\makebox(1,1){$1$}}
\put(3,1){\makebox(1,1){$1$}}
\put(4,1){\makebox(1,1){$2$}}
\put(5,1){\makebox(1,1){$2$}}
\put(0,0){\makebox(1,1){$\ol{1}$}}
\end{picture}
\quad = \quad
\setlength{\unitlength}{5mm}
\begin{picture}(6,2)(0,0.8)
\multiput(0,0)(1,0){3}{\line(0,1){2}}
\multiput(3,1)(1,0){4}{\line(0,1){1}}
\multiput(0,1)(0,1){2}{\line(1,0){6}}
\put(0,0){\line(1,0){2}}
\put(0,1){\makebox(1,1){$0$}}
\put(1,1){\makebox(1,1){$1$}}
\put(2,1){\makebox(1,1){$1$}}
\put(3,1){\makebox(1,1){$1$}}
\put(4,1){\makebox(1,1){$2$}}
\put(5,1){\makebox(1,1){$2$}}
\put(0,0){\makebox(1,1){$\ol{1}$}}
\put(1,0){\makebox(1,1){$\ol{0}$}}
\end{picture}
\\
2 &\rightarrow 
\setlength{\unitlength}{5mm}
\begin{picture}(6,2)(0,0.8)
\multiput(0,0)(1,0){3}{\line(0,1){2}}
\multiput(3,1)(1,0){4}{\line(0,1){1}}
\multiput(0,1)(0,1){2}{\line(1,0){6}}
\put(0,0){\line(1,0){2}}
\put(0,1){\makebox(1,1){$0$}}
\put(1,1){\makebox(1,1){$1$}}
\put(2,1){\makebox(1,1){$1$}}
\put(3,1){\makebox(1,1){$1$}}
\put(4,1){\makebox(1,1){$2$}}
\put(5,1){\makebox(1,1){$2$}}
\put(0,0){\makebox(1,1){$\ol{1}$}}
\put(1,0){\makebox(1,1){$\ol{0}$}}
\end{picture}
\quad = \quad
\setlength{\unitlength}{5mm}
\begin{picture}(6,2)(0,0.8)
\multiput(0,0)(1,0){4}{\line(0,1){2}}
\multiput(4,1)(1,0){3}{\line(0,1){1}}
\multiput(0,1)(0,1){2}{\line(1,0){6}}
\put(0,0){\line(1,0){3}}
\put(0,1){\makebox(1,1){$0$}}
\put(1,1){\makebox(1,1){$1$}}
\put(2,1){\makebox(1,1){$1$}}
\put(3,1){\makebox(1,1){$1$}}
\put(4,1){\makebox(1,1){$2$}}
\put(5,1){\makebox(1,1){$2$}}
\put(0,0){\makebox(1,1){$2$}}
\put(1,0){\makebox(1,1){$\ol{1}$}}
\put(2,0){\makebox(1,1){$\ol{0}$}}
\end{picture}
\\
2 &\rightarrow 
\setlength{\unitlength}{5mm}
\begin{picture}(6,2)(0,0.8)
\multiput(0,0)(1,0){4}{\line(0,1){2}}
\multiput(4,1)(1,0){3}{\line(0,1){1}}
\multiput(0,1)(0,1){2}{\line(1,0){6}}
\put(0,0){\line(1,0){3}}
\put(0,1){\makebox(1,1){$0$}}
\put(1,1){\makebox(1,1){$1$}}
\put(2,1){\makebox(1,1){$1$}}
\put(3,1){\makebox(1,1){$1$}}
\put(4,1){\makebox(1,1){$2$}}
\put(5,1){\makebox(1,1){$2$}}
\put(0,0){\makebox(1,1){$2$}}
\put(1,0){\makebox(1,1){$\ol{1}$}}
\put(2,0){\makebox(1,1){$\ol{0}$}}
\end{picture}
\quad = \quad
\setlength{\unitlength}{5mm}
\begin{picture}(6,2)(0,0.8)
\multiput(0,0)(1,0){5}{\line(0,1){2}}
\multiput(5,1)(1,0){2}{\line(0,1){1}}
\multiput(0,1)(0,1){2}{\line(1,0){6}}
\put(0,0){\line(1,0){4}}
\put(0,1){\makebox(1,1){$0$}}
\put(1,1){\makebox(1,1){$0$}}
\put(2,1){\makebox(1,1){$1$}}
\put(3,1){\makebox(1,1){$1$}}
\put(4,1){\makebox(1,1){$2$}}
\put(5,1){\makebox(1,1){$2$}}
\put(0,0){\makebox(1,1){$2$}}
\put(1,0){\makebox(1,1){$2$}}
\put(2,0){\makebox(1,1){$\ol{0}$}}
\put(3,0){\makebox(1,1){$\ol{0}$}}
\end{picture}
\end{align*}
\vskip3ex
\noindent
The reverse bumping procedure goes as follows.
\begin{align*}
T^{(0)} &=
\setlength{\unitlength}{5mm}
\begin{picture}(6,2)(0,0.8)
\multiput(0,0)(1,0){5}{\line(0,1){2}}
\multiput(5,1)(1,0){2}{\line(0,1){1}}
\multiput(0,1)(0,1){2}{\line(1,0){6}}
\put(0,0){\line(1,0){4}}
\put(0,1){\makebox(1,1){$0$}}
\put(1,1){\makebox(1,1){$0$}}
\put(2,1){\makebox(1,1){$1$}}
\put(3,1){\makebox(1,1){$1$}}
\put(4,1){\makebox(1,1){$2$}}
\put(5,1){\makebox(1,1){$2$}}
\put(0,0){\makebox(1,1){$2$}}
\put(1,0){\makebox(1,1){$2$}}
\put(2,0){\makebox(1,1){$\ol{0}$}}
\put(3,0){\makebox(1,1){$\ol{0}$}}
\end{picture}
& \\
T^{(1)} &=
\setlength{\unitlength}{5mm}
\begin{picture}(6,2)(0,0.8)
\multiput(0,0)(1,0){5}{\line(0,1){2}}
\put(5,1){\line(0,1){1}}
\multiput(0,1)(0,1){2}{\line(1,0){5}}
\put(0,0){\line(1,0){4}}
\put(0,1){\makebox(1,1){$0$}}
\put(1,1){\makebox(1,1){$1$}}
\put(2,1){\makebox(1,1){$1$}}
\put(3,1){\makebox(1,1){$2$}}
\put(4,1){\makebox(1,1){$2$}}
\put(0,0){\makebox(1,1){$2$}}
\put(1,0){\makebox(1,1){$2$}}
\put(2,0){\makebox(1,1){$\ol{0}$}}
\put(3,0){\makebox(1,1){$\ol{0}$}}
\end{picture}
&,w_1 = 0 \\
T^{(2)} &=
\setlength{\unitlength}{5mm}
\begin{picture}(6,2)(0,0.8)
\multiput(0,0)(1,0){5}{\line(0,1){2}}
\multiput(0,0)(0,1){3}{\line(1,0){4}}
\put(0,1){\makebox(1,1){$1$}}
\put(1,1){\makebox(1,1){$1$}}
\put(2,1){\makebox(1,1){$2$}}
\put(3,1){\makebox(1,1){$2$}}
\put(0,0){\makebox(1,1){$2$}}
\put(1,0){\makebox(1,1){$2$}}
\put(2,0){\makebox(1,1){$\ol{0}$}}
\put(3,0){\makebox(1,1){$\ol{0}$}}
\end{picture}
&,w_2 = 0 \\
T^{(3)} &=
\setlength{\unitlength}{5mm}
\begin{picture}(6,2)(0,0.8)
\multiput(0,0)(1,0){4}{\line(0,1){2}}
\put(4,1){\line(0,1){1}}
\multiput(0,1)(0,1){2}{\line(1,0){4}}
\put(0,0){\line(1,0){3}}
\put(0,1){\makebox(1,1){$1$}}
\put(1,1){\makebox(1,1){$1$}}
\put(2,1){\makebox(1,1){$2$}}
\put(3,1){\makebox(1,1){$2$}}
\put(0,0){\makebox(1,1){$2$}}
\put(1,0){\makebox(1,1){$\ol{0}$}}
\put(2,0){\makebox(1,1){$\ol{0}$}}
\end{picture}
&,w_3 = 2 \\
T^{(4)} &=
\setlength{\unitlength}{5mm}
\begin{picture}(6,2)(0,0.8)
\multiput(0,0)(1,0){3}{\line(0,1){2}}
\multiput(3,1)(1,0){2}{\line(0,1){1}}
\multiput(0,1)(0,1){2}{\line(1,0){4}}
\put(0,0){\line(1,0){2}}
\put(0,1){\makebox(1,1){$1$}}
\put(1,1){\makebox(1,1){$1$}}
\put(2,1){\makebox(1,1){$2$}}
\put(3,1){\makebox(1,1){$2$}}
\put(0,0){\makebox(1,1){$\ol{0}$}}
\put(1,0){\makebox(1,1){$\ol{0}$}}
\end{picture}
&,w_4 = 2 \\
T^{(5)} &=
\setlength{\unitlength}{5mm}
\begin{picture}(6,2)(0,0.8)
\multiput(0,0)(1,0){2}{\line(0,1){2}}
\multiput(2,1)(1,0){3}{\line(0,1){1}}
\multiput(0,1)(0,1){2}{\line(1,0){4}}
\put(0,0){\line(1,0){1}}
\put(0,1){\makebox(1,1){$1$}}
\put(1,1){\makebox(1,1){$1$}}
\put(2,1){\makebox(1,1){$2$}}
\put(3,1){\makebox(1,1){$2$}}
\put(0,0){\makebox(1,1){$\ol{0}$}}
\end{picture}
&,w_5 = \ol{0} \\
T^{(6)} &=
\setlength{\unitlength}{5mm}
\begin{picture}(6,2)(0,0.3)
\multiput(0,0)(1,0){5}{\line(0,1){1}}
\multiput(0,0)(0,1){2}{\line(1,0){4}}
\put(0,0){\makebox(1,1){$1$}}
\put(1,0){\makebox(1,1){$1$}}
\put(2,0){\makebox(1,1){$2$}}
\put(3,0){\makebox(1,1){$2$}}
\end{picture}
&, w_6 = \ol{0}
\end{align*}
Thus we obtained the right hand side.
We assign $H_{B_3,B_2}=0$ to this element since we have $l'=6, k'=4$ and $m=4$
in this case.
\end{example}
\begin{example}
$B_3 \otimes B_2 \simeq B_2 \otimes B_3$ for $D^{(2)}_3$.
\begin{displaymath}
\begin{array}{ccccccc}
\setlength{\unitlength}{5mm}
\begin{picture}(3,1)(0,0.3)
\multiput(1,0)(1,0){2}{\line(0,1){1}}
\multiput(1,0)(0,1){2}{\line(1,0){1}}
\put(1,0){\makebox(1,1){$1$}}
\end{picture}
& \otimes & 
\setlength{\unitlength}{5mm}
\begin{picture}(2,1)(0,0.3)
\multiput(0,0)(1,0){3}{\line(0,1){1}}
\multiput(0,0)(0,1){2}{\line(1,0){2}}
\put(0,0){\makebox(1,1){$2$}}
\put(1,0){\makebox(1,1){$\circ$}}
\end{picture}
& \stackrel{\sim}{\mapsto}  &
\setlength{\unitlength}{5mm}
\begin{picture}(2,1)(0,0.3)
\multiput(0,0)(1,0){3}{\line(0,1){1}}
\multiput(0,0)(0,1){2}{\line(1,0){2}}
\put(0,0){\makebox(1,1){$1$}}
\put(1,0){\makebox(1,1){$1$}}
\end{picture}
& \otimes & 
\setlength{\unitlength}{5mm}
\begin{picture}(3,1)(0,0.3)
\multiput(0,0)(1,0){4}{\line(0,1){1}}
\multiput(0,0)(0,1){2}{\line(1,0){3}}
\put(0,0){\makebox(1,1){$2$}}
\put(1,0){\makebox(1,1){$\circ$}}
\put(2,0){\makebox(1,1){$\ol{1}$}}
\end{picture}
\\
& & & & & & \\
\setlength{\unitlength}{5mm}
\begin{picture}(3,1)(0,0.3)
\multiput(0.5,0)(1,0){3}{\line(0,1){1}}
\multiput(0.5,0)(0,1){2}{\line(1,0){2}}
\put(0.5,0){\makebox(1,1){$1$}}
\put(1.5,0){\makebox(1,1){$\circ$}}
\end{picture}
& \otimes & 
\setlength{\unitlength}{5mm}
\begin{picture}(2,1)(0,0.3)
\multiput(0.5,0)(1,0){2}{\line(0,1){1}}
\multiput(0.5,0)(0,1){2}{\line(1,0){1}}
\put(0.5,0){\makebox(1,1){$2$}}
\end{picture}
& \stackrel{\sim}{\mapsto}  &
\setlength{\unitlength}{5mm}
\begin{picture}(2,1)(0,0.3)
\multiput(0.5,0)(1,0){2}{\line(0,1){1}}
\multiput(0.5,0)(0,1){2}{\line(1,0){1}}
\put(0.5,0){\makebox(1,1){$1$}}
\end{picture}
& \otimes & 
\setlength{\unitlength}{5mm}
\begin{picture}(3,1)(0,0.3)
\multiput(0.5,0)(1,0){3}{\line(0,1){1}}
\multiput(0.5,0)(0,1){2}{\line(1,0){2}}
\put(0.5,0){\makebox(1,1){$2$}}
\put(1.5,0){\makebox(1,1){$\circ$}}
\end{picture}
\\
& & & & & & \\
\setlength{\unitlength}{5mm}
\begin{picture}(3,1)(0,0.3)
\multiput(0,0)(1,0){4}{\line(0,1){1}}
\multiput(0,0)(0,1){2}{\line(1,0){3}}
\put(0,0){\makebox(1,1){$1$}}
\put(1,0){\makebox(1,1){$1$}}
\put(2,0){\makebox(1,1){$\circ$}}
\end{picture}
& \otimes & 
\setlength{\unitlength}{5mm}
\begin{picture}(2,1)(0,0.3)
\multiput(0,0)(1,0){3}{\line(0,1){1}}
\multiput(0,0)(0,1){2}{\line(1,0){2}}
\put(0,0){\makebox(1,1){$2$}}
\put(1,0){\makebox(1,1){$\ol{1}$}}
\end{picture}
& \stackrel{\sim}{\mapsto}  &
\setlength{\unitlength}{5mm}
\begin{picture}(2,1)(0,0.3)
\multiput(0,0)(1,0){3}{\line(0,1){1}}
\multiput(0,0)(0,1){2}{\line(1,0){2}}
\put(0,0){\makebox(1,1){$1$}}
\put(1,0){\makebox(1,1){$\circ$}}
\end{picture}
& \otimes & 
\setlength{\unitlength}{5mm}
\begin{picture}(3,1)(0,0.3)
\multiput(1,0)(1,0){2}{\line(0,1){1}}
\multiput(1,0)(0,1){2}{\line(1,0){1}}
\put(1,0){\makebox(1,1){$2$}}
\end{picture}
\end{array}
\end{displaymath}
Here we have picked up three samples.
According to the rule of the type $B$ column insertion in \cite{B2} we obtain
\begin{math}
(\ol{1} \longrightarrow 
\setlength{\unitlength}{5mm}
\begin{picture}(1,1)(0,0.3)
\multiput(0,0)(1,0){2}{\line(0,1){1}}
\multiput(0,0)(0,1){2}{\line(1,0){1}}
\put(0,0){\makebox(1,1){$1$}}
\end{picture}
) = \emptyset
\end{math}.
In this rule we find that
both sides of the
above three mappings give a common tableau
\begin{math}
\setlength{\unitlength}{3mm}
\begin{picture}(2,2)(0,0.5)
\multiput(0,0)(1,0){2}{\line(0,1){2}}
\put(2,1){\line(0,1){1}}
\multiput(0,1)(0,1){2}{\line(1,0){2}}
\put(0,0){\line(1,0){1}}
\put(0,1){\makebox(1,1){${\scriptstyle 1}$}}
\put(1,1){\makebox(1,1){${\scriptstyle \circ}$}}
\put(0,0){\makebox(1,1){${\scriptstyle 2}$}}
\end{picture}
\end{math}.

Theorem \ref{th:main1} asserts that we are able to determine
the isomorphism of $U'_q(D^{(2)}_3)$ crystals by means of the tableau $\T (\omega (b))$.
The above three mappings are embedded into the following
mappings in
$B_6 \otimes B_4 \simeq B_4 \otimes B_6$ for the $U'_q(C^{(1)}_2)$ 
crystals.
\begin{displaymath}
\begin{array}{ccccccc}
\setlength{\unitlength}{5mm}
\begin{picture}(6,1)(0,0.3)
\multiput(0,0)(1,0){7}{\line(0,1){1}}
\multiput(0,0)(0,1){2}{\line(1,0){6}}
\put(0,0){\makebox(1,1){$0$}}
\put(1,0){\makebox(1,1){$0$}}
\put(2,0){\makebox(1,1){$1$}}
\put(3,0){\makebox(1,1){$1$}}
\put(4,0){\makebox(1,1){$\ol{0}$}}
\put(5,0){\makebox(1,1){$\ol{0}$}}
\end{picture}
& \otimes & 
\setlength{\unitlength}{5mm}
\begin{picture}(4,1)(0,0.3)
\multiput(0,0)(1,0){5}{\line(0,1){1}}
\multiput(0,0)(0,1){2}{\line(1,0){4}}
\put(0,0){\makebox(1,1){$2$}}
\put(1,0){\makebox(1,1){$2$}}
\put(2,0){\makebox(1,1){$2$}}
\put(3,0){\makebox(1,1){$\ol{2}$}}
\end{picture}
& \stackrel{\sim}{\mapsto}  &
\setlength{\unitlength}{5mm}
\begin{picture}(4,1)(0,0.3)
\multiput(0,0)(1,0){5}{\line(0,1){1}}
\multiput(0,0)(0,1){2}{\line(1,0){4}}
\put(0,0){\makebox(1,1){$1$}}
\put(1,0){\makebox(1,1){$1$}}
\put(2,0){\makebox(1,1){$1$}}
\put(3,0){\makebox(1,1){$1$}}
\end{picture}
& \otimes & 
\setlength{\unitlength}{5mm}
\begin{picture}(6,1)(0,0.3)
\multiput(0,0)(1,0){7}{\line(0,1){1}}
\multiput(0,0)(0,1){2}{\line(1,0){6}}
\put(0,0){\makebox(1,1){$2$}}
\put(1,0){\makebox(1,1){$2$}}
\put(2,0){\makebox(1,1){$2$}}
\put(3,0){\makebox(1,1){$\ol{2}$}}
\put(4,0){\makebox(1,1){$\ol{1}$}}
\put(5,0){\makebox(1,1){$\ol{1}$}}
\end{picture}
\\
& & & & & & \\
\setlength{\unitlength}{5mm}
\begin{picture}(6,1)(0,0.3)
\multiput(0,0)(1,0){7}{\line(0,1){1}}
\multiput(0,0)(0,1){2}{\line(1,0){6}}
\put(0,0){\makebox(1,1){$0$}}
\put(1,0){\makebox(1,1){$1$}}
\put(2,0){\makebox(1,1){$1$}}
\put(3,0){\makebox(1,1){$2$}}
\put(4,0){\makebox(1,1){$\ol{2}$}}
\put(5,0){\makebox(1,1){$\ol{0}$}}
\end{picture}
& \otimes & 
\setlength{\unitlength}{5mm}
\begin{picture}(4,1)(0,0.3)
\multiput(0,0)(1,0){5}{\line(0,1){1}}
\multiput(0,0)(0,1){2}{\line(1,0){4}}
\put(0,0){\makebox(1,1){$0$}}
\put(1,0){\makebox(1,1){$2$}}
\put(2,0){\makebox(1,1){$2$}}
\put(3,0){\makebox(1,1){$\ol{0}$}}
\end{picture}
& \stackrel{\sim}{\mapsto}  &
\setlength{\unitlength}{5mm}
\begin{picture}(4,1)(0,0.3)
\multiput(0,0)(1,0){5}{\line(0,1){1}}
\multiput(0,0)(0,1){2}{\line(1,0){4}}
\put(0,0){\makebox(1,1){$0$}}
\put(1,0){\makebox(1,1){$1$}}
\put(2,0){\makebox(1,1){$1$}}
\put(3,0){\makebox(1,1){$\ol{0}$}}
\end{picture}
& \otimes & 
\setlength{\unitlength}{5mm}
\begin{picture}(6,1)(0,0.3)
\multiput(0,0)(1,0){7}{\line(0,1){1}}
\multiput(0,0)(0,1){2}{\line(1,0){6}}
\put(0,0){\makebox(1,1){$0$}}
\put(1,0){\makebox(1,1){$2$}}
\put(2,0){\makebox(1,1){$2$}}
\put(3,0){\makebox(1,1){$2$}}
\put(4,0){\makebox(1,1){$\ol{2}$}}
\put(5,0){\makebox(1,1){$\ol{0}$}}
\end{picture}
\\
& & & & & & \\
\setlength{\unitlength}{5mm}
\begin{picture}(6,1)(0,0.3)
\multiput(0,0)(1,0){7}{\line(0,1){1}}
\multiput(0,0)(0,1){2}{\line(1,0){6}}
\put(0,0){\makebox(1,1){$1$}}
\put(1,0){\makebox(1,1){$1$}}
\put(2,0){\makebox(1,1){$1$}}
\put(3,0){\makebox(1,1){$1$}}
\put(4,0){\makebox(1,1){$2$}}
\put(5,0){\makebox(1,1){$\ol{2}$}}
\end{picture}
& \otimes & 
\setlength{\unitlength}{5mm}
\begin{picture}(4,1)(0,0.3)
\multiput(0,0)(1,0){5}{\line(0,1){1}}
\multiput(0,0)(0,1){2}{\line(1,0){4}}
\put(0,0){\makebox(1,1){$2$}}
\put(1,0){\makebox(1,1){$2$}}
\put(2,0){\makebox(1,1){$\ol{1}$}}
\put(3,0){\makebox(1,1){$\ol{1}$}}
\end{picture}
& \stackrel{\sim}{\mapsto}  &
\setlength{\unitlength}{5mm}
\begin{picture}(4,1)(0,0.3)
\multiput(0,0)(1,0){5}{\line(0,1){1}}
\multiput(0,0)(0,1){2}{\line(1,0){4}}
\put(0,0){\makebox(1,1){$1$}}
\put(1,0){\makebox(1,1){$1$}}
\put(2,0){\makebox(1,1){$2$}}
\put(3,0){\makebox(1,1){$\ol{2}$}}
\end{picture}
& \otimes & 
\setlength{\unitlength}{5mm}
\begin{picture}(6,1)(0,0.3)
\multiput(0,0)(1,0){7}{\line(0,1){1}}
\multiput(0,0)(0,1){2}{\line(1,0){6}}
\put(0,0){\makebox(1,1){$0$}}
\put(1,0){\makebox(1,1){$0$}}
\put(2,0){\makebox(1,1){$2$}}
\put(3,0){\makebox(1,1){$2$}}
\put(4,0){\makebox(1,1){$\ol{0}$}}
\put(5,0){\makebox(1,1){$\ol{0}$}}
\end{picture}
\end{array}
\end{displaymath}
The both sides of
the first mapping give the tableau
\begin{math}
\setlength{\unitlength}{3mm}
\begin{picture}(6,2)(0,0.5)
\multiput(0,0)(1,0){5}{\line(0,1){2}}
\multiput(5,1)(1,0){2}{\line(0,1){1}}
\multiput(0,1)(0,1){2}{\line(1,0){6}}
\put(0,0){\line(1,0){4}}
\put(0,1){\makebox(1,1){${\scriptstyle 0}$}}
\put(1,1){\makebox(1,1){${\scriptstyle 0}$}}
\put(2,1){\makebox(1,1){${\scriptstyle 1}$}}
\put(3,1){\makebox(1,1){${\scriptstyle 1}$}}
\put(4,1){\makebox(1,1){${\scriptstyle \ol{0}}$}}
\put(5,1){\makebox(1,1){${\scriptstyle \ol{0}}$}}
\put(0,0){\makebox(1,1){${\scriptstyle 2}$}}
\put(1,0){\makebox(1,1){${\scriptstyle 2}$}}
\put(2,0){\makebox(1,1){${\scriptstyle 2}$}}
\put(3,0){\makebox(1,1){${\scriptstyle \ol{2}}$}}
\end{picture}
\end{math}.
By deleting a $0,\ol{0}$ pair, 
those of the second one give the tableau
\begin{math}
\setlength{\unitlength}{3mm}
\begin{picture}(4,2)(0,0.5)
\multiput(0,0)(1,0){3}{\line(0,1){2}}
\multiput(3,1)(1,0){2}{\line(0,1){1}}
\multiput(0,1)(0,1){2}{\line(1,0){4}}
\put(0,0){\line(1,0){2}}
\put(0,1){\makebox(1,1){${\scriptstyle 1}$}}
\put(1,1){\makebox(1,1){${\scriptstyle 1}$}}
\put(2,1){\makebox(1,1){${\scriptstyle 2}$}}
\put(3,1){\makebox(1,1){${\scriptstyle \ol{2}}$}}
\put(0,0){\makebox(1,1){${\scriptstyle 2}$}}
\put(1,0){\makebox(1,1){${\scriptstyle 2}$}}
\end{picture}
\end{math}.
Those of the third one give the tableau
\begin{math}
\setlength{\unitlength}{3mm}
\begin{picture}(6,2)(0,0.5)
\multiput(0,0)(1,0){5}{\line(0,1){2}}
\multiput(5,1)(1,0){2}{\line(0,1){1}}
\multiput(0,1)(0,1){2}{\line(1,0){6}}
\put(0,0){\line(1,0){4}}
\put(0,1){\makebox(1,1){${\scriptstyle 0}$}}
\put(1,1){\makebox(1,1){${\scriptstyle 0}$}}
\put(2,1){\makebox(1,1){${\scriptstyle 1}$}}
\put(3,1){\makebox(1,1){${\scriptstyle 1}$}}
\put(4,1){\makebox(1,1){${\scriptstyle 2}$}}
\put(5,1){\makebox(1,1){${\scriptstyle \ol{2}}$}}
\put(0,0){\makebox(1,1){${\scriptstyle 2}$}}
\put(1,0){\makebox(1,1){${\scriptstyle 2}$}}
\put(2,0){\makebox(1,1){${\scriptstyle \ol{0}}$}}
\put(3,0){\makebox(1,1){${\scriptstyle \ol{0}}$}}
\end{picture}
\end{math}.
They are distinct.
The right hand side is uniquely
determined from the left hand side.
\end{example}
\section{\mathversion{bold}$U_q'(B_n^{(1)})$ and $U_q'(D_n^{(1)})$ crystal cases}
\label{sec:bd}
\subsection{\mathversion{bold}Definitions : $U_q'(B_n^{(1)})$ crystal case}
\label{subsec:typeB}
Given a positive integer $l$,
let us denote by $B_{l}$ the $U_q'(B_n^{(1)})$ crystal
defined in \cite{KKM}.
As a set $B_{l}$ reads
$$
B_{l} = \left\{(
x_1,\ldots, x_n,x_\circ,\overline{x}_n,\ldots,\overline{x}_1) \Biggm|
x_\circ=\mbox{$0$ or $1$}, x_i, \overline{x}_i \in \Z_{\ge 0},
x_\circ+\sum_{i=1}^n(x_i + \overline{x}_i) = l \right\}.
$$
For its crystal structure see \cite{KKM}.
$B_{l}$ is isomorphic to
$B(l \La_1)$ as a $U_q(B_n)$ crystal.
We depict the element
$b= (x_1, \ldots, x_n, x_\circ,\overline{x}_n,\ldots,\overline{x}_1) \in
B_{l}$ by the tableau
\begin{displaymath}
\mathcal{T}(b)=\overbrace{\fbox{$\vphantom{\ol{1}} 1 \cd 1$}}^{x_1}\!
\fbox{$\vphantom{\ol{1}}\cd$}\!
\overbrace{\fbox{$\vphantom{\ol{1}}n \cd n$}}^{x_n}\!
\overbrace{\fbox{$\vphantom{\ol{1}}\hphantom{1}\circ\hphantom{1}$}}^{x_\circ}\!
\overbrace{\fbox{$\vphantom{\ol{1}}\ol{n} \cd \ol{n}$}}^{\ol{x}_n}\!
\fbox{$\vphantom{\ol{1}}\cd$}\!
\overbrace{\fbox{$\ol{1} \cd \ol{1}$}}^{\ol{x}_1}.
\end{displaymath}
The length of this one-row tableau is equal to $l$, namely
$x_\circ+\sum_{i=1}^n(x_i + \overline{x}_i) =l$.
\subsection{\mathversion{bold}Definitions : $U_q'(D_n^{(1)})$ crystal case}
\label{subsec:typeD}
Given a positive integer $l$,
let us denote by $B_{l}$ the $U_q'(D_n^{(1)})$ crystal
defined in \cite{KKM}.
As a set $B_{l}$ reads
$$
B_{l} = \left\{(
x_1,\ldots, x_n,\overline{x}_n,\ldots,\overline{x}_1) \Biggm|
\mbox{$x_n=0$ or $\overline{x}_n=0$}, x_i, \overline{x}_i \in \Z_{\ge 0},
\sum_{i=1}^n(x_i + \overline{x}_i) = l \right\}.
$$
For its crystal structure see \cite{KKM}.
$B_{l}$ is isomorphic to
$B(l \La_1)$ as a $U_q(D_n)$ crystal.
We depict the element
$b= (x_1, \ldots, x_n,\overline{x}_n,\ldots,\overline{x}_1) \in
B_{l}$ by the tableau
\begin{displaymath}
{\mathcal T} (b)=\overbrace{\fbox{$\vphantom{\ol{1}} 1 \cd 1$}}^{x_1}\!
\fbox{$\vphantom{\ol{1}}\cd$}\!
\overbrace{\fbox{$\vphantom{\ol{1}}n \cd n$}}^{x_n}\!
\overbrace{\fbox{$\vphantom{\ol{1}}\ol{n} \cd \ol{n}$}}^{\ol{x}_n}\!
\fbox{$\vphantom{\ol{1}}\cd$}\!
\overbrace{\fbox{$\ol{1} \cd \ol{1}$}}^{\ol{x}_1}.
\end{displaymath}
The length of this one-row tableau is equal to $l$, namely
$\sum_{i=1}^n(x_i + \overline{x}_i) =l$.
\subsection{\mathversion{bold}Column insertion and inverse insertion for
$B_n$}
\label{subsec:cib}
Set an alphabet $\mathcal{X}=\mathcal{A} \sqcup \{\circ\} \sqcup \bar{\mathcal{A}},\,
\mathcal{A}=\{ 1,\dots,n\}$ and
$\bar{\mathcal{A}}=\{\overline{1},\dots,\overline{n}\}$,
with the total order
$1 < 2 < \dots < n < \circ < \overline{n} < \dots < \overline{2} < \overline{1}$.
\subsubsection{Semistandard $B$ tableaux}
\label{subsubsec:ssbt}
Let us consider a {\em semistandard $B$ tableaux}
made by the letters from this alphabet.
We follow \cite{KN} for its definition.
We present the definition here, but restrict ourselves to
special cases that are sufficient for our purpose.
Namely we consider only
those tableaux that have no more than two rows
in their shapes. 
Thus they have the forms as in (\ref{eq:pictureofsst}),
with the letters inside the boxes are now chosen from the alphabet 
given as above.
The letters obey the conditions
(\ref{eq:notdecrease}) and the absence of the $(x,x)$-configuration
(\ref{eq:absenceofxxconf}) where we now assume $1 \leq x <n$.
They also obey the following conditions:
\begin{equation}
\alpha_a < \beta_a \quad \mbox{or} \quad (\alpha_a,\beta_a) = (\circ,\circ),
\end{equation}
\begin{equation}
\label{eq:absenceofoneonebar}
(\alpha_a,\beta_a) \ne (1,\ol{1}),
\end{equation}
\begin{equation}
\label{eq:absenceofnnconf}
(\alpha_a,\beta_{a+1}) \ne (n,\ol{n}),(n,\circ),(\circ,\circ),(\circ,\ol{n}).
\end{equation}
The last conditions (\ref{eq:absenceofnnconf}) are referred to as the
absence of the $(n,n)$-configurations.
\subsubsection{\mathversion{bold}Column insertion for 
$B_n$ \cite{B2}}
\label{subsubsec:insb}
We give a partial list of patterns of column insertions on the 
semistandard $B$ tableaux that are sufficient for our purpose.
For the alphabet $\mathcal{X}$, we follow the convention
that Greek letters $ \alpha, \beta, \ldots $ belong to
$\mathcal{X}$ while Latin
letters $x,y,\ldots$ (resp. $\overline{x},\overline{y},\ldots$) belong to
$\mathcal{A}$ (resp. $\bar{\mathcal{A}}$).
For the interpretation of the pictorial equations in the list,
see the remarks in Section \ref{subsubsec:insc}.
Note that this list does not exhaust all cases.
It does not contain, for instance, the insertion
\begin{math}
\setlength{\unitlength}{3mm}
\begin{picture}(3,2)(0,0.3)
\multiput(0,0)(1,0){2}{\line(0,1){1}}
\multiput(0,0)(0,1){2}{\line(1,0){1}}
\put(0,0){\makebox(1,1){${\scriptstyle \circ}$}}
\put(1,0){\makebox(1,1){${\scriptstyle \rightarrow}$}}
\multiput(2,0)(1,0){2}{\line(0,1){2}}
\multiput(2,0)(0,1){3}{\line(1,0){1}}
\put(2,0){\makebox(1,1){${\scriptstyle \circ}$}}
\put(2,1){\makebox(1,1){${\scriptstyle \circ}$}}
\end{picture}
\end{math}.
As far as this case is concerned, we see that
in Rule \ref{rule:typeB} neither $\mathcal{T} (b_1)$ nor $\mathcal{T} (b_2)$
has more than one $\circ$'s.
Thus we do not encounter a situation where more than two $\circ$'s
appear in the procedure. 
(See Proposition \ref{pr:atmosttworows}.)

\noindent
\ToBOX{A0}{\alpha}%
	\raisebox{1.25mm}{,}

\noindent
\ToDOMINO{A1}{\alpha}{\beta}{\alpha}{\beta}%
	\raisebox{4mm}{if $\alpha < \beta$ or $(\alpha,\beta)=(\circ,\circ)$,}

\noindent
\ToYOKODOMINO{B0}{\beta}{\alpha}{\beta}{\alpha}%
	\raisebox{1.25mm}{if $\alpha \le \beta$ and $(\alpha,\beta) \ne (\circ,\circ)$,}

\noindent
\ToHOOKnn{B1}{\alpha}{\gamma}{\beta}{\gamma}{\alpha}{\beta}%
	\raisebox{4mm}{if $\alpha < \beta \leq \gamma$ and $(\alpha,\gamma) \ne (x,\overline{x})$
	and $(\beta,\gamma) \ne (\circ,\circ)$,}

\noindent
\ToHOOKnn{B2}{\beta}{\gamma}{\alpha}{\beta}{\alpha}{\gamma}
	\raisebox{4mm}{if $\alpha \leq \beta < \gamma$ and $(\alpha,\gamma) \ne (x,\overline{x})$
	and $(\alpha,\beta) \ne (\circ,\circ)$,}

\noindent
\ToHOOKnn{B3}{\circ}{\overline{x}}{\circ}{\overline{x}}{\circ}{\circ}
	\raisebox{4mm}{,}

\noindent
\ToHOOKnn{B4}{\circ}{\circ}{x}{\circ}{x}{\circ}
	\raisebox{4mm}{,}

\noindent
\ToHOOKll{B5}{x}{\overline{x}}{\beta}{\overline{x\!-\!1}}{x\!-\!1}{\beta}%
	\raisebox{4mm}{if $x \le \beta \le \overline{x}$ and $x\ne 1$,}

\noindent
\ToHOOKln{B6}{\beta}{\overline{x}}{x}{\beta}{x\!+\!1}{\overline{x\!+\!1}}%
	\raisebox{4mm}{if $x < \beta < \overline{x}$ and $x\ne n$,}

\noindent
\ToHOOKnn{B7}{\circ}{\overline{n}}{n}{\overline{n}}{n}{\circ}%
	\raisebox{4mm}{.}

\noindent
\subsubsection{Column insertion and $U_q (B_n)$ crystal morphism}
To illustrate Proposition \ref{pr:morphgen}
let us check a morphism of the $U_q(B_3)$ crystal
$B(\Lambda_2) \ot B(\Lambda_1)$ by taking an example.
Let $\psi$ be the map that is similarly defined
as in Section \ref{subsubsec:insc} for type $C$ case.
\begin{example}
\label{ex:morB}
\begin{displaymath}
\begin{CD}
\setlength{\unitlength}{5mm}
\begin{picture}(3,2)(0,0.3)
\multiput(0,0)(1,0){2}{\line(0,1){2}}
\multiput(0,0)(0,1){3}{\line(1,0){1}}
\put(0,1){\makebox(1,1){$\circ$}}
\put(0,0){\makebox(1,1){$\ol{3}$}}
\put(1,0.5){\makebox(1,1){$\otimes$}}
\put(2,0.5){\makebox(1,1){$\circ$}}
\multiput(2,0.5)(1,0){2}{\line(0,1){1}}
\multiput(2,0.5)(0,1){2}{\line(1,0){1}}
\end{picture}
@>\text{$\etd_3$}>>
\setlength{\unitlength}{5mm}
\begin{picture}(3,2)(0,0.3)
\multiput(0,0)(1,0){2}{\line(0,1){2}}
\multiput(0,0)(0,1){3}{\line(1,0){1}}
\put(0,1){\makebox(1,1){$\circ$}}
\put(0,0){\makebox(1,1){$\ol{3}$}}
\put(1,0.5){\makebox(1,1){$\otimes$}}
\put(2,0.5){\makebox(1,1){$3$}}
\multiput(2,0.5)(1,0){2}{\line(0,1){1}}
\multiput(2,0.5)(0,1){2}{\line(1,0){1}}
\end{picture}
@>\text{$\etd_3$}>>
\setlength{\unitlength}{5mm}
\begin{picture}(3,2)(0,0.3)
\multiput(0,0)(1,0){2}{\line(0,1){2}}
\multiput(0,0)(0,1){3}{\line(1,0){1}}
\put(0,1){\makebox(1,1){$\circ$}}
\put(0,0){\makebox(1,1){$\circ$}}
\put(1,0.5){\makebox(1,1){$\otimes$}}
\put(2,0.5){\makebox(1,1){$3$}}
\multiput(2,0.5)(1,0){2}{\line(0,1){1}}
\multiput(2,0.5)(0,1){2}{\line(1,0){1}}
\end{picture}
@>\text{$\etd_3$}>>
\setlength{\unitlength}{5mm}
\begin{picture}(3,2)(0,0.3)
\multiput(0,0)(1,0){2}{\line(0,1){2}}
\multiput(0,0)(0,1){3}{\line(1,0){1}}
\put(0,1){\makebox(1,1){$3$}}
\put(0,0){\makebox(1,1){$\circ$}}
\put(1,0.5){\makebox(1,1){$\otimes$}}
\put(2,0.5){\makebox(1,1){$3$}}
\multiput(2,0.5)(1,0){2}{\line(0,1){1}}
\multiput(2,0.5)(0,1){2}{\line(1,0){1}}
\end{picture}
\\
@VV\text{$\psi$}V @VV\text{$\psi$}V @VV\text{$\psi$}V @VV\text{$\psi$}V \\
\setlength{\unitlength}{5mm}
\begin{picture}(2,2)(0,0.3)
\multiput(0,0)(1,0){2}{\line(0,1){2}}
\put(0,0){\line(1,0){1}}
\multiput(0,1)(0,1){2}{\line(1,0){2}}
\put(2,1){\line(0,1){1}}
\put(0,1){\makebox(1,1){$\circ$}}\put(1,1){\makebox(1,1){$\ol{3}$}}
\put(0,0){\makebox(1,1){$\circ$}}
\end{picture}
@>\text{$\etd_3$}>>
\setlength{\unitlength}{5mm}
\begin{picture}(2,2)(0,0.3)
\multiput(0,0)(1,0){2}{\line(0,1){2}}
\put(0,0){\line(1,0){1}}
\multiput(0,1)(0,1){2}{\line(1,0){2}}
\put(2,1){\line(0,1){1}}
\put(0,1){\makebox(1,1){$3$}}\put(1,1){\makebox(1,1){$\ol{3}$}}
\put(0,0){\makebox(1,1){$\circ$}}
\end{picture}
@>\text{$\etd_3$}>>
\setlength{\unitlength}{5mm}
\begin{picture}(2,2)(0,0.3)
\multiput(0,0)(1,0){2}{\line(0,1){2}}
\put(0,0){\line(1,0){1}}
\multiput(0,1)(0,1){2}{\line(1,0){2}}
\put(2,1){\line(0,1){1}}
\put(0,1){\makebox(1,1){$3$}}\put(1,1){\makebox(1,1){$\circ$}}
\put(0,0){\makebox(1,1){$\circ$}}
\end{picture}
@>\text{$\etd_3$}>>
\setlength{\unitlength}{5mm}
\begin{picture}(2,2)(0,0.3)
\multiput(0,0)(1,0){2}{\line(0,1){2}}
\put(0,0){\line(1,0){1}}
\multiput(0,1)(0,1){2}{\line(1,0){2}}
\put(2,1){\line(0,1){1}}
\put(0,1){\makebox(1,1){$3$}}\put(1,1){\makebox(1,1){$3$}}
\put(0,0){\makebox(1,1){$\circ$}}
\end{picture}
\end{CD}
\end{displaymath}
\vskip3ex
\noindent
Here the $\psi$'s are given by Case B3, B7, B4 and B2 column insertions,
respectively from left to right.
\end{example}
\subsubsection{\mathversion{bold}Inverse insertion for 
$B_n$ \cite{B2}}
\label{subsubsec:invinsb}
We give a list of inverse column insertions on 
semistandard $B$ tableaux that are sufficient for our purpose.
For the interpretation of the pictorial equations in the list,
see the remarks in Section
\ref{subsubsec:invinsc}.
\par
\noindent
\FromYOKODOMINO{C0}{\beta}{\alpha}{\beta}{\alpha}%
	\raisebox{1.25mm}{if $\alpha \le \beta$ and $(\alpha,\beta) \ne (\circ,\circ)$,}

\noindent
\FromHOOKnn{C1}{\gamma}{\alpha}{\beta}{\alpha}{\gamma}{\beta}%
	\raisebox{4mm}{if $\alpha < \beta \leq \gamma$ and $(\alpha,\gamma) \ne (x,\overline{x})$
	and $(\beta,\gamma) \ne (\circ,\circ)$,}

\noindent
\FromHOOKnn{C2}{\beta}{\alpha}{\gamma}{\beta}{\gamma}{\alpha}%
	\raisebox{4mm}{if $\alpha \leq \beta < \gamma$ and $(\alpha,\gamma) \ne (x,\overline{x})$
	and $(\alpha,\beta) \ne (\circ,\circ)$,}

\noindent
\FromHOOKnn{C3}{\overline{x}}{\circ}{\circ}{\circ}{\overline{x}}{\circ}
	\raisebox{4mm}{,}

\noindent
\FromHOOKnn{C4}{\circ}{x}{\circ}{\circ}{\circ}{x}
	\raisebox{4mm}{,}

\noindent
\FromHOOKnl{C5}{\overline{x}}{x}{\beta}{x\!+\!1}{\overline{x\!+\!1}}{\beta}%
	\raisebox{4mm}{if $x < \beta < \overline{x}$ and $x\ne n$,}

\noindent
\FromHOOKll{C6}{\beta}{x}{\overline{x}}{\beta}{\overline{x\!-\!1}}{x\!-\!1}%
	\raisebox{4mm}{if $x \le \beta \le \overline{x}$ and $x\ne 1$,}

\noindent
\FromHOOKnn{C7}{\overline{n}}{n}{\circ}{\circ}{\overline{n}}{n}%
	\raisebox{4mm}{.}
\subsubsection{\mathversion{bold} Column bumping lemma for 
$B_n$}
\label{subsubsec:cblB}
The aim of this subsection is to give a simple result on
successive insertions of two letters into a 
tableau (Corollary \ref{cor:bcblxxx}).
This result will be
used in the proof of the main theorem (Theorem \ref{th:main3}).
This corollary follows from Lemma \ref{lem:bcblxx}.
This lemma is (a special case of) the {\em column bumping lemma},
whose claim is almost the same as that for the original lemma
for the usual tableaux (\cite{F}, Exercise 3 of
Appendix A).

We restrict ourselves to the situation where by column insertions there only appear 
semistandard $B$ tableaux with at most two rows.
We consider a column insertion of a letter $\alpha$ into a tableau $T$.
We insert the $\alpha$ into the leftmost column of $T$.
According to the rules, the $\alpha$ is set in the column, and bump a
letter if possible.
The bumped letter is then inserted in the right column.
The procedure continues until we come to Case A0 or A1.

When a letter is inserted in a tableau, we can define a {\em bumping route}.
It is a collection of boxes in the new tableau that has those letters 
set by the insertion.
In each column there is at most one such box.
Thus we regard the bumping route as a path that goes from the left to 
the right.
In the classification of the column insertions, we regard 
that the inserted letter is set in the first row in Cases A0, B0, B2, B4, B6 and
B7, and that it is set in the second row in the other cases.
\begin{example}
Here we give an example of column insertion and its resulting
bumping route in a $B_3$ tableau.
\begin{displaymath}
\setlength{\unitlength}{5mm}
\begin{picture}(4,4)(0,1)
\multiput(0,2)(1,0){5}{\line(0,1){2}}
\multiput(0,2)(0,1){3}{\line(1,0){4}}
\put(0,3){\makebox(1,1){$1$}}
\put(1,3){\makebox(1,1){$2$}}
\put(2,3){\makebox(1,1){$\circ$}}
\put(3,3){\makebox(1,1){$\ol{3}$}}
\put(0,2){\makebox(1,1){$3$}}
\put(1,2){\makebox(1,1){$3$}}
\put(2,2){\makebox(1,1){$\ol{3}$}}
\put(3,2){\makebox(1,1){$\ol{2}$}}
\put(0,1){\makebox(1,1){$\uparrow$}}
\put(0,0){\makebox(1,1){$2$}}
\end{picture}
\quad
\Rightarrow
\quad
\begin{picture}(4,4)(0,1)
\multiput(0,2)(1,0){5}{\line(0,1){2}}
\multiput(0,2)(0,1){3}{\line(1,0){4}}
\put(0,3){\makebox(1,1){$1$}}
\put(1,3){\makebox(1,1){$2$}}
\put(2,3){\makebox(1,1){$\circ$}}
\put(3,3){\makebox(1,1){$\ol{3}$}}
\put(0,2){\makebox(1,1){$2$}}
\put(0,2){\makebox(1,1){$\bigcirc$}}
\put(1,2){\makebox(1,1){$3$}}
\put(2,2){\makebox(1,1){$\ol{3}$}}
\put(3,2){\makebox(1,1){$\ol{2}$}}
\put(1,1){\makebox(1,1){$\uparrow$}}
\put(1,0){\makebox(1,1){$3$}}
\end{picture}
\quad
\Rightarrow
\quad
\begin{picture}(4,4)(0,1)
\multiput(0,2)(1,0){5}{\line(0,1){2}}
\multiput(0,2)(0,1){3}{\line(1,0){4}}
\put(0,3){\makebox(1,1){$1$}}
\put(1,3){\makebox(1,1){$2$}}
\put(2,3){\makebox(1,1){$\circ$}}
\put(3,3){\makebox(1,1){$\ol{3}$}}
\put(0,2){\makebox(1,1){$2$}}
\put(0,2){\makebox(1,1){$\bigcirc$}}
\put(1,2){\makebox(1,1){$3$}}
\put(1,2){\makebox(1,1){$\bigcirc$}}
\put(2,2){\makebox(1,1){$\ol{3}$}}
\put(3,2){\makebox(1,1){$\ol{2}$}}
\put(2,1){\makebox(1,1){$\uparrow$}}
\put(2,0){\makebox(1,1){$3$}}
\end{picture}
\quad
\Rightarrow
\quad
\begin{picture}(4,4)(0,1)
\multiput(0,2)(1,0){5}{\line(0,1){2}}
\multiput(0,2)(0,1){3}{\line(1,0){4}}
\put(0,3){\makebox(1,1){$1$}}
\put(1,3){\makebox(1,1){$2$}}
\put(2,3){\makebox(1,1){$3$}}
\put(2,3){\makebox(1,1){$\bigcirc$}}
\put(3,3){\makebox(1,1){$\ol{3}$}}
\put(0,2){\makebox(1,1){$2$}}
\put(0,2){\makebox(1,1){$\bigcirc$}}
\put(1,2){\makebox(1,1){$3$}}
\put(1,2){\makebox(1,1){$\bigcirc$}}
\put(2,2){\makebox(1,1){$\circ$}}
\put(3,2){\makebox(1,1){$\ol{2}$}}
\put(3,1){\makebox(1,1){$\uparrow$}}
\put(3,0){\makebox(1,1){$\ol{3}$}}
\end{picture}
\end{displaymath}
\begin{displaymath}
\quad
\Rightarrow
\quad
\setlength{\unitlength}{5mm}
\begin{picture}(5,4)(0,1)
\multiput(0,2)(1,0){5}{\line(0,1){2}}
\multiput(0,2)(0,1){3}{\line(1,0){4}}
\put(0,3){\makebox(1,1){$1$}}
\put(1,3){\makebox(1,1){$2$}}
\put(2,3){\makebox(1,1){$3$}}
\put(2,3){\makebox(1,1){$\bigcirc$}}
\put(3,3){\makebox(1,1){$\ol{3}$}}
\put(3,3){\makebox(1,1){$\bigcirc$}}
\put(0,2){\makebox(1,1){$2$}}
\put(0,2){\makebox(1,1){$\bigcirc$}}
\put(1,2){\makebox(1,1){$3$}}
\put(1,2){\makebox(1,1){$\bigcirc$}}
\put(2,2){\makebox(1,1){$\circ$}}
\put(3,2){\makebox(1,1){$\ol{2}$}}
\put(4,1){\makebox(1,1){$\uparrow$}}
\put(4,0){\makebox(1,1){$\ol{3}$}}
\end{picture}
\quad
\Rightarrow
\quad
\begin{picture}(5,4)(0,1)
\multiput(0,2)(1,0){5}{\line(0,1){2}}
\put(5,3){\line(0,1){1}}
\multiput(0,3)(0,1){2}{\line(1,0){5}}
\put(0,2){\line(1,0){4}}
\put(0,3){\makebox(1,1){$1$}}
\put(1,3){\makebox(1,1){$2$}}
\put(2,3){\makebox(1,1){$3$}}
\put(2,3){\makebox(1,1){$\bigcirc$}}
\put(3,3){\makebox(1,1){$\ol{3}$}}
\put(3,3){\makebox(1,1){$\bigcirc$}}
\put(0,2){\makebox(1,1){$2$}}
\put(0,2){\makebox(1,1){$\bigcirc$}}
\put(1,2){\makebox(1,1){$3$}}
\put(1,2){\makebox(1,1){$\bigcirc$}}
\put(2,2){\makebox(1,1){$\circ$}}
\put(3,2){\makebox(1,1){$\ol{2}$}}
\put(4,3){\makebox(1,1){$\ol{3}$}}
\put(4,3){\makebox(1,1){$\bigcirc$}}
\end{picture}
\end{displaymath}
\end{example}

\begin{lemma}
\label{lem:bcblx}
The bumping route does not move down.
\end{lemma}
\begin{proof}
It suffices to consider the bumping processes occurring
on pairs of neighboring columns
in the tableau.
Our strategy is as follows.
We are to show that if in the left column the inserted letter sets into
the first row then the same occurs in the right column as well.
Let us classify the situations on the neighboring columns into five cases.
\begin{enumerate}
\item
Suppose that in the following column insertion Case
B0 has occurred in the first column.
\begin{equation}
\setlength{\unitlength}{5mm}
\begin{picture}(5,1.4)(0,0.3)
\multiput(0,0)(1,0){2}{\line(0,1){1}}
\multiput(0,0)(0,1){2}{\line(1,0){1}}
\multiput(2.5,0)(1,0){3}{\line(0,1){1}}
\multiput(2.5,0)(0,1){2}{\line(1,0){2}}
\put(0,0){\makebox(1,1){$\alpha$}}
\put(1,0){\makebox(1.5,1){$\longrightarrow$}}
\put(2.5,0){\makebox(1,1){$\beta$}}
\put(3.5,0){\makebox(1,1){$\gamma$}}
\end{picture}
\end{equation}
Then in the second column Case B0 occurs and Case A1 does not happen.

The semistandard condition for $B$ tableau imposes that 
$(\beta,\gamma) \ne (\circ,\circ)$ and $\beta \leq \gamma$.
Thus if $\beta$ is bumped out from the left column, it certainly bumps 
$\gamma$ out of the right column.
\item
Suppose that in the following column insertion one of the Cases
B2, B4 or B6 has occurred in the first column.
\begin{equation}
\setlength{\unitlength}{5mm}
\begin{picture}(2.5,2.4)(0,0.3)
\multiput(0,0)(1,0){2}{\line(0,1){1}}
\multiput(0,0)(0,1){2}{\line(1,0){1}}
\put(0,0){\makebox(1,1){$\alpha$}}
\put(1,0){\makebox(1.5,1){$\longrightarrow$}}
\multiput(2.5,1)(0,1){2}{\line(1,0){2}}
\multiput(2.5,0)(1,0){2}{\line(0,1){2}}
\put(2.5,0){\line(1,0){1}}
\put(4.5,1){\line(0,1){1}}
\put(2.5,0){\makebox(1,1){$\delta$}}
\put(2.5,1){\makebox(1,1){$\beta$}}
\put(3.5,1){\makebox(1,1){$\gamma$}}
\end{picture}
\end{equation}
Then in the second column Case B0 occurs and Case A1 does not happen.

The reason is as follows.
Whichever one of the B2, B4 or B6 may have occurred in the first column,
the letter bumped out from the first column is always $\beta$.
And again we have
the semistandard condition between $\beta$ and $\gamma$.

\item
In the following column insertion Case
B7 occurs in the first column.
\begin{equation}
\setlength{\unitlength}{5mm}
\begin{picture}(2.5,2.4)(0,0.3)
\multiput(0,0)(1,0){2}{\line(0,1){1}}
\multiput(0,0)(0,1){2}{\line(1,0){1}}
\put(0,0){\makebox(1,1){$n$}}
\put(1,0){\makebox(1.5,1){$\longrightarrow$}}
\multiput(2.5,1)(0,1){2}{\line(1,0){2}}
\multiput(2.5,0)(1,0){2}{\line(0,1){2}}
\put(2.5,0){\line(1,0){1}}
\put(4.5,1){\line(0,1){1}}
\put(2.5,0){\makebox(1,1){$\ol{n}$}}
\put(2.5,1){\makebox(1,1){$\circ$}}
\put(3.5,1){\makebox(1,1){$\gamma$}}
\end{picture}
\end{equation}
Then in the second column Case B0 occurs and Case A1 does not happen.

The letter bumped out from the first column is $\ol{n}$.
Due to the semistandard condition we have $\gamma \geq \ol{n}$, hence the claim
follows.
\item
Suppose that in the following column insertion one of the Cases
B2, B4 or B6 has occurred in the first column.
\begin{equation}
\setlength{\unitlength}{5mm}
\begin{picture}(2.5,2.4)(0,0.3)
\multiput(0,0)(1,0){2}{\line(0,1){1}}
\multiput(0,0)(0,1){2}{\line(1,0){1}}
\multiput(2.5,0)(1,0){3}{\line(0,1){1}}
\multiput(2.5,0)(0,1){2}{\line(1,0){2}}
\put(0,0){\makebox(1,1){$\alpha$}}
\put(1,0){\makebox(1.5,1){$\longrightarrow$}}
\multiput(2.5,0)(0,1){3}{\line(1,0){2}}
\multiput(2.5,0)(1,0){3}{\line(0,1){2}}
\put(2.5,0){\makebox(1,1){$\delta$}}
\put(2.5,1){\makebox(1,1){$\beta$}}
\put(3.5,0){\makebox(1,1){$\varepsilon$}}
\put(3.5,1){\makebox(1,1){$\gamma$}}
\end{picture}
\end{equation}
Then in the second column Cases B1, B3 and B5 do not happen.

The reason is as follows.
The letter bumped out from the first column is always $\beta$.

Since $\beta \leq \gamma$, B1 does not happen.
Since $(\beta,\gamma) \ne (\circ,\circ)$, B3 does not happen.
B5 does not happen since $(\beta,\gamma,\veps) \ne (x,x,\ol{x})$, i.e. due to
the absence of the $(x,x)$-configuration (\ref{eq:absenceofxxconf}).

\item
In the following column insertion Case
B7 occurs in the first column.
\begin{equation}
\setlength{\unitlength}{5mm}
\begin{picture}(2.5,2.4)(0,0.3)
\multiput(0,0)(1,0){2}{\line(0,1){1}}
\multiput(0,0)(0,1){2}{\line(1,0){1}}
\multiput(2.5,0)(1,0){3}{\line(0,1){1}}
\multiput(2.5,0)(0,1){2}{\line(1,0){2}}
\put(0,0){\makebox(1,1){$n$}}
\put(1,0){\makebox(1.5,1){$\longrightarrow$}}
\multiput(2.5,0)(0,1){3}{\line(1,0){2}}
\multiput(2.5,0)(1,0){3}{\line(0,1){2}}
\put(2.5,0){\makebox(1,1){$\ol{n}$}}
\put(2.5,1){\makebox(1,1){$\circ$}}
\put(3.5,0){\makebox(1,1){$\varepsilon$}}
\put(3.5,1){\makebox(1,1){$\gamma$}}
\end{picture}
\end{equation}
Then in the second column Cases B1, B3 and B5 do not happen.

The letter bumped out from the first column is $\ol{n}$.
Due to the semistandard condition we have $\gamma \geq \ol{n}$,
hence the claim follows.
\end{enumerate}
\end{proof}
\begin{lemma}
\label{lem:bcblxx}
Let $\alpha' \leq \alpha$ and $(\alpha,\alpha') \ne (\circ,\circ)$.
Let $R$ be the bumping route that is made when $\alpha$ is inserted into $T$,
and $R'$ be the bumping route that is made when 
$\alpha'$ is inserted into $\left( \alpha \longrightarrow T \right)$.
Then $R'$ does not lie below $R$.
\end{lemma}
\begin{proof}
First we consider the case where the bumping route
lies only in the first row.
Suppose that, when $\alpha$ was inserted
into the tableau $T$, it was set in the first row
in the first column.
We are to show that when $\alpha'$ is inserted, 
it will be also set in the first row in the first column.
If $T$ is an empty set (resp.~has only one row), the insertion of
$\alpha$ should have been  A0 (resp.~B0).
In either case we have B0 when $\alpha'$ is inserted, hence the claim is true.
Suppose $T$ has two rows.
By assumption B2, B4, B6 or B7 has occurred when $\alpha$ was inserted.
We see that,
if B4, B6 or B7 has occurred, 
then B2 will occur when $\alpha'$ is inserted.
Thus it is enough to show that 
if B2 has occurred,
then B1, B3 or B5 does not happen when $\alpha'$ is inserted.
Since $\alpha' \leq \alpha$, B1 does not happen.
Since $(\alpha,\alpha') \ne (\circ,\circ)$, B3 does not happen.
B5 does not happen, since the first column does not have the entry
${x \atop \overline{x}}$  as the result of B2 type insertion of
$\alpha$.

Second we consider the case where the bumping route $R$ lies across 
the first and the second rows.
Suppose that from the leftmost column to
the $(i-1)$-th column the bumping route
lies in the second row, and from
the $i$-th column to the rightmost column it lies in the first row.
Let us call the position of the vertical line between the
$(i-1)$-th and the $i$-th columns the {\em crossing point} of $R$.
It is unique due to Lemma \ref{lem:bcblx}.
We call an analogous position of $R'$ its crossing point.
We are to show that the crossing point of $R'$ does not locate
strictly right to the crossing point of $R$.
Let the situation around the crossing point of $R$ be
\begin{equation}
\label{eq:crptb}
\setlength{\unitlength}{5mm}
\begin{picture}(4,2.4)(0,0.3)
\multiput(1,0)(1,0){3}{\line(0,1){2}}
\multiput(1,0)(0,1){3}{\line(1,0){2}}
\multiput(0.1,0)(0,1){3}{
\multiput(0,0)(0.2,0){5}{\line(1,0){0.1}}}
\multiput(3,0)(0,1){3}{
\multiput(0,0)(0.2,0){5}{\line(1,0){0.1}}}
\put(1,0){\makebox(1,1){$\xi$}}
\put(2,1){\makebox(1,1){$\eta$}}
\end{picture}
\quad
\mbox{or}
\quad
\setlength{\unitlength}{5mm}
\begin{picture}(4,2.4)(0,0.3)
\multiput(1,0)(1,0){2}{\line(0,1){2}}
\put(3,1){\line(0,1){1}}
\multiput(1,1)(0,1){2}{\line(1,0){2}}
\put(1,0){\line(1,0){1}}
\multiput(0.1,0)(0,1){3}{
\multiput(0,0)(0.2,0){5}{\line(1,0){0.1}}}
\multiput(3,1)(0,1){2}{
\multiput(0,0)(0.2,0){5}{\line(1,0){0.1}}}
\put(1,0){\makebox(1,1){$\xi$}}
\put(2,1){\makebox(1,1){$\eta$}}
\end{picture}
\,
\mbox{.}
\end{equation}
While the insertion of $\alpha$ that led to these configurations,
let $\eta'$ be the letter that was bumped out from the 
left column. 

Claim 1: $\xi \leq \eta' \le \eta$ and $(\xi,\eta) \ne (\circ,\circ)$.
To see this note that
in the left column,
B1, B3 or B5 has occurred when $\alpha$ was inserted.
We have $\xi \leq \eta'$ and $(\xi,\eta') \ne (\circ,\circ)$ (B1), 
or $\xi < \eta'$ (B3, B5).
In the right column
A0, B0, B2, B4, B6 or B7 has subsequently occurred.
We have $\eta' = \eta$ (A0, B0, B2, B4, B7), or $\eta' < \eta$ (B6).
In any case we have $\xi \le \eta' \le \eta$ and $(\xi,\eta) \ne (\circ,\circ)$.

Claim 2: In (\ref{eq:crptb}) the following configurations do not exist.
\begin{equation}
\setlength{\unitlength}{5mm}
\begin{picture}(4,2.4)(0,0.3)
\multiput(1,0)(1,0){3}{\line(0,1){2}}
\multiput(1,0)(0,1){3}{\line(1,0){2}}
\multiput(0.1,0)(0,1){3}{
\multiput(0,0)(0.2,0){5}{\line(1,0){0.1}}}
\multiput(3,0)(0,1){3}{
\multiput(0,0)(0.2,0){5}{\line(1,0){0.1}}}
\put(1,0){\makebox(1,1){$\ol{x}$}}
\put(1,1){\makebox(1,1){$x$}}
\put(2,1){\makebox(1,1){$\ol{x}$}}
\end{picture}
\,
\mbox{,}
\,
\setlength{\unitlength}{5mm}
\begin{picture}(4,2.4)(0,0.3)
\multiput(1,0)(1,0){2}{\line(0,1){2}}
\put(3,1){\line(0,1){1}}
\multiput(1,1)(0,1){2}{\line(1,0){2}}
\put(1,0){\line(1,0){1}}
\multiput(0.1,0)(0,1){3}{
\multiput(0,0)(0.2,0){5}{\line(1,0){0.1}}}
\multiput(3,1)(0,1){2}{
\multiput(0,0)(0.2,0){5}{\line(1,0){0.1}}}
\put(1,0){\makebox(1,1){$\ol{x}$}}
\put(1,1){\makebox(1,1){$x$}}
\put(2,1){\makebox(1,1){$\ol{x}$}}
\end{picture}
\,
\mbox{,}
\,
\begin{picture}(4,2.4)(0,0.3)
\multiput(1,0)(1,0){3}{\line(0,1){2}}
\multiput(1,0)(0,1){3}{\line(1,0){2}}
\multiput(0.1,0)(0,1){3}{
\multiput(0,0)(0.2,0){5}{\line(1,0){0.1}}}
\multiput(3,0)(0,1){3}{
\multiput(0,0)(0.2,0){5}{\line(1,0){0.1}}}
\put(1,0){\makebox(1,1){$x$}}
\put(2,0){\makebox(1,1){$\ol{x}$}}
\put(2,1){\makebox(1,1){$x$}}
\end{picture}
\,
\mbox{,}
\,
\begin{picture}(4,2.4)(0,0.3)
\multiput(1,0)(1,0){3}{\line(0,1){2}}
\multiput(1,0)(0,1){3}{\line(1,0){2}}
\multiput(0.1,0)(0,1){3}{
\multiput(0,0)(0.2,0){5}{\line(1,0){0.1}}}
\multiput(3,0)(0,1){3}{
\multiput(0,0)(0.2,0){5}{\line(1,0){0.1}}}
\put(1,0){\makebox(1,1){$\circ$}}
\put(2,1){\makebox(1,1){$\circ$}}
\end{picture}
\,
\mbox{or}
\,
\setlength{\unitlength}{5mm}
\begin{picture}(4,2.4)(0,0.3)
\multiput(1,0)(1,0){2}{\line(0,1){2}}
\put(3,1){\line(0,1){1}}
\multiput(1,1)(0,1){2}{\line(1,0){2}}
\put(1,0){\line(1,0){1}}
\multiput(0.1,0)(0,1){3}{
\multiput(0,0)(0.2,0){5}{\line(1,0){0.1}}}
\multiput(3,1)(0,1){2}{
\multiput(0,0)(0.2,0){5}{\line(1,0){0.1}}}
\put(1,0){\makebox(1,1){$\circ$}}
\put(2,1){\makebox(1,1){$\circ$}}
\end{picture}
\,
\mbox{.}
\label{eq:bforbiddenconfs}
\end{equation}
Due to Claim 1,
the first and the second configurations can exist only if
B1 with $\alpha= x, \gamma=\eta'$ happens in the left column 
and $\xi = \eta'= \eta = \overline{x}$.
But $(\alpha,\gamma) = (x,\overline{x})$ is not compatible with B1.
The third configuration can exist only if 
B6 happens in the right column and $\xi= \eta'= \eta= x$ by Claim 1.
But we see from the proof of Claim 1 that B6 actually happens only when $\eta' < \eta$.
The fourth and the fifth are forbidden since 
$(\xi,\eta) \neq (\circ,\circ)$ by Claim 1.
Claim 2 is proved.

Let the situation around the crossing point of $R$ be one of (\ref{eq:crptb})
excluding (\ref{eq:bforbiddenconfs}).
When inserting $\alpha'$, 
suppose in the left column of the crossing point,  B1, B3 or B5 has occurred.
Let $\xi'$ be the letter bumped out therefrom.

Claim 3: $\xi' \leq \eta$ and $(\xi',\eta) \ne (\circ,\circ)$.
We divide the check into two cases.
a) If B1 or B3 has occurred in the left column, we have $\xi' = \xi$.
Thus the assertion follows from Claim 1.
b) If B5 has occurred, the left column had the entry ${x \atop \ol{x}}$ and 
we have $\xi' = \ol{x-1}$, $\xi = \ol{x}$.
Claim 1 tells $\xi = \ol{x} \le \eta$, and Claim 2 does $\eta \neq \ol{x}$.
Therefore we have $\xi' = \ol{x-1} \le \eta$.
$(\xi',\eta) \ne (\circ,\circ)$ is obvious.
Claim 3 is proved.

Now we are ready to finish the proof of the main assertion.
Assume the same situation as Claim 3.
We should verify that A1, B1, B3 and B5 do not occur in the right column.
Claim 3 immediately prohibits A1, B1 and B3 in the right column.
Suppose that B5 happens in the right column.
It means that 
$\eta \in \{1,\ldots n\}$, $\xi' \geq \eta$ and the right column had the entry
${\eta \atop \ol{\eta}}$.
Since $\xi' \leq \eta$ by Claim 3, we find $\xi' = \eta$, therefore
$\xi' \in \{1,\ldots, n\}$.
Such $\xi'$ can be bumped out from  B1 process only in the left column 
and not from B3 or B5.
It follows that $\xi' = \xi$.
This leads to the third configuration in (\ref{eq:bforbiddenconfs}), 
hence a contradiction.

Finally we consider the case where the bumping route $R$
lies only in the second row.
If $R'$ lies below $R$ the tableau should have 
more than two rows, which is prohibited by
Proposition \ref{pr:atmosttworows}.
\end{proof}
\begin{coro}
\label{cor:bcblxxx}
Let $\alpha' \leq \alpha$ and $(\alpha,\alpha') \ne (\circ,\circ)$.
Suppose that a new box is added at the end of the first row
when $\alpha$ is inserted into $T$.
Then a new box is added also at the end of the first row
when $\alpha'$ is inserted into $\left( \alpha \longrightarrow T \right)$.
\end{coro}
\subsection{\mathversion{bold}Column insertion and inverse insertion for
$D_n$}
\label{subsec:cid}
Set an alphabet $\mathcal{X}=\mathcal{A} \sqcup \bar{\mathcal{A}},\,
\mathcal{A}=\{ 1,\dots,n\}$ and
$\bar{\mathcal{A}}=\{\overline{1},\dots,\overline{n}\}$,
with the partial order
$1 < 2 < \dots < {n \atop \ol{n}} < \dots < \overline{2} < \overline{1}$.
\subsubsection{Semistandard $D$ tableaux}
\label{subsubsec:ssdt}
Let us consider a {\em semistandard $D$ tableau}
made by the letters from this alphabet.
We follow \cite{KN} for its definition.
We present the definition here, but restrict ourselves to
special cases that are sufficient for our purpose.
Namely we consider only
those tableaux that have no more than two rows
in their shapes. 
Thus they have the forms as in (\ref{eq:pictureofsst}),
with the letters inside the boxes being chosen from the alphabet 
given as above.
The letters obey the conditions
(\ref{eq:notdecrease}),\footnote{Note that there is no order between $n$ and $\ol{n}$.} 
(\ref{eq:absenceofoneonebar})
 and the absence of the $(x,x)$-configuration
(\ref{eq:absenceofxxconf}) where we now assume $1 \leq x <n$.
They also obey the following conditions:
\begin{equation}
\alpha_a < \beta_a \quad \mbox{or} \quad (\alpha_a,\beta_a) = (n,\ol{n})
\quad \mbox{or} \quad (\alpha_a,\beta_a) = (\ol{n},n),
\end{equation}
\begin{equation}
(\alpha_a,\alpha_{a+1},\beta_a,\beta_{a+1}) \ne (n-1,n,n,\ol{n-1}),
(n-1,\ol{n},\ol{n},\ol{n-1}),
\end{equation}
\begin{equation}
\label{eq:absenceofnnconfd}
(\alpha_a,\beta_{a+1}) \ne (n,n),(n,\ol{n}),(\ol{n},n),(\ol{n},\ol{n}).
\end{equation}
The conditions (\ref{eq:absenceofnnconfd}) are referred to as the
absence of the $(n,n)$-configurations.

\subsubsection{\mathversion{bold}Column insertion for 
$D_n$ \cite{B2}}
\label{subsubsec:insd}
We give a list of column insertions on 
semistandard $D$ tableaux that are sufficient for our purpose.
For the alphabet $\mathcal{X}$, we follow the convention
that Greek letters $ \alpha, \beta, \ldots $ belong to
$\mathcal{X}$ while Latin
letters $x,y,\ldots$ (resp. $\overline{x},\overline{y},\ldots$) belong to
$\mathcal{A}$ (resp. $\bar{\mathcal{A}}$).
For the interpretation of the pictorial equations in the list,
see the remarks in Section
\ref{subsubsec:insc}.
Note that this list does not exhaust all cases.
It does not contain, for instance, the insertion
\begin{math}
\setlength{\unitlength}{3mm}
\begin{picture}(3,2)(0,0.3)
\multiput(0,0)(1,0){2}{\line(0,1){1}}
\multiput(0,0)(0,1){2}{\line(1,0){1}}
\put(0,0){\makebox(1,1){${\scriptstyle n}$}}
\put(1,0){\makebox(1,1){${\scriptstyle \rightarrow}$}}
\multiput(2,0)(1,0){2}{\line(0,1){2}}
\multiput(2,0)(0,1){3}{\line(1,0){1}}
\put(2,0){\makebox(1,1){${\scriptstyle \ol{n}}$}}
\put(2,1){\makebox(1,1){${\scriptstyle n}$}}
\end{picture}
\end{math}.
(See Proposition \ref{pr:atmosttworows}.)

\noindent
\ToBOX{A0}{\alpha}
	\raisebox{1.25mm}{,}

\noindent
\ToDOMINO{A1}{\alpha}{\beta}{\alpha}{\beta}%
	\raisebox{4mm}{if $\alpha < \beta$ or $(\alpha,\beta)=(n,\overline{n})$ or $(\overline{n},n)$,}

\noindent
\ToYOKODOMINO{B0}{\beta}{\alpha}{\beta}{\alpha}%
	\raisebox{1.25mm}{if $\alpha \le \beta$,}

\noindent
\ToHOOKnn{B1}{\alpha}{\gamma}{\beta}{\gamma}{\alpha}{\beta}%
	\raisebox{4mm}{if $\alpha < \beta \leq \gamma$ and $(\alpha,\gamma) \ne (x,\overline{x})$,}

\noindent
\ToHOOKnn{B2}{\beta}{\gamma}{\alpha}{\beta}{\alpha}{\gamma}
	\raisebox{4mm}{if $\alpha \leq \beta < \gamma$ and $(\alpha,\gamma) \ne (x,\overline{x})$,}

\noindent
\ToHOOKll{B3}{x}{\overline{x}}{\beta}{\overline{x\!-\!1}}{x\!-\!1}{\beta}%
	\raisebox{4mm}{if $x \le \beta \le \overline{x}$ and $x\ne n,1\,$,}

\noindent
\ToHOOKln{B4}{\beta}{\overline{x}}{x}{\beta}{x\!+\!1}{\overline{x\!+\!1}}%
	\raisebox{4mm}{if $x < \beta < \overline{x}$ and $x\ne n\!-\!1,n\,$,}

\noindent
\ToHOOKnn{B5}{\mu_1}{\mu_2}{x}{\mu_1}{x}{\mu_2}%
	\raisebox{4mm}{if $(\mu_1,\mu_2) = (n,\overline{n})$ or $(\overline{n},n)$ and
	$x \ne n$,}

\noindent
\ToHOOKnn{B6}{\mu_1}{\overline{x}}{\mu_2}{\overline{x}}{\mu_1}{\mu_2}%
	\raisebox{4mm}{if $(\mu_1,\mu_2) = (n,\overline{n})$ or $(\overline{n},n)$ and
	$\overline{x} \ne \overline{n}$,}

\noindent
\ToHOOKllnn{B7}{\mu}{\overline{n\!-\!1}}{n\!-\!1}{\mu}{\mu}{\overline{\mu}}%
	\raisebox{4mm}{if $\mu = n$ or $\overline{n}\;\;
	(\overline{\mu}:=n$ if $\mu=\overline{n}$),}

\noindent
\ToHOOKll{B8}{\mu_1}{\mu_2}{\mu_2}{\overline{n\!-\!1}}{n\!-\!1}{\mu_2}%
	\raisebox{4mm}{if $(\mu_1,\mu_2) = (n,\overline{n})$ or $(\overline{n},n)$.}

\noindent
\subsubsection{Column insertion and $U_q (D_n)$ crystal morphism}
To illustrate Proposition \ref{pr:morphgen}
let us check a morphism of the $U_q(D_4)$ crystal
$B(\Lambda_2) \ot B(\Lambda_1)$ by taking two examples.
Let $\psi$ be the map that is similarly defined
as in Section \ref{subsubsec:insc} for type $C$ case.
\begin{example}
\label{ex:morD1}
\begin{displaymath}
\begin{CD}
\setlength{\unitlength}{5mm}
\begin{picture}(3,2)(0,0.3)
\multiput(0,0)(1,0){2}{\line(0,1){2}}
\multiput(0,0)(0,1){3}{\line(1,0){1}}
\put(0,1){\makebox(1,1){$\ol{4}$}}
\put(0,0){\makebox(1,1){$4$}}
\put(1,0.5){\makebox(1,1){$\otimes$}}
\put(2,0.5){\makebox(1,1){$3$}}
\multiput(2,0.5)(1,0){2}{\line(0,1){1}}
\multiput(2,0.5)(0,1){2}{\line(1,0){1}}
\end{picture}
@>\text{$\ftd_4$}>>
\setlength{\unitlength}{5mm}
\begin{picture}(3,2)(0,0.3)
\multiput(0,0)(1,0){2}{\line(0,1){2}}
\multiput(0,0)(0,1){3}{\line(1,0){1}}
\put(0,1){\makebox(1,1){$\ol{4}$}}
\put(0,0){\makebox(1,1){$\ol{3}$}}
\put(1,0.5){\makebox(1,1){$\otimes$}}
\put(2,0.5){\makebox(1,1){$3$}}
\multiput(2,0.5)(1,0){2}{\line(0,1){1}}
\multiput(2,0.5)(0,1){2}{\line(1,0){1}}
\end{picture}
@>\text{$\ftd_4$}>>
\setlength{\unitlength}{5mm}
\begin{picture}(3,2)(0,0.3)
\multiput(0,0)(1,0){2}{\line(0,1){2}}
\multiput(0,0)(0,1){3}{\line(1,0){1}}
\put(0,1){\makebox(1,1){$\ol{4}$}}
\put(0,0){\makebox(1,1){$\ol{3}$}}
\put(1,0.5){\makebox(1,1){$\otimes$}}
\put(2,0.5){\makebox(1,1){$\ol{4}$}}
\multiput(2,0.5)(1,0){2}{\line(0,1){1}}
\multiput(2,0.5)(0,1){2}{\line(1,0){1}}
\end{picture}
\\
@VV\text{$\psi$}V @VV\text{$\psi$}V @VV\text{$\psi$}V \\
\setlength{\unitlength}{5mm}
\begin{picture}(2,2)(0,0.3)
\multiput(0,0)(1,0){2}{\line(0,1){2}}
\put(0,0){\line(1,0){1}}
\multiput(0,1)(0,1){2}{\line(1,0){2}}
\put(2,1){\line(0,1){1}}
\put(0,1){\makebox(1,1){$3$}}\put(1,1){\makebox(1,1){$\ol{4}$}}
\put(0,0){\makebox(1,1){$4$}}
\end{picture}
@>\text{$\ftd_4$}>>
\setlength{\unitlength}{5mm}
\begin{picture}(2,2)(0,0.3)
\multiput(0,0)(1,0){2}{\line(0,1){2}}
\put(0,0){\line(1,0){1}}
\multiput(0,1)(0,1){2}{\line(1,0){2}}
\put(2,1){\line(0,1){1}}
\put(0,1){\makebox(1,1){$\ol{4}$}}\put(1,1){\makebox(1,1){$\ol{4}$}}
\put(0,0){\makebox(1,1){$4$}}
\end{picture}
@>\text{$\ftd_4$}>>
\setlength{\unitlength}{5mm}
\begin{picture}(2,2)(0,0.3)
\multiput(0,0)(1,0){2}{\line(0,1){2}}
\put(0,0){\line(1,0){1}}
\multiput(0,1)(0,1){2}{\line(1,0){2}}
\put(2,1){\line(0,1){1}}
\put(0,1){\makebox(1,1){$\ol{4}$}}\put(1,1){\makebox(1,1){$\ol{4}$}}
\put(0,0){\makebox(1,1){$\ol{3}$}}
\end{picture}
\end{CD}
\end{displaymath}
\vskip3ex
\noindent
Here the $\psi$'s are given by Case B5, B7 and B2 column insertions,
respectively.
\end{example}
\begin{example}
\label{ex:morD2}
\begin{displaymath}
\begin{CD}
\setlength{\unitlength}{5mm}
\begin{picture}(3,2)(0,0.3)
\multiput(0,0)(1,0){2}{\line(0,1){2}}
\multiput(0,0)(0,1){3}{\line(1,0){1}}
\put(0,1){\makebox(1,1){$\ol{4}$}}
\put(0,0){\makebox(1,1){$\ol{3}$}}
\put(1,0.5){\makebox(1,1){$\otimes$}}
\put(2,0.5){\makebox(1,1){$4$}}
\multiput(2,0.5)(1,0){2}{\line(0,1){1}}
\multiput(2,0.5)(0,1){2}{\line(1,0){1}}
\end{picture}
@>\text{$\etd_4$}>>
\setlength{\unitlength}{5mm}
\begin{picture}(3,2)(0,0.3)
\multiput(0,0)(1,0){2}{\line(0,1){2}}
\multiput(0,0)(0,1){3}{\line(1,0){1}}
\put(0,1){\makebox(1,1){$\ol{4}$}}
\put(0,0){\makebox(1,1){$4$}}
\put(1,0.5){\makebox(1,1){$\otimes$}}
\put(2,0.5){\makebox(1,1){$4$}}
\multiput(2,0.5)(1,0){2}{\line(0,1){1}}
\multiput(2,0.5)(0,1){2}{\line(1,0){1}}
\end{picture}
@>\text{$\etd_4$}>>
\setlength{\unitlength}{5mm}
\begin{picture}(3,2)(0,0.3)
\multiput(0,0)(1,0){2}{\line(0,1){2}}
\multiput(0,0)(0,1){3}{\line(1,0){1}}
\put(0,1){\makebox(1,1){$3$}}
\put(0,0){\makebox(1,1){$4$}}
\put(1,0.5){\makebox(1,1){$\otimes$}}
\put(2,0.5){\makebox(1,1){$4$}}
\multiput(2,0.5)(1,0){2}{\line(0,1){1}}
\multiput(2,0.5)(0,1){2}{\line(1,0){1}}
\end{picture}
\\
@VV\text{$\psi$}V @VV\text{$\psi$}V @VV\text{$\psi$}V \\
\setlength{\unitlength}{5mm}
\begin{picture}(2,2)(0,0.3)
\multiput(0,0)(1,0){2}{\line(0,1){2}}
\put(0,0){\line(1,0){1}}
\multiput(0,1)(0,1){2}{\line(1,0){2}}
\put(2,1){\line(0,1){1}}
\put(0,1){\makebox(1,1){$\ol{4}$}}\put(1,1){\makebox(1,1){$\ol{3}$}}
\put(0,0){\makebox(1,1){$4$}}
\end{picture}
@>\text{$\etd_4$}>>
\setlength{\unitlength}{5mm}
\begin{picture}(2,2)(0,0.3)
\multiput(0,0)(1,0){2}{\line(0,1){2}}
\put(0,0){\line(1,0){1}}
\multiput(0,1)(0,1){2}{\line(1,0){2}}
\put(2,1){\line(0,1){1}}
\put(0,1){\makebox(1,1){$3$}}\put(1,1){\makebox(1,1){$\ol{3}$}}
\put(0,0){\makebox(1,1){$4$}}
\end{picture}
@>\text{$\etd_4$}>>
\setlength{\unitlength}{5mm}
\begin{picture}(2,2)(0,0.3)
\multiput(0,0)(1,0){2}{\line(0,1){2}}
\put(0,0){\line(1,0){1}}
\multiput(0,1)(0,1){2}{\line(1,0){2}}
\put(2,1){\line(0,1){1}}
\put(0,1){\makebox(1,1){$3$}}\put(1,1){\makebox(1,1){$4$}}
\put(0,0){\makebox(1,1){$4$}}
\end{picture}
\end{CD}
\end{displaymath}
\vskip3ex
\noindent
Here the $\psi$'s are given by Case B6, B8 and B1 column insertions,
respectively.
\end{example}
\subsubsection{\mathversion{bold} Inverse insertion for 
$D_n$ \cite{B2}}
\label{subsubsec:invinsd}
We give a list of inverse column insertions on 
semistandard $D$ tableaux that are sufficient for our purpose.
For the interpretation of the pictorial equations in the list,
see the remarks in Section
\ref{subsubsec:invinsc}.
\par
\noindent
\FromYOKODOMINO{C0}{\beta}{\alpha}{\beta}{\alpha}%
	\raisebox{1.25mm}{if $\alpha \le \beta$,}

\noindent
\FromHOOKnn{C1}{\gamma}{\alpha}{\beta}{\alpha}{\gamma}{\beta}%
	\raisebox{4mm}{if $\alpha < \beta \leq \gamma$ and $(\alpha,\gamma) \ne (x,\overline{x})$,}

\noindent
\FromHOOKnn{C2}{\beta}{\alpha}{\gamma}{\beta}{\gamma}{\alpha}%
	\raisebox{4mm}{if $\alpha \leq \beta < \gamma$ and $(\alpha,\gamma) \ne (x,\overline{x})$,}

\noindent
\FromHOOKnl{C3}{\overline{x}}{x}{\beta}{x\!+\!1}{\overline{x\!+\!1}}{\beta}%
	\raisebox{4mm}{if $x < \beta < \overline{x}$ and $x\ne n\!-\!1,n\,$,}

\noindent
\FromHOOKll{C4}{\beta}{x}{\overline{x}}{\beta}{\overline{x\!-\!1}}{x\!-\!1}%
	\raisebox{4mm}{if $x \le \beta \le \overline{x}$ and $x\ne n,1\,$,}

\noindent
\FromHOOKnn{C5}{\mu_1}{x}{\mu_2}{\mu_1}{\mu_2}{x}%
	\raisebox{4mm}{if $(\mu_1,\mu_2) = (n,\overline{n})$ or $(\overline{n},n)$ and
	$x \ne n$}

\noindent
\FromHOOKnn{C6}{\overline{x}}{\mu_1}{\mu_2}{\mu_1}{\overline{x}}{\mu_2}%
	\raisebox{4mm}{if $(\mu_1,\mu_2) = (n,\overline{n})$ or $(\overline{n},n)$ and
	$\overline{x} \ne \overline{n}$,}

\noindent
\FromHOOKll{C7}{\mu_1}{\mu_1}{\mu_2}{\mu_1}{\overline{n\!-\!1}}{n\!-\!1}%
	\raisebox{4mm}{if $(\mu_1,\mu_2) = (n,\overline{n})$ or $(\overline{n},n)$,}

\noindent
\FromHOOKllnn{C8}{\overline{n\!-\!1}}{n\!-\!1}{\mu}{\overline{\mu}}{\mu}{\mu}%
	\raisebox{4mm}{if $\mu = n$ or $\overline{n} \;\;
	 (\overline{\mu}:=n$ if $\mu=\overline{n}$).}
\subsubsection{\mathversion{bold} Column bumping lemma for 
$D_n$}
\label{subsubsec:cblD}
The aim of this subsection is to give a simple result on
successive insertions of two letters into a 
tableau (Corollary \ref{cor:cblxxx}).
This result will be
used in the proof of the main theorem (Theorem \ref{th:main3}).
This corollary follows from the column bumping lemma
(Lemma \ref{lem:cblxx}).

We restrict ourselves to the situation where by column insertions there only appear 
semistandard $D$ tableaux with at most two rows.
In the classification of the column insertions, we regard 
that the inserted letter is set in the first row in 
Cases A0, B0, B2, B4, B5 and
B7, and that it is set in the second row in the other cases.
Then the bumping route is defined in the same way as
in section~\ref{subsubsec:cblB}.
\begin{lemma}
\label{lem:cblx}
The bumping route does not move down.
\end{lemma}
\begin{proof}
It is enough to consider the following three cases.
\begin{enumerate}
\item
Suppose that in the following column insertion Case
B0 has occurred in the first column.
\begin{equation}
\setlength{\unitlength}{5mm}
\begin{picture}(5,1.4)(0,0.3)
\multiput(0,0)(1,0){2}{\line(0,1){1}}
\multiput(0,0)(0,1){2}{\line(1,0){1}}
\multiput(2.5,0)(1,0){3}{\line(0,1){1}}
\multiput(2.5,0)(0,1){2}{\line(1,0){2}}
\put(0,0){\makebox(1,1){$\alpha$}}
\put(1,0){\makebox(1.5,1){$\longrightarrow$}}
\put(2.5,0){\makebox(1,1){$\beta$}}
\put(3.5,0){\makebox(1,1){$\gamma$}}
\end{picture}
\end{equation}
Then in the second column Case B0 occurs and Case A1 does not happen.

The semistandard condition for $D$ tableau imposes that 
$(\beta,\gamma) \ne (n,\ol{n}), (\ol{n},n)$ and $\beta \leq \gamma$.
Thus if $\beta$ is bumped out from the left column, it certainly bumps 
$\gamma$ out of the right column.
\item
Suppose that in the following column insertion one of the Cases
B2, B4, B5 or B7 has occurred in the first column.
\begin{equation}
\setlength{\unitlength}{5mm}
\begin{picture}(2.5,2.4)(0,0.3)
\multiput(0,0)(1,0){2}{\line(0,1){1}}
\multiput(0,0)(0,1){2}{\line(1,0){1}}
\put(0,0){\makebox(1,1){$\alpha$}}
\put(1,0){\makebox(1.5,1){$\longrightarrow$}}
\multiput(2.5,1)(0,1){2}{\line(1,0){2}}
\multiput(2.5,0)(1,0){2}{\line(0,1){2}}
\put(2.5,0){\line(1,0){1}}
\put(4.5,1){\line(0,1){1}}
\put(2.5,0){\makebox(1,1){$\delta$}}
\put(2.5,1){\makebox(1,1){$\beta$}}
\put(3.5,1){\makebox(1,1){$\gamma$}}
\end{picture}
\end{equation}
Then in the second column Case B0 occurs and Case A1 does not happen.

The reason is as follows.
Whichever one of the B2, B4, B5 or B7 may have occurred in the first column,
the letter bumped out from the first column is always $\beta$. And we have
the semistandard condition between $\beta$ and $\gamma$.
\item
Suppose that in the following column insertion one of the Cases
B2, B4, B5 or B7 has occurred in the first column.
\begin{equation}
\setlength{\unitlength}{5mm}
\begin{picture}(2.5,2.4)(0,0.3)
\multiput(0,0)(1,0){2}{\line(0,1){1}}
\multiput(0,0)(0,1){2}{\line(1,0){1}}
\multiput(2.5,0)(1,0){3}{\line(0,1){1}}
\multiput(2.5,0)(0,1){2}{\line(1,0){2}}
\put(0,0){\makebox(1,1){$\alpha$}}
\put(1,0){\makebox(1.5,1){$\longrightarrow$}}
\multiput(2.5,0)(0,1){3}{\line(1,0){2}}
\multiput(2.5,0)(1,0){3}{\line(0,1){2}}
\put(2.5,0){\makebox(1,1){$\delta$}}
\put(2.5,1){\makebox(1,1){$\beta$}}
\put(3.5,0){\makebox(1,1){$\varepsilon$}}
\put(3.5,1){\makebox(1,1){$\gamma$}}
\end{picture}
\end{equation}
Then in the second column Cases B1, B3, B6 and B8 do not happen.

The reason is as follows.
As in the previous case
the letter bumped out from the first column is always $\beta$.

Since $\beta \leq \gamma$, B1 does not happen.
Since $(\beta,\gamma) \ne (n,\ol{n}), (\ol{n},n)$, B6 and B8 do not happen.
B3 does not happen since $(\beta,\gamma,\veps) \ne (x,x,\ol{x})$, i.e. due to
the absence of the $(x,x)$-configuration (\ref{eq:absenceofxxconf}).
\end{enumerate}
\end{proof}
\begin{lemma}
\label{lem:cblxx}
Let $\alpha' \leq \alpha$, 
in particular $(\alpha,\alpha') \ne (n,\ol{n}),(\ol{n},n)$.
Let $R$ be the bumping route that is made when $\alpha$ is inserted into $T$,
and $R'$ be the bumping route that is made when 
$\alpha'$ is inserted into $\left( \alpha \longrightarrow T \right)$.
Then $R'$ does not lie below $R$.
\end{lemma}
\begin{proof}
First we consider the case where the bumping route
lies only in the first row.
Suppose that, when $\alpha$ was inserted
into the tableau $T$, it was set in the first row
in the first column.
We are to show that when $\alpha'$ is inserted, 
it will be also set in the first row in the first column.
If $T$ is an empty set (resp.~has only one row), the insertion of
$\alpha$ should have been  A0 (resp.~B0).
In either case we have B0 when $\alpha'$ is inserted, hence the claim is true.
Suppose $T$ has two rows.
By assumption B2, B4, B5 or B7 has occurred when $\alpha$ was inserted.
We see that;
a) If B7 has occurred, 
then B5 will occur when $\alpha'$ is inserted;
b) If B5 or B4 has occurred, 
then B2 will occur when $\alpha'$ is inserted.
Thus it is enough to show that 
if B2 has occurred,
then B1, B3, B6 and B8 do not happen when $\alpha'$ is inserted.
Since $\alpha' \leq \alpha$, B1 does not happen.
Since $(\alpha , \alpha') \ne (n,\ol{n}),(\ol{n},n)$, B6 and B8 do not happen.
B3 does not happen, since the first column does not have the entry
${x \atop \overline{x}}$  as the result of B2 type insertion of
$\alpha$.

Second we consider the case where the bumping route $R$ lies across 
the first and the second rows.
Suppose that from the leftmost column to
the $(i-1)$-th column the bumping route
lies in the second row, and from
the $i$-th column to the rightmost column it lies in the first row.
As in the type $B$ case
let us call the position of the vertical line between the
$(i-1)$-th and the $i$-th columns the crossing point of $R$.
It is unique due to Lemma \ref{lem:cblx}.
We call an analogous position of $R'$ its crossing point.
We are to show that the crossing point of $R'$ does not locate
strictly right to the crossing point of $R$.
Let the situation around the crossing point of $R$ be
\begin{equation}
\label{eq:crptd}
\setlength{\unitlength}{5mm}
\begin{picture}(4,2.4)(0,0.3)
\multiput(1,0)(1,0){3}{\line(0,1){2}}
\multiput(1,0)(0,1){3}{\line(1,0){2}}
\multiput(0.1,0)(0,1){3}{
\multiput(0,0)(0.2,0){5}{\line(1,0){0.1}}}
\multiput(3,0)(0,1){3}{
\multiput(0,0)(0.2,0){5}{\line(1,0){0.1}}}
\put(1,0){\makebox(1,1){$\xi$}}
\put(2,1){\makebox(1,1){$\eta$}}
\end{picture}
\quad
\mbox{or}
\quad
\setlength{\unitlength}{5mm}
\begin{picture}(4,2.4)(0,0.3)
\multiput(1,0)(1,0){2}{\line(0,1){2}}
\put(3,1){\line(0,1){1}}
\multiput(1,1)(0,1){2}{\line(1,0){2}}
\put(1,0){\line(1,0){1}}
\multiput(0.1,0)(0,1){3}{
\multiput(0,0)(0.2,0){5}{\line(1,0){0.1}}}
\multiput(3,1)(0,1){2}{
\multiput(0,0)(0.2,0){5}{\line(1,0){0.1}}}
\put(1,0){\makebox(1,1){$\xi$}}
\put(2,1){\makebox(1,1){$\eta$}}
\end{picture}
\,
\mbox{.}
\end{equation}
While the insertion of $\alpha$ that led to these configurations,
let $\eta'$ be the letter that was bumped out from the 
left column. 

Claim 1: $\xi \leq \eta$ and $(\xi,\eta) \ne (n,\ol{n}),(\ol{n},n)$.
To see this note that
in the left column,
B1, B3, B6 or B8 has occurred when $\alpha$ was inserted.
We have $\xi \leq \eta'$ (B1) or $\xi < \eta'$ (B3, B6, B8).
In the right column
A0, B0, B2, B4, B5 or B7 has subsequently occurred.
We have $\eta' = \eta$ (A0, B0, B2, B5), or $\eta' < \eta$ (B4, B7).
In any case we have $\xi \leq \eta$ and $(\xi,\eta) 
\ne (n,\ol{n}),(\ol{n},n)$.

Claim 2: In (\ref{eq:crptd}) the following configurations do not exist.
\begin{equation}
\setlength{\unitlength}{5mm}
\begin{picture}(4,2.4)(0,0.3)
\multiput(1,0)(1,0){3}{\line(0,1){2}}
\multiput(1,0)(0,1){3}{\line(1,0){2}}
\multiput(0.1,0)(0,1){3}{
\multiput(0,0)(0.2,0){5}{\line(1,0){0.1}}}
\multiput(3,0)(0,1){3}{
\multiput(0,0)(0.2,0){5}{\line(1,0){0.1}}}
\put(1,0){\makebox(1,1){$\ol{x}$}}
\put(1,1){\makebox(1,1){$x$}}
\put(2,1){\makebox(1,1){$\ol{x}$}}
\end{picture}
\quad
\mbox{,}
\quad
\setlength{\unitlength}{5mm}
\begin{picture}(4,2.4)(0,0.3)
\multiput(1,0)(1,0){2}{\line(0,1){2}}
\put(3,1){\line(0,1){1}}
\multiput(1,1)(0,1){2}{\line(1,0){2}}
\put(1,0){\line(1,0){1}}
\multiput(0.1,0)(0,1){3}{
\multiput(0,0)(0.2,0){5}{\line(1,0){0.1}}}
\multiput(3,1)(0,1){2}{
\multiput(0,0)(0.2,0){5}{\line(1,0){0.1}}}
\put(1,0){\makebox(1,1){$\ol{x}$}}
\put(1,1){\makebox(1,1){$x$}}
\put(2,1){\makebox(1,1){$\ol{x}$}}
\end{picture}
\quad
\mbox{,}
\quad
\begin{picture}(4,2.4)(0,0.3)
\multiput(1,0)(1,0){3}{\line(0,1){2}}
\multiput(1,0)(0,1){3}{\line(1,0){2}}
\multiput(0.1,0)(0,1){3}{
\multiput(0,0)(0.2,0){5}{\line(1,0){0.1}}}
\multiput(3,0)(0,1){3}{
\multiput(0,0)(0.2,0){5}{\line(1,0){0.1}}}
\put(1,0){\makebox(1,1){$x$}}
\put(2,0){\makebox(1,1){$\ol{x}$}}
\put(2,1){\makebox(1,1){$x$}}
\end{picture}
\quad
\mbox{or}
\quad
\begin{picture}(4,2.4)(0,0.3)
\multiput(1,0)(1,0){3}{\line(0,1){2}}
\multiput(1,0)(0,1){3}{\line(1,0){2}}
\multiput(0.1,0)(0,1){3}{
\multiput(0,0)(0.2,0){5}{\line(1,0){0.1}}}
\multiput(3,0)(0,1){3}{
\multiput(0,0)(0.2,0){5}{\line(1,0){0.1}}}
\put(1,0){\makebox(1,1){$\ol{n}$}}
\put(2,0){\makebox(1,1){$n$}}
\put(2,1){\makebox(1,1){$\ol{n}$}}
\end{picture}
\,
\mbox{.}
\label{eq:forbiddenconfs}
\end{equation}
Due to Claim 1,
the first and the second configurations can exist only if
B1 with $\alpha= x, \gamma=\eta'$ happens in the left column 
and $\xi = \eta'= \eta = \overline{x}$.
But $(\alpha,\gamma) = (x,\overline{x})$ is not compatible with B1.
The third (resp.~fourth) configuration can exist only if 
B4 (resp.~B7) happens in the right column and $\xi= \eta'= \eta= x
\mbox{(resp.~$=\ol{n}$)}$ by Claim 1.
But we see from the proof of Claim 1 that B4 (resp.~B7) actually happens only when $\eta' < \eta$.
Claim 2 is proved.

Let the situation around the crossing 
point of $R$ be one of (\ref{eq:crptd})
excluding (\ref{eq:forbiddenconfs}).
When inserting $\alpha'$, 
suppose in the left column of the crossing point,  
B1, B3, B6 or B8 has occurred.
Let $\xi'$ be the letter bumped out therefrom.

Claim 3: $\xi' \leq \eta$ and $(\xi',\eta) \ne (n,\ol{n}),(\ol{n},n)$.
We divide the check into three cases.
a) If B1 or B6 has occurred in the left column, we have $\xi' = \xi$.
Thus the assertion follows from Claim 1.
b) If B3 has occurred, the left column had the entry ${x \atop \ol{x}}$ and 
we have $\xi' = \ol{x-1}$, $\xi = \ol{x}$.
Claim 1 tells $\xi = \ol{x} \le \eta$, and Claim 2 does $\eta \neq \ol{x}$.
Therefore we have $\xi' = \ol{x-1} \le \eta$.
c) If B8 has occurred we have $\xi' = \ol{n-1}$ for 
either $\xi = n$ or $\xi = \ol{n}$.
If $\xi = n \mbox{ (resp. $\ol{n}$)}$ the entry on the left of $\eta$
was $\ol{n} \mbox{ (resp. $n$)}$, 
therefore $\eta \geq \ol{n} \mbox{ (resp. $n$)}$.
On the other hand Claim 1 tells $(\xi,\eta) \ne (n,\ol{n}),(\ol{n},n)$.
Thus we have $\eta \geq \ol{n-1}$.
Claim 3 is proved.

Now we are ready to finish the proof of the main assertion.
Assume the same situation as Claim 3.
We should verify that A1, B1, B3, B6 and B8 do not 
occur in the right column.
Claim 3 immediately prohibits A1, B1, B6 and B8 in the right column.
Suppose that B3 happens in the right column.
It means that 
$\eta \in \{1,\ldots n\}$, $\xi' \geq \eta$ and the 
right column had the entry
${\eta \atop \ol{\eta}}$.
Since $\xi' \leq \eta$ by Claim 3, we find $\xi' = \eta$, therefore
$\xi' \in \{1,\ldots, n\}$.
Such $\xi'$ can be bumped out from  B1 process only in the left column 
and not from B3, B6 or B8.
It follows that $\xi' = \xi$.
This leads to the third configuration in (\ref{eq:forbiddenconfs}), 
hence a contradiction.

Finally we consider the case where the bumping route $R$
lies only in the second row.
If $R'$ lies below $R$ the tableau should have 
more than two rows, which is prohibited by
Proposition \ref{pr:atmosttworows}.
\end{proof}
\begin{coro}
\label{cor:cblxxx}
Let $\alpha' \leq \alpha$, in particular $(\alpha,\alpha') \ne
(n,\ol{n}),(\ol{n},n)$.
Suppose that a new box is added at the end of the first row
when $\alpha$ is inserted into $T$.
Then a new box is added also at the end of the first row
when $\alpha'$ is inserted into $\left( \alpha \longrightarrow T \right)$.
\end{coro}

\subsection{\mathversion{bold}Main theorem : 
$B^{(1)}_n$ and $D^{(1)}_n$ cases}
\label{subsec:isoruletypeb}
Given $b_1 \otimes b_2 \in B_{l} \otimes B_{k}$,
we define the element
$b'_2 \otimes b'_1 \in B_{k} \otimes B_{l}$
and $l',k', m \in \Z_{\ge 0}$ by the following rule.

\begin{rules}\label{rule:typeB}
\hfill\par\noindent
Set $z = \min(\sharp\,\fbx{1} \text{ in }{\mathcal T}(b_1),\,
\sharp\,\fbx{\ol{1}} \text{ in }{\mathcal T}(b_2))$.
Thus ${\mathcal T}(b_1)$ and ${\mathcal T}(b_2)$ can be depicted by
\[
{\mathcal T}(b_1) =
\setlength{\unitlength}{5mm}
\begin{picture}(7.5,1.4)(0,0.3)
\multiput(0,0)(0,1){2}{\line(1,0){7}}
\put(0,0){\line(0,1){1}}
\put(3,0){\line(0,1){1}}
\put(7,0){\line(0,1){1}}
\put(3,0){\makebox(4,1){$T_*$}}
\put(0,0){\makebox(3,1){$1\cdots 1$}}
\put(0,0.9){\makebox(3,1){$z$}}
\end{picture},\;\;
{\mathcal T}(b_2) =
\setlength{\unitlength}{5mm}
\begin{picture}(6.5,1.4)(0,0.3)
\multiput(0,0)(0,1){2}{\line(1,0){6}}
\put(0,0){\line(0,1){1}}
\put(0.9,0){\line(0,1){1}}
\put(2,0){\line(0,1){1}}
\put(3,0){\line(0,1){1}}
\put(6,0){\line(0,1){1}}
\put(0,0){\makebox(1,1){$v_{1}$}}
\put(1,0){\makebox(1,1){$\cdots$}}
\put(2,0){\makebox(1,1){$v_{k'}$}}
\put(3,0){\makebox(3,1){$\ol{1}\cdots\ol{1}$}}
\put(3,0.9){\makebox(3,1){$z$}}
\end{picture}.
\]
%
Let $l' = l-z$ and $k' = k-z$,
hence $T_*$ is a one-row tableau with length $l'$.
Operate the column insertions and define
\begin{equation}
\label{eq:prodtab}
T^{(0)} := (v_1 \longrightarrow ( \cdots ( v_{k'-1} \longrightarrow ( 
v_{k'} \longrightarrow T_* ) ) \cdots ) ).
\end{equation}
It has the form (See Proposition \ref{pr:atmosttworows}.):
\begin{equation}
\label{eq:prodtabpic}
\setlength{\unitlength}{5mm}
\begin{picture}(20,4)
\put(5,1.5){\makebox(3,1){$T^{(0)}=$}}
\put(8,1){\line(1,0){3.5}}
\put(8,2){\line(1,0){9}}
\put(8,3){\line(1,0){9}}
\put(8,1){\line(0,1){2}}
\put(11.5,1){\line(0,1){1}} 
\put(12.5,2){\line(0,1){1}} 
\put(17,2){\line(0,1){1}}
\put(12.5,2){\makebox(4.5,1){$i_{m+1} \;\cdots\; i_{l'}$}}
\put(8,1){\makebox(3,1){$\;\;i_1 \cdots i_m$}}
\put(8.5,2){\makebox(3,1){$\;\;j_1 \cdots\cdots j_{k'}$}}
\end{picture}
\end{equation}
where $m$ is  the length of the second row, hence that of the first
row is $l'+k'-m$. ($0 \le m \le k'$.)

Next we bump out  $l'$ letters from
the tableau $T^{(0)}$ by the reverse bumping
algorithm.
For the boxes containing $i_{l'}, i_{l'-1}, \ldots, i_1$ in the above
tableau, we do it first for $i_{l'}$ then $i_{l'-1}$ and so on.
Correspondingly, let $w_{1}$ be the first letter that is bumped out from
the leftmost column and $w_2$ be the second  and so on.
Denote by $T^{(i)}$  the resulting tableau when $w_i$ is bumped out
($1 \le i \le l'$).
Now $b'_1 \in B_l$ and $b'_2 \in B_k$ are uniquely specified by
\[
{\mathcal T}(b'_2) =
\setlength{\unitlength}{5mm}
\begin{picture}(6.5,1.4)(0,0.3)
\multiput(0,0)(0,1){2}{\line(1,0){6}}
\put(0,0){\line(0,1){1}}
\put(3,0){\line(0,1){1}}
\put(6,0){\line(0,1){1}}
\put(3,0){\makebox(3,1){$T^{(l')}$}}
\put(0,0){\makebox(3,1){$1\cdots 1$}}
\put(0,0.9){\makebox(3,1){$z$}}
\end{picture},\;\;
{\mathcal T}(b'_1) =
\setlength{\unitlength}{5mm}
\begin{picture}(7.5,1.4)(0,0.3)
\multiput(0,0)(0,1){2}{\line(1,0){7}}
\put(0,0){\line(0,1){1}}
\put(1.25,0){\line(0,1){1}}
\put(2.75,0){\line(0,1){1}}
\put(4,0){\line(0,1){1}}
\put(7,0){\line(0,1){1}}
\put(0,0){\makebox(1.25,1){$w_{1}$}}
\put(1.25,0){\makebox(1.5,1){$\cdots$}}
\put(2.75,0){\makebox(1.25,1){$w_{l'}$}}
\put(4,0){\makebox(3,1){$\ol{1}\cdots\ol{1}$}}
\put(4,0.9){\makebox(3,1){$z$}}
\end{picture}.
\]
\end{rules}
\hfill
(End of the Rule)
\vskip3ex
We normalize the energy function as $H_{B_l B_k}(b_1 \otimes b_2)=0$
for 
\begin{math}
\mathcal{T}(b_1) =
\setlength{\unitlength}{5mm}
\begin{picture}(3,1.5)(0,0.3)
\multiput(0,0)(0,1){2}{\line(1,0){3}}
\multiput(0,0)(3,0){2}{\line(0,1){1}}
\put(0,0){\makebox(3,1){$1\cdots 1$}}
\put(0,1){\makebox(3,0.5){$\scriptstyle l$}}
\end{picture}
\end{math}
and
\begin{math}
\mathcal{T}(b_2) =
\setlength{\unitlength}{5mm}
\begin{picture}(3,1.5)(0,0.3)
\multiput(0,0)(0,1){2}{\line(1,0){3}}
\multiput(0,0)(3,0){2}{\line(0,1){1}}
\put(0,0){\makebox(3,1){$\ol{1}\cdots \ol{1}$}}
\put(0,1){\makebox(3,0.5){$\scriptstyle k$}}
\end{picture}
\end{math}
irrespective of $l < k$ or $l \ge k$.
Our main result for $U'_q(B^{(1)}_n)$ and $U'_q(D^{(1)}_n)$ is the following.
\begin{theorem}\label{th:main3}
Given $b_1 \ot b_2 \in B_l \ot B_k$, find $b'_2 \ot b'_1 \in B_k \ot B_l$
and $l', k', m$ by Rule \ref{rule:typeB} 
with type $B$ (resp.~type $D$) insertion.
Let $\iota: B_l \ot B_k \stackrel{\sim}{\rightarrow} B_k \ot B_l$ be the isomorphism of
$U'_q(B^{(1)}_n)$ (resp.~ $U'_q(D^{(1)}_n)$) crystal.
Then we have
\begin{align*}
\iota(b_1\otimes b_2)& = b'_2 \otimes b'_1,\\
H_{B_l B_k}(b_1 \otimes b_2) &= 2\min(l',k')- m.
\end{align*}
\end{theorem}
%
\noindent
Before giving a proof of this theorem
we present two propositions associated with Rule \ref{rule:typeB}.

Let the product tableau ${\mathcal T}(b_1) *{\mathcal T}(b_2)$
be given by the $T^{(0)}$ 
in eq.~(\ref{eq:prodtab}) in Rule \ref{rule:typeB}.
We assume that it is indeed a (semistandard $B$ or $D$) tableau.
\begin{proposition}
\label{pr:atmosttworows}
The product tableau ${\mathcal T}(b_1) *{\mathcal T}(b_2)$ 
made by (\ref{eq:prodtab}) has no more than two rows.
\end{proposition}
\begin{proof}
Let $T_{\bullet}$ be a tableau that appears in 
the intermediate process of the sequence of the
column insertions (\ref{eq:prodtab}).
Assume that $T_{\bullet}$ has two rows.
We denote by $\alpha$ the letter which we are going to insert 
into $T_{\bullet}$ in the next step of the sequence, and denote by $\beta$
the letter which resides in the second row 
of the leftmost column of $T_{\bullet}$.
It suffices to show that the $\alpha$ does not
make a new box in the third row in the leftmost column.
In other words it suffices to show that $\alpha \leq \beta$ and
$(\alpha,\beta) \ne (\circ,\circ)$ 
(resp.~and in particular $(\alpha,\beta) \ne (n,\ol{n}),(\ol{n},n)$)
in $B_n$ (resp.~$D_n$) case.

Let us first consider $B_n$ case.
We divide the proof in two cases: (i) 
$\beta = \circ$ (ii) $\beta \ne \circ$.
In case (i) either this $\beta=\circ$ was
originally contained in ${\mathcal T}(b_2)$ or 
this $\beta=\circ$ was made by Case B7 in section \ref{subsubsec:insb}.
In any case we see $\alpha \leq n$ (thus $\alpha < \beta$) because of the 
original arrangement of the letters
in ${\mathcal T}(b_2)$.
(Note that ${\mathcal T}(b_2)$ did not have more than one $\circ$s.)
In case (ii) either this $\beta$ was
originally contained in ${\mathcal T}(b_2)$ or 
this $\beta$ is an $\ol{x+1}$ which had
originally been an $\ol{x}$ in ${\mathcal T}(b_2)$ and then
transformed into $\ol{x+1}$ by Case B6 in section \ref{subsubsec:insb}.
In any case we see $\alpha \leq \beta$
and $(\alpha,\beta) \ne (\circ,\circ)$.

Second we consider $D_n$ case.
We divide the proof in two cases: (i) $\beta = n, \ol{n}$ (ii) 
$\beta \ne n, \ol{n}$.
In case (i) either this $\beta=n, \ol{n}$ was
originally contained in ${\mathcal T}(b_2)$ or 
this $\beta$ was made by Case B7 in section \ref{subsubsec:insd}.
In any case we see $\alpha \leq \beta$,
in particular $(\alpha,\beta) \ne (n,\ol{n}),(\ol{n},n)$,
because of the original arrangement of the letters
in ${\mathcal T}(b_2)$.
(Note that ${\mathcal T}(b_2)$ 
did not contain $n$ and $\ol{n}$ simultaneously.)
In case (ii) either this $\beta$ was
originally contained in ${\mathcal T}(b_2)$ or 
this $\beta$ is an $\ol{x+1}$ which had 
originally been an $\ol{x}$ in ${\mathcal T}(b_2)$ and then
transformed into $\ol{x+1}$ by Case B4 in section \ref{subsubsec:insd}.
In any case we see $\alpha \leq \beta$
and $(\alpha,\beta) \ne (n,\ol{n}),(\ol{n},n)$.
\end{proof}

Let $\geh = B^{(1)}_n \mbox{ or } D^{(1)}_n$
and $\ol{\geh} = B_n \mbox{ or } D_n$.
By neglecting zero arrows, the crystal graph of 
$B_l \ot B_k$
decomposes into $U_q(\ol{\geh})$ crystals,
where only arrows with indices $i=1,\ldots,n$ remain.
Let us regard $b_1 \in B_l$ as an element of 
$U_q(\ol{\geh})$ crystal $B(l \Lambda_1)$, and
regard $b_2 \in B_k$ as an element of $B(k \Lambda_1)$.
Then $b_1 \ot b_2$ is regarded as an element of $U_q(\ol{\geh})$ crystal 
$B(l \Lambda_1) \ot B(k \Lambda_1)$.
On the other hand the tableau
${\mathcal T}(b_1) *{\mathcal T}(b_2)$ 
specifies an element of $B(\lambda)$
which we shall denote by $b_1 * b_2$,
where $B(\lambda)$ is a $U_q(\ol{\geh})$ crystal
that appears in the decomposition of $B(l \Lambda_1) \ot B(k \Lambda_1)$.
\begin{proposition}
\label{pr:crysmorpb}
The map $\psi : b_1 \ot b_2 \mapsto b_1 * b_2$ is a
$U_q(\ol{\geh})$ crystal morphism, i.e. the 
actions of $\etd_i$ and $\ftd_i$
for $i=1,\ldots,n$ 
commute with the map $\psi$.
\end{proposition}
\noindent
This proposition is a special case of Proposition~\ref{pr:morphgen}.
Note that, although
we have removed the $z$ pairs of $1$'s and $\ol{1}$'s from the tableaux
by hand, this elimination of the letters
is also a part of this rule
of column insertions (i.e. 
\begin{math}
(\ol{1} \longrightarrow 
\setlength{\unitlength}{5mm}
\begin{picture}(1,1)(0,0.3)
\multiput(0,0)(1,0){2}{\line(0,1){1}}
\multiput(0,0)(0,1){2}{\line(1,0){1}}
\put(0,0){\makebox(1,1){$1$}}
\end{picture}
) = \emptyset
\end{math}),
followed by the sliding (jeu de taquin) rules \cite{B1,B2}.
%

\par\noindent
\begin{proof}[Proof of Theorem \ref{th:main3}]
First we consider the isomorphism.
We are to show:
\begin{enumerate}
	\item If $b_1 \ot b_2$ is mapped to $b'_2 \ot b'_1$ under
	the isomorphism, then the product tableau ${\mathcal T}(b_1) *
	{\mathcal T}(b_2)$ is equal to the product tableau
	${\mathcal T}(b'_2) * {\mathcal T}(b'_1)$.
	\item If $k$ and $l$ are specified, we can recover 
	${\mathcal T}(b'_2) $ and ${\mathcal T}(b'_1)$ from their
	product tableau by using the algorithm shown in Rule \ref{rule:typeB}.
	In other words, we can retrieve them by assuming
	the arrangement of the locations $i_{l'},\ldots,i_1$ of 
	the boxes in (\ref{eq:prodtabpic}) from 
	which we start the reverse bumpings.
\end{enumerate}
Claim 2 is verified from
Corollary \ref{cor:bcblxxx} or \ref{cor:cblxxx}.
We consider Claim 1 in the following.
The value of the energy function value will be settled at the same time.

Thanks to the $U_q(\ol{\geh})$ crystal 
morphism (Proposition \ref{pr:crysmorpb}),
it suffices to prove the theorem
for any element in each connected component of the
$U_q(\ol{\geh})$ crystal.
We take up the $U_q(\ol{\geh})$ highest weight element 
as such a particular element.
There is a special extreme $U_q(\ol{\geh})$ highest weight element
\begin{equation}
\iota : (l,0,\ldots,0) \ot (k,0,\ldots,0) 
 \stackrel{\sim}{\mapsto} 
(k,0,\ldots,0) \ot (l,0,\ldots,0),
\label{eqn:ultrahighest}
\end{equation}
wherein we find that they are obviously mapped to each other
under the $U'_q(\geh)$ isomorphism,
and that the image of the map is also obviously
obtained by Rule \ref{rule:typeB}.
Let us assume  $l \geq k$.
(The other case can be treated in a similar way.)
Suppose that $b_1 \ot b_2 \in B_l \ot B_k$
is a $U_q(\ol{\geh})$  highest element.
In general, it has the form:
$$b_1\ot b_2 = (l,0,\ldots,0) \ot (x_1,x_2,0,\ldots,0,\ol{x}_1),$$
%
where $x_1, x_2$ and $\ol{x}_1$ are arbitrary
as long as $k = x_1 + x_2 + \ol{x}_1$.
We are to obtain its image under the isomorphism.
Applying
\begin{eqnarray*}
&&\et{0}^{\ol{x}_1} \et{2}^{x_2+\ol{x}_1} \cdots
\et{n-1}^{x_2+\ol{x}_1} \et{n}^{2x_2+2\ol{x}_1}
\et{n-1}^{x_2+\ol{x}_1}
\cdots \et{2}^{x_2+\ol{x}_1} \et{0}^{x_2+\ol{x}_1}
\; \mbox{for $\geh = B^{(1)}_n$}\\
&&\et{0}^{\ol{x}_1} \et{2}^{x_2+\ol{x}_1} \cdots
\et{n-1}^{x_2+\ol{x}_1} \et{n}^{x_2+\ol{x}_1}
\et{n-2}^{x_2+\ol{x}_1}
\cdots \et{2}^{x_2+\ol{x}_1} \et{0}^{x_2+\ol{x}_1}
\; \mbox{for $\geh = D^{(1)}_n$}
\end{eqnarray*}
to the both sides of (\ref{eqn:ultrahighest}), 
we find
\begin{displaymath}
\iota :
(l,0,\ldots,0) \ot (x_1,x_2,0,\ldots,0,\ol{x}_1)
\stackrel{\sim}{\mapsto} 
(k,0,\ldots,0) \ot (x_1',x_2,0,\ldots,0,\ol{x}_1).
\end{displaymath}
Here  $x_1' = l-x_2-\ol{x}_1$.
In the course of the application of $\tilde{e}_i$'s,
the value of the energy function has changed as
$$
H\left((l,0,\ldots,0) \ot (x_1,x_2,0,\ldots,0,\ol{x}_1)\right) =
H\left((l,0,\ldots,0) \ot (k,0,\ldots,0)\right) - x_2 - 2 \ol{x}_1.
$$
(We have omitted the subscripts of the energy function.)
Thus according to our normalization we have
$H(b_1 \ot b_2)=2(k- \ol{x}_1)-x_2 $.
(Note that the $z$ in Rule \ref{rule:typeB} 
is now equal to $\ol{x}_1$, hence
we have $k' = k-\ol{x}_1$.)
On the other hand
for this highest element
the column insertions lead to a common tableau
\begin{displaymath}
\setlength{\unitlength}{5mm}
\begin{picture}(7,2)
\put(0,0.5){\makebox(3,1){$T^{(0)}=$}}
\put(3,2){\line(1,0){4}}
\put(3,1){\line(1,0){4}}
\put(3,0){\line(1,0){3}}
\put(3,0){\line(0,1){2}}
\put(6,0){\line(0,1){1}}
\put(7,1){\line(0,1){1}}
\put(3,1){\makebox(4,1){$1\cdots\cdots 1$}}
\put(3,0){\makebox(3,1){$2\cdots 2$}}
\end{picture}
\end{displaymath}
whose second row has length $x_2$
(and first row has the length  $l+k-x_2-2\ol{x}_1$).
This completes the proof.
\end{proof}
\subsection{Examples}
\label{subsec:exBD}
\begin{example}
$B_5 \otimes B_3 \simeq B_3 \otimes B_5$ for $B^{(1)}_5$.
\begin{displaymath}
\begin{array}{ccccccc}
\setlength{\unitlength}{5mm}
\begin{picture}(5,1)(0,0.3)
\multiput(0,0)(1,0){6}{\line(0,1){1}}
\multiput(0,0)(0,1){2}{\line(1,0){5}}
\put(0,0){\makebox(1,1){$5$}}
\put(1,0){\makebox(1,1){$5$}}
\put(2,0){\makebox(1,1){$\circ$}}
\put(3,0){\makebox(1,1){$\ol{5}$}}
\put(4,0){\makebox(1,1){$\ol{5}$}}
\end{picture}
& \otimes & 
\setlength{\unitlength}{5mm}
\begin{picture}(3,1)(0,0.3)
\multiput(0,0)(1,0){4}{\line(0,1){1}}
\multiput(0,0)(0,1){2}{\line(1,0){3}}
\put(0,0){\makebox(1,1){$5$}}
\put(1,0){\makebox(1,1){$\circ$}}
\put(2,0){\makebox(1,1){$\ol{5}$}}
\end{picture}
& \stackrel{\sim}{\mapsto}  &
\setlength{\unitlength}{5mm}
\begin{picture}(3,1)(0,0.3)
\multiput(0,0)(1,0){4}{\line(0,1){1}}
\multiput(0,0)(0,1){2}{\line(1,0){3}}
\put(0,0){\makebox(1,1){$5$}}
\put(1,0){\makebox(1,1){$\circ$}}
\put(2,0){\makebox(1,1){$\ol{5}$}}
\end{picture}
& \otimes & 
\setlength{\unitlength}{5mm}
\begin{picture}(5,1)(0,0.3)
\multiput(0,0)(1,0){6}{\line(0,1){1}}
\multiput(0,0)(0,1){2}{\line(1,0){5}}
\put(0,0){\makebox(1,1){$5$}}
\put(1,0){\makebox(1,1){$5$}}
\put(2,0){\makebox(1,1){$\circ$}}
\put(3,0){\makebox(1,1){$\ol{5}$}}
\put(4,0){\makebox(1,1){$\ol{5}$}}
\end{picture}
\\
& & & & & & \\
\setlength{\unitlength}{5mm}
\begin{picture}(5,1)(0,0.3)
\multiput(0,0)(1,0){6}{\line(0,1){1}}
\multiput(0,0)(0,1){2}{\line(1,0){5}}
\put(0,0){\makebox(1,1){$5$}}
\put(1,0){\makebox(1,1){$5$}}
\put(2,0){\makebox(1,1){$\ol{5}$}}
\put(3,0){\makebox(1,1){$\ol{5}$}}
\put(4,0){\makebox(1,1){$\ol{5}$}}
\end{picture}
& \otimes & 
\setlength{\unitlength}{5mm}
\begin{picture}(3,1)(0,0.3)
\multiput(0,0)(1,0){4}{\line(0,1){1}}
\multiput(0,0)(0,1){2}{\line(1,0){3}}
\put(0,0){\makebox(1,1){$4$}}
\put(1,0){\makebox(1,1){$4$}}
\put(2,0){\makebox(1,1){$\circ$}}
\end{picture}
& \stackrel{\sim}{\mapsto}  &
\setlength{\unitlength}{5mm}
\begin{picture}(3,1)(0,0.3)
\multiput(0,0)(1,0){4}{\line(0,1){1}}
\multiput(0,0)(0,1){2}{\line(1,0){3}}
\put(0,0){\makebox(1,1){$\circ$}}
\put(1,0){\makebox(1,1){$\ol{5}$}}
\put(2,0){\makebox(1,1){$\ol{5}$}}
\end{picture}
& \otimes & 
\setlength{\unitlength}{5mm}
\begin{picture}(5,1)(0,0.3)
\multiput(0,0)(1,0){6}{\line(0,1){1}}
\multiput(0,0)(0,1){2}{\line(1,0){5}}
\put(0,0){\makebox(1,1){$4$}}
\put(1,0){\makebox(1,1){$4$}}
\put(2,0){\makebox(1,1){$5$}}
\put(3,0){\makebox(1,1){$5$}}
\put(4,0){\makebox(1,1){$\ol{5}$}}
\end{picture}
\\
& & & & & & \\
\setlength{\unitlength}{5mm}
\begin{picture}(5,1)(0,0.3)
\multiput(0,0)(1,0){6}{\line(0,1){1}}
\multiput(0,0)(0,1){2}{\line(1,0){5}}
\put(0,0){\makebox(1,1){$1$}}
\put(1,0){\makebox(1,1){$1$}}
\put(2,0){\makebox(1,1){$\circ$}}
\put(3,0){\makebox(1,1){$\ol{5}$}}
\put(4,0){\makebox(1,1){$\ol{5}$}}
\end{picture}
& \otimes & 
\setlength{\unitlength}{5mm}
\begin{picture}(3,1)(0,0.3)
\multiput(0,0)(1,0){4}{\line(0,1){1}}
\multiput(0,0)(0,1){2}{\line(1,0){3}}
\put(0,0){\makebox(1,1){$\circ$}}
\put(1,0){\makebox(1,1){$\ol{1}$}}
\put(2,0){\makebox(1,1){$\ol{1}$}}
\end{picture}
& \stackrel{\sim}{\mapsto}  &
\setlength{\unitlength}{5mm}
\begin{picture}(3,1)(0,0.3)
\multiput(0,0)(1,0){4}{\line(0,1){1}}
\multiput(0,0)(0,1){2}{\line(1,0){3}}
\put(0,0){\makebox(1,1){$1$}}
\put(1,0){\makebox(1,1){$1$}}
\put(2,0){\makebox(1,1){$\circ$}}
\end{picture}
& \otimes & 
\setlength{\unitlength}{5mm}
\begin{picture}(5,1)(0,0.3)
\multiput(0,0)(1,0){6}{\line(0,1){1}}
\multiput(0,0)(0,1){2}{\line(1,0){5}}
\put(0,0){\makebox(1,1){$\circ$}}
\put(1,0){\makebox(1,1){$\ol{5}$}}
\put(2,0){\makebox(1,1){$\ol{5}$}}
\put(3,0){\makebox(1,1){$\ol{1}$}}
\put(4,0){\makebox(1,1){$\ol{1}$}}
\end{picture}
\end{array}
\end{displaymath}
Here we have picked up three samples that are specific to type $B$.
The values of the energy function are 
assigned to be 3, 5 and 1, respectively.

Let us illustrate in more detail the procedure of Rule \ref{rule:typeB}
by taking the first example.
{}From the left hand side we proceed the column insertions as follows.
\begin{align*}
\ol{5} &\rightarrow 
\setlength{\unitlength}{5mm}
\begin{picture}(5,1)(0,0.3)
\multiput(0,0)(1,0){6}{\line(0,1){1}}
\multiput(0,0)(0,1){2}{\line(1,0){5}}
\put(0,0){\makebox(1,1){$5$}}
\put(1,0){\makebox(1,1){$5$}}
\put(2,0){\makebox(1,1){$\circ$}}
\put(3,0){\makebox(1,1){$\ol{5}$}}
\put(4,0){\makebox(1,1){$\ol{5}$}}
\end{picture}
\quad = \quad
\setlength{\unitlength}{5mm}
\begin{picture}(5,2)(0,0.8)
\multiput(0,0)(1,0){2}{\line(0,1){2}}
\multiput(2,1)(1,0){4}{\line(0,1){1}}
\multiput(0,1)(0,1){2}{\line(1,0){5}}
\put(0,0){\line(1,0){1}}
\put(0,1){\makebox(1,1){$5$}}
\put(1,1){\makebox(1,1){$5$}}
\put(2,1){\makebox(1,1){$\circ$}}
\put(3,1){\makebox(1,1){$\ol{5}$}}
\put(4,1){\makebox(1,1){$\ol{5}$}}
\put(0,0){\makebox(1,1){$\ol{5}$}}
\end{picture}
\\
\circ &\rightarrow 
\setlength{\unitlength}{5mm}
\begin{picture}(5,2)(0,0.8)
\multiput(0,0)(1,0){2}{\line(0,1){2}}
\multiput(2,1)(1,0){4}{\line(0,1){1}}
\multiput(0,1)(0,1){2}{\line(1,0){5}}
\put(0,0){\line(1,0){1}}
\put(0,1){\makebox(1,1){$5$}}
\put(1,1){\makebox(1,1){$5$}}
\put(2,1){\makebox(1,1){$\circ$}}
\put(3,1){\makebox(1,1){$\ol{5}$}}
\put(4,1){\makebox(1,1){$\ol{5}$}}
\put(0,0){\makebox(1,1){$\ol{5}$}}
\end{picture}
\quad = \quad
\setlength{\unitlength}{5mm}
\begin{picture}(5,2)(0,0.8)
\multiput(0,0)(1,0){3}{\line(0,1){2}}
\multiput(3,1)(1,0){3}{\line(0,1){1}}
\multiput(0,1)(0,1){2}{\line(1,0){5}}
\put(0,0){\line(1,0){2}}
\put(0,1){\makebox(1,1){$4$}}
\put(1,1){\makebox(1,1){$5$}}
\put(2,1){\makebox(1,1){$\circ$}}
\put(3,1){\makebox(1,1){$\ol{5}$}}
\put(4,1){\makebox(1,1){$\ol{5}$}}
\put(0,0){\makebox(1,1){$\circ$}}
\put(1,0){\makebox(1,1){$\ol{4}$}}
\end{picture}
\\
5 &\rightarrow 
\setlength{\unitlength}{5mm}
\begin{picture}(5,2)(0,0.8)
\multiput(0,0)(1,0){3}{\line(0,1){2}}
\multiput(3,1)(1,0){3}{\line(0,1){1}}
\multiput(0,1)(0,1){2}{\line(1,0){5}}
\put(0,0){\line(1,0){2}}
\put(0,1){\makebox(1,1){$4$}}
\put(1,1){\makebox(1,1){$5$}}
\put(2,1){\makebox(1,1){$\circ$}}
\put(3,1){\makebox(1,1){$\ol{5}$}}
\put(4,1){\makebox(1,1){$\ol{5}$}}
\put(0,0){\makebox(1,1){$\circ$}}
\put(1,0){\makebox(1,1){$\ol{4}$}}
\end{picture}
\quad = \quad
\setlength{\unitlength}{5mm}
\begin{picture}(5,2)(0,0.8)
\multiput(0,0)(1,0){4}{\line(0,1){2}}
\multiput(4,1)(1,0){2}{\line(0,1){1}}
\multiput(0,1)(0,1){2}{\line(1,0){5}}
\put(0,0){\line(1,0){3}}
\put(0,1){\makebox(1,1){$4$}}
\put(1,1){\makebox(1,1){$5$}}
\put(2,1){\makebox(1,1){$\circ$}}
\put(3,1){\makebox(1,1){$\ol{5}$}}
\put(4,1){\makebox(1,1){$\ol{5}$}}
\put(0,0){\makebox(1,1){$5$}}
\put(1,0){\makebox(1,1){$\circ$}}
\put(2,0){\makebox(1,1){$\ol{4}$}}
\end{picture}
\end{align*}
\vskip3ex
\noindent
The reverse bumping procedure goes as follows.
\begin{align*}
T^{(0)} &=
\setlength{\unitlength}{5mm}
\begin{picture}(5,2)(0,0.8)
\multiput(0,0)(1,0){4}{\line(0,1){2}}
\multiput(4,1)(1,0){2}{\line(0,1){1}}
\multiput(0,1)(0,1){2}{\line(1,0){5}}
\put(0,0){\line(1,0){3}}
\put(0,1){\makebox(1,1){$4$}}
\put(1,1){\makebox(1,1){$5$}}
\put(2,1){\makebox(1,1){$\circ$}}
\put(3,1){\makebox(1,1){$\ol{5}$}}
\put(4,1){\makebox(1,1){$\ol{5}$}}
\put(0,0){\makebox(1,1){$5$}}
\put(1,0){\makebox(1,1){$\circ$}}
\put(2,0){\makebox(1,1){$\ol{4}$}}
\end{picture}
& \\
T^{(1)} &=
\setlength{\unitlength}{5mm}
\begin{picture}(5,2)(0,0.8)
\multiput(0,0)(1,0){4}{\line(0,1){2}}
\put(4,1){\line(0,1){1}}
\multiput(0,1)(0,1){2}{\line(1,0){4}}
\put(0,0){\line(1,0){3}}
\put(0,1){\makebox(1,1){$4$}}
\put(1,1){\makebox(1,1){$\circ$}}
\put(2,1){\makebox(1,1){$\ol{5}$}}
\put(3,1){\makebox(1,1){$\ol{5}$}}
\put(0,0){\makebox(1,1){$5$}}
\put(1,0){\makebox(1,1){$\circ$}}
\put(2,0){\makebox(1,1){$\ol{4}$}}
\end{picture}
&, w_1 = 5 \\
T^{(2)} &=
\setlength{\unitlength}{5mm}
\begin{picture}(5,2)(0,0.8)
\multiput(0,0)(1,0){4}{\line(0,1){2}}
\multiput(0,0)(0,1){3}{\line(1,0){3}}
\put(0,1){\makebox(1,1){$4$}}
\put(1,1){\makebox(1,1){$\circ$}}
\put(2,1){\makebox(1,1){$\ol{5}$}}
\put(0,0){\makebox(1,1){$\circ$}}
\put(1,0){\makebox(1,1){$\ol{5}$}}
\put(2,0){\makebox(1,1){$\ol{4}$}}
\end{picture}
&, w_2 = 5 \\
T^{(3)} &=
\setlength{\unitlength}{5mm}
\begin{picture}(5,2)(0,0.8)
\multiput(0,0)(1,0){3}{\line(0,1){2}}
\put(3,1){\line(0,1){1}}
\multiput(0,1)(0,1){2}{\line(1,0){3}}
\put(0,0){\line(1,0){2}}
\put(0,1){\makebox(1,1){$4$}}
\put(1,1){\makebox(1,1){$\circ$}}
\put(2,1){\makebox(1,1){$\ol{5}$}}
\put(0,0){\makebox(1,1){$\ol{5}$}}
\put(1,0){\makebox(1,1){$\ol{4}$}}
\end{picture}
&, w_3 = \circ \\
T^{(4)} &=
\setlength{\unitlength}{5mm}
\begin{picture}(5,2)(0,0.8)
\multiput(0,0)(1,0){2}{\line(0,1){2}}
\multiput(2,1)(1,0){2}{\line(0,1){1}}
\multiput(0,1)(0,1){2}{\line(1,0){3}}
\put(0,0){\line(1,0){1}}
\put(0,1){\makebox(1,1){$5$}}
\put(1,1){\makebox(1,1){$\circ$}}
\put(2,1){\makebox(1,1){$\ol{5}$}}
\put(0,0){\makebox(1,1){$\ol{5}$}}
\end{picture}
&, w_4 = \ol{5} \\
T^{(5)} &=
\setlength{\unitlength}{5mm}
\begin{picture}(5,2)(0,0.3)
\multiput(0,0)(1,0){4}{\line(0,1){1}}
\multiput(0,0)(0,1){2}{\line(1,0){3}}
\put(0,0){\makebox(1,1){$5$}}
\put(1,0){\makebox(1,1){$\circ$}}
\put(2,0){\makebox(1,1){$\ol{5}$}}
\end{picture}
&, w_5 = \ol{5}
\end{align*}
Thus we obtained the right hand side.
We have $H_{B_5,B_3}=3$, since $l'=5, k'=3$ and $m=3$.
\end{example}
\begin{example}
$B_2 \otimes B_1 \simeq B_1 \otimes B_2$ for $D^{(1)}_5$.
\begin{displaymath}
\begin{array}{ccccccc}
\setlength{\unitlength}{5mm}
\begin{picture}(2,1)(0,0.3)
\multiput(0,0)(1,0){3}{\line(0,1){1}}
\multiput(0,0)(0,1){2}{\line(1,0){2}}
\put(0,0){\makebox(1,1){$4$}}
\put(1,0){\makebox(1,1){$\ol{4}$}}
\end{picture}
& \otimes & 
\setlength{\unitlength}{5mm}
\begin{picture}(1,1)(0,0.3)
\multiput(0,0)(1,0){2}{\line(0,1){1}}
\multiput(0,0)(0,1){2}{\line(1,0){1}}
\put(0,0){\makebox(1,1){$5$}}
\end{picture}
& \stackrel{\sim}{\mapsto}  &
\setlength{\unitlength}{5mm}
\begin{picture}(1,1)(0,0.3)
\multiput(0,0)(1,0){2}{\line(0,1){1}}
\multiput(0,0)(0,1){2}{\line(1,0){1}}
\put(0,0){\makebox(1,1){$\ol{5}$}}
\end{picture}
& \otimes & 
\setlength{\unitlength}{5mm}
\begin{picture}(2,1)(0,0.3)
\multiput(0,0)(1,0){3}{\line(0,1){1}}
\multiput(0,0)(0,1){2}{\line(1,0){2}}
\put(0,0){\makebox(1,1){$5$}}
\put(1,0){\makebox(1,1){$5$}}
\end{picture}
\\
& & & & & & \\
\setlength{\unitlength}{5mm}
\begin{picture}(2,1)(0,0.3)
\multiput(0,0)(1,0){3}{\line(0,1){1}}
\multiput(0,0)(0,1){2}{\line(1,0){2}}
\put(0,0){\makebox(1,1){$\ol{5}$}}
\put(1,0){\makebox(1,1){$\ol{5}$}}
\end{picture}
& \otimes & 
\setlength{\unitlength}{5mm}
\begin{picture}(1,1)(0,0.3)
\multiput(0,0)(1,0){2}{\line(0,1){1}}
\multiput(0,0)(0,1){2}{\line(1,0){1}}
\put(0,0){\makebox(1,1){$5$}}
\end{picture}
& \stackrel{\sim}{\mapsto}  &
\setlength{\unitlength}{5mm}
\begin{picture}(1,1)(0,0.3)
\multiput(0,0)(1,0){2}{\line(0,1){1}}
\multiput(0,0)(0,1){2}{\line(1,0){1}}
\put(0,0){\makebox(1,1){$\ol{5}$}}
\end{picture}
& \otimes & 
\setlength{\unitlength}{5mm}
\begin{picture}(2,1)(0,0.3)
\multiput(0,0)(1,0){3}{\line(0,1){1}}
\multiput(0,0)(0,1){2}{\line(1,0){2}}
\put(0,0){\makebox(1,1){$4$}}
\put(1,0){\makebox(1,1){$\ol{4}$}}
\end{picture}
\end{array}
\end{displaymath}
Here we have picked up two samples that are specific to type $D$.
\end{example}
\vspace{0.4cm}
\noindent
{\bf Acknowledgements} \hspace{0.1cm}
It is a pleasure to thank T.H. Baker for many helpful discussions
and correspondence.

\end{document}